\documentclass[a4paper]{amsart}
\usepackage{verbatim}
\usepackage[dvips]{color}
\usepackage{amssymb}

\author{Giuseppe Mingione}
\address{Dipartimento di Matematica, Universit\`a di Parma\\
Viale G.~P.~Usberti 53/a, Campus, 43100 Parma, Italy; e-mail:
giuseppe.mingione@unipr.it.}

\allowdisplaybreaks

\newtheorem{theorem}{Theorem}[section]
\newtheorem{prop}{Proposition}[section]
\newtheorem{lemma}{Lemma}[section]
\newtheorem{cor}{Corollary}[section]
\theoremstyle{definition}
\newtheorem{definition}{Definition}
\newtheorem{remark}{Remark}[section]

\numberwithin{equation}{section}

\newcommand{\Dim}{\textnormal{dim}_{\mathcal{H}}}
\newcommand{\CC}{\subset\hskip-.2em\subset }
\newcommand{\rif}[1]{(\ref{#1})}

\newcommand{\mmm}[1]{(s^2+|z|^2)^{\frac{p-#1}{2}}}

\newcommand{\BMO}{\textnormal{BMO}}
\newcommand{\VMO}{\textnormal{VMO}}
\newcommand{\tu}{\tilde{u}}

\def\eqn#1$$#2$${\begin{equation}\label#1#2\end{equation}}
\def\charfn_#1{{\raise1.2pt\hbox{$\chi
_{\kern-1pt\lower3pt\hbox{{$\scriptstyle#1$}}}$}}}

\def\qq1{q_*}
\def\q2{q_{**}}
\def\dist{\operatorname{dist}}

\def\ep{\varepsilon}
\def\en{\mathbb N}
\def\er{\mathbb R}

\def\loc{\operatorname{loc}}

\def\ma{\mathbb R^{nN}}

\newdimen\vintbar
\def\vint{-\kern-\vintbar\int}

\def\B{\mathcal B}

\def\D{\mathcal D}

\def\F{\mathcal F}

\def\N{\mathcal N}

\def\0{\boldsymbol 0}

\newcommand{\aint}{-\!\!\!\!\!\! \int}

\newcommand{\Wp}{W^{1,p}}

\newcommand{\ratio}{L/\nu}

\newcommand{\MM}{\mathcal{M}}
\newcommand{\Qin}{Q_{\textnormal{inn}}}
\newcommand{\Qou}{Q_{\textnormal{out}}}
\newcommand{\dro}{(\varrho-t)}

\newcommand{\divo}{\textnormal{div}}

 \newcommand{\mean}[1]{-\hskip-1.08em\int_{#1}}

\newcommand{\omp}{\Omega^{\prime}}
\newcommand{\ompp}{\Omega^{\prime \prime}}
\newcommand{\distc}{\textnormal{dist}(\omp, \partial \ompp)}
\newcommand{\M}{M^*}
\newcommand{\Ma}{\er^{N\times n}}

\newcommand{\cu}{{\rm B}}

\newcommand{\trif}[1] {\textnormal{\rif{#1}}}

\newcommand{\VVd}{|V(Du)-V(Dv)|^2}
\newcommand{\VVq}{|V(Du)-V(Dv)|^{2q/p}}
\newcommand{\QQ}{2B_0}

\newtoks\by
\newtoks\paper
\newtoks\book
\newtoks\jour
\newtoks\yr
\newtoks\pages
\newtoks\vol
\newtoks\publ
\def\et{ \& }
\def\name[#1, #2]{#1 #2}
\def\ota{{\hbox{\bf ???}}}
\def\cLear{\by=\ota\paper=\ota\book=\ota\jour=\ota\yr=\ota
\pages=\ota\vol=\ota\publ=\ota}
\def\endpaper{\the\by, \textit{\the\paper},
{\the\jour} \textbf{\the\vol} (\the\yr), \the\pages.\cLear}
\def\endbook{\the\by, \textit{\the\book},
\the\publ, \the\yr.\cLear}
\def\endpap{\the\by, \textit{\the\paper}, \the\jour.\cLear}
\def\endproc{\the\by, \textit{\the\paper}, \the\book, \the\publ,
\the\yr, \the\pages.\cLear}

\title[Calder\'on-Zygmund estimates and measure
data] {The Calder\'on-Zygmund theory for elliptic problems with
measure data}

\begin{document}
Ann.~SNS Pisa, Cl.~Sci.~V, Vol.~6 (2007), 195-261

 \setcounter{tocdepth}{1}

\begin{abstract}
We consider non-linear elliptic equations having a measure in the
right hand side, of the type $ \divo \ a(x,Du)=\mu, $
 and prove differentiability and integrability
results for solutions. New estimates in Marcinkiewicz spaces are
also given, and the impact of the measure datum density properties
on the regularity of solutions is analyzed in order to build a
suitable Calder\'on-Zygmund theory for the problem. All the
regularity results presented in this paper are provided together
with explicit local a priori estimates.
\end{abstract}
 \maketitle
\centerline{{\em To the memory of Vic Mizel, mathematician and
gentleman}}
 \tableofcontents
\section{Introduction and results}
Let us consider the following Dirichlet problem: \eqn{Dir1}
$$
\left\{
    \begin{array}{cc}
    -\divo \ a(x,Du)=\mu & \qquad \mbox{in $\Omega$}\\
        u= 0&\qquad \mbox{on $\partial\Omega$.}
\end{array}\right.
$$
Here we assume that $\Omega \subset \er^n$ is a bounded domain,
$\mu$ is a signed Radon measure with finite total variation
$|\mu|(\Omega)< \infty$, and $a\colon \Omega \times \er^n\to \er^n$
is a Carath\`eodory vector field satisfying the following standard
monotonicity and Lipschitz assumptions: \eqn{asp}
$$
\left\{
    \begin{array}{c}
   \nu (s^2+|z_1|^2+|z_2|^2)^{\frac{p-2}{2}}|z_2-z_1|^2 \leq
\langle
a(x,z_2)-a(x,z_1),z_2-z_1\rangle  \\ \\
    | a(x,z_2)-a(x,z_1)| \leq L
(s^2+|z_1|^2+|z_2|^2)^{\frac{p-2}{2}}|z_2-z_1|\\ \\
|a(x,0)|\leq Ls^{p-1} \;,
    \end{array}
    \right.
$$
for every $z_1,z_2\in \er^n$, $x \in \Omega$. Here, and in the rest
of the paper, when referring to the structural
 properties of $a$, and in particular to \rif{asp}, we shall always assume
 \eqn{datadata}
$$p\geq 2, \qquad n\geq 2, \qquad 0< \nu \leq L, \qquad s \geq 0\;.$$
The measure $\mu$ will be considered as defined on the whole $\er^n$
by simply letting $|\mu|(\er^n\setminus \Omega)=0$. At certain
stages, we shall also require the following Lipschitz continuity
assumption on the map $x \mapsto a(x,z)$: \eqn{lipi}
$$
|a(x,z)-a(x_{0},z)| \leq
    L|x-x_{0}|\mmm{1},\qquad \forall \ x,x_0 \in \Omega, \ z \in
    \er^n\;.
$$
Assumptions \rif{asp} are modeled on the basic example \eqn{plap}
$$
-\divo [c(x)(s^2+|Du|^2)^{\frac{p-2}{2}}Du]= \mu, \qquad  \nu \leq
c(x)\leq L\;,
$$
which is indeed covered here. When $s=0$ and $c(x)\equiv 1$ we have
the familiar $p$-Laplacean operator on the left-hand side\eqn{pmeas}
$$-\triangle_p u = -\divo (|Du|^{p-2}Du)=\mu\;.$$

For the problem \rif{Dir1} in the rest of the paper we shall adopt
the following distributional-like notion of solution, compare with
\cite{BG1} for instance.
\begin{definition}\label{solsol}
A solution $u$ to the problem \rif{Dir1} under assumptions
\rif{asp}, is a function $u \in W^{1,1}_0(\Omega)$ such that
$a(x,Du) \in L^{1}(\Omega,\er^n)$ and \eqn{desol}
$$
\int_{\Omega} a(x,Du)D\varphi \, dx =\int_{\Omega} \varphi\, d \mu,
\qquad \textnormal{for every}\  \varphi \in C^{\infty}_c(\Omega)\;.
$$
\end{definition}
The existence of such a solution is usually obtained combining a
priori estimates with a suitable approximation scheme \cite{BG1,
 DHM2, Dallaglio}, see also Section 5 below. The same approach is followed
 here and therefore in the
rest of the paper when talking about regularity we shall refer to
that of Solutions Obtained as Limits of Approximations (SOLA)
\cite{boccardo, Dallaglio}, and {\em we shall actually
simultaneously obtain existence and regularity results}. Here we
just want to recall that uniqueness of solutions to \rif{Dir1} in
the sense of Definition \ref{solsol} generally fails \cite{Serrin},
and a main open problem of the theory is identifying a suitable
functional class where {\em a unique solution} can be defined and
found. In this respect many possible definitions have been proposed,
and technically demanding attempts have been made: for this we refer
for instance to \cite{elenco, BOG, KiXu, Ra2}, and to the references
therein. Nevertheless, a general uniqueness theory is still missing
except for $p=2$ or $p=n$ \cite{boccardo, DHM2, GIS}; in particular
we refer to the paper \cite{DMOP} for a rather comperhensive
discussion about the uniqueness problem, and measure data problems
in general. We shall not discuss uniqueness problems any further,
our aims here being quite different: we are mainly interested in
{\em a priori regularity estimates}. For the same reason, we shall
confine ourselves to distributional solutions as defined in
\rif{solsol}, while the results we are going to propose could be
approached also for other notions of solutions: entropy ones, for
instance.

The study of problem \rif{Dir1} began with the fundamental work of
Littman\et Stampacchia\et Weinberger \cite{LSW, stamp}, who defined
solutions in a duality sense in the case of linear equations with
measurable coefficients: $a_i(x,z)\equiv \tilde{a}_{ij}(x)z_j$. When
referring to Definition \ref{solsol}, the existence theory for the
general quasi linear Leray-Lions type operators in \rif{Dir1}$_1$
has been established in the by now classical paper of Boccardo\et
Gallou\"et \cite{BG1}, who proved the existence of a solution $u$ to
problem \rif{Dir1} such that \eqn{boccardo}
$$
Du \in L^{q}(\Omega,\er^n), \qquad \textnormal{for every}\  q <b \
\mbox{when}  \ p\leq n\;,
$$
where \eqn{boccardoex}
$$
b: = \frac{n(p-1)}{n-1}\;.
$$
 Dolzmann\et Hungerb\"uhler\et M\"uller were able to prove the same
 result for a
large class of systems including the $p$-Laplacean one \cite{DHM1,
DHM2}. Inclusion \rif{boccardo} is optimal in the scale of Lebesgue
spaces, see Section 11.1, as $Du \not \in L^b$ in general. Anyway
\rif{boccardo} can be sharpened using Marcinkiewicz spaces
\cite{elenco, DHM2}, see \rif{dema} below, since \eqn{Marbase}
$$Du \in \MM^b(\Omega, \er^n)\;.$$
When $p>n$ instead, $\mu$ belongs to $W^{-1,p'}$, that is the dual
of $W^{1,p}$, and the existence of a unique solution in the natural
space $W^{1,p}_0(\Omega)$ follows by standard duality methods
\cite{lions}. Related regularity results for the equation
\rif{pmeas} with a non-negative measure $\mu$ were given by
Lindqvist \cite{Lind}, in connection to the notion of
$p$-superharmonic functions; see also \cite{HKM} for a fairly
comprehensive treatment of this subject. Related estimates and
problems, using various techniques, are in \cite{delvecchio, FR,HKM,
KiL,Ta}.

\subsection {General measures}Up to now, regularity results in $L^q$
spaces of the type in \rif{boccardo}-\rif{Marbase} are the only ones
available in the literature. {\em One of the aims of this paper is
to give the first higher regularity results for the gradient of
solutions}, in particular estimating the oscillations of the
gradient rather than its size. Let us focus for simplicity on the
case $p=2$, looking at \rif{boccardo} from a different viewpoint,
considering $\triangle u =f$. In this case {\em the standard
Calder\'on-Zygmund theory} \cite{GT} asserts \eqn{SCZ}
$$f \in L^{1+\ep} \Longrightarrow Du  \in W^{1,1+\ep} \qquad \qquad \textnormal{for every}\  \ep >0\;.$$
Using Sobolev's embedding theorem we have in particular $Du \in
L^{n/(n-1)}$, that is, the limit case of \rif{boccardo}. This does
not hold when $\ep = 0$, since the inclusion $Du \in W^{1,1}$
generally fails. So, one could interpret \rif{boccardo} {\em as the
trace of a potentially existent Calder\'on-Zygmund theory below the
limit case $W^{1,1}$.} Indeed we have
\begin{theorem}[of Calder\'on-Zygmund type]\label{main1}
Under the assumptions \trif{asp} and \trif{lipi} with $p\leq n$,
there exists a solution $u \in W^{1,1}_0(\Omega)$ to the problem
\trif{Dir1} such that \eqn{pipi}
$$
Du \in W^{\frac{1-\ep}{p-1},p-1}_{\loc}(\Omega,\er^n) \;,
$$
for every $\ep \in (0,1)$, and in particular \eqn{duedue}
 $$Du \in W^{1-\ep,1}_{\loc}(\Omega,\er^n)\,, \qquad \qquad \mbox{when}\quad p=2\;.$$
More in general \eqn{fra1}
$$
Du \in W^{\frac{\sigma(q)-\ep}{q},q}_{\loc}(\Omega,\er^n)\;,
$$
for every $\ep \in (0,\sigma(q))$, where \eqn{ssii}
$$
 p-1 \leq q <
   \frac{n(p-1)}{n-1}=b,  \hspace{1cm}  \ \sigma(q):=
    n- \frac{q(n-1)}{p-1}=n(1-q/b) \;,
$$and $b$ is in
\trif{boccardoex}.
\end{theorem}
In other words, in \rif{duedue} we ``almost have" second derivatives
of $u$; see anyway \rif{secondeu} below and comments after
\rif{plala}. Explicit local estimates are actually available:
\begin{theorem}[Calder\'on-Zygmund estimates I]\label{main1es}
Under the assumptions and notations of Theorem \ref{main1}, let $q
\in [p-1,b)$ and $\sigma \in (0,\sigma(q))$. There exists a constant
$c \equiv c(n,p,\ratio, q,\sigma)$ such that for every ball $B_R\CC
\Omega$ of radius $R>0$ it holds \begin{eqnarray} \int_{B_{R/2}}
\int_{B_{R/2}} \! \frac{|Du(x) - Du(y) |^{q}}{|x-y|^{n+\sigma}} \ dx
\ dy  & \leq & \frac{c}{R^\sigma}\int_{B_R} (|Du|^q+s^q)\,
dx\nonumber \\ && \qquad +
cR^{\sigma(q)-\sigma}[|\mu|(\overline{B_R})]
^{\frac{q}{p-1}}\;.\label{stimaloc1}
\end{eqnarray}
Moreover, for any open subset $\Omega'\CC \Omega$ the local estimate
\eqn{apest}
$$
\int_{\Omega} \! |Du|^{q}\, dx + \int_{\Omega'} \int_{\Omega'} \!
\frac{|Du(x) - Du(y) |^{q}}{|x-y|^{n+\sigma}} \ dx \ dy \leq
c[|\mu|(\Omega)] ^{\frac{q}{p-1}}+cs^q|\Omega|
$$
holds, where $c \equiv c(n,p,\ratio,q,\sigma, \dist(\Omega',
\partial \Omega),\Omega)$.
\end{theorem}
Therefore, {\em it is possible to establish an optimal
Calder\'on-Zygmund theory for non-linear elliptic problems with
measure data, provided the right Sobolev spaces are considered:}
that is, the fractional ones; see the definition in
\rif{defrazionari}. Fractional Sobolev spaces are an essential tool
in modern analysis in that they provide the natural intermediate
scale to state optimal regularity results, and to show the
persistence of certain assertions up to the so called ``limit
cases". Inclusion \rif{duedue} is an instance of this situation, and
the comparison between \rif{duedue} and \rif{SCZ} tells us that
Calder\'on-Zygmund theory does not have $W^{1,1}$ as an end-point,
but it continues below $W^{1,1}$. Inclusions \rif{pipi}-\rif{fra1}
are sharp for every choice of the couple $(q,\sigma(q))$ in
\rif{ssii} as $Du \not\in W^{\sigma(q)/q,q}_{\loc}$ in general; note
that $(p-1,\sigma(p-1))=(p-1,1)$. On the converse, inclusion
\rif{fra1} admits \rif{boccardo} as a corollary, at least in a local
fashion; see Section 11.1.~When $p\not=2$, as $\ep \searrow 0$ we do
not approach an integer fractional differentiability exponent in
\rif{pipi}, as for \rif{duedue}, but only $1/(p-1)$. This is not a
surprise: even for the model case \eqn{plala}
$$\triangle_p u=0\;,$$ the existence of second derivatives of
$W^{1,p}$-solutions is not clear due to the degeneracy of the
problem, while fractional derivatives naturally appear \cite{Min1}
\eqn{frazioni}
$$ u \in W^{1,p}\Longrightarrow Du \in W^{2/p, p}\;. $$
On the other hand, a classical result going back to K.~Uhlenbeck
\cite{U} asserts that although $Du$ may be not differentiable for
\rif{plala}, {\em certain natural non-linear expressions of the
gradient still are} (in T.~Iwaniec's words \cite{donofrio}). Indeed,
defining
$$
V(Du):=(s^2+|Du|^2)^{\frac{p-2}{4}}Du\;,
$$
then under assumptions \rif{asp} we have that $V(Du) \in
W^{1,2}_{\loc}(\Omega,\mathbb R^n)$ for any $W^{1,p}$-solution to
\eqn{omogenea}
$$ \divo \ a(Du)=0\;;
$$ see Lemma \ref{revbase} below, and under stronger
assumptions also \cite{Ham, Manp, Manth}, \cite{G}, Chapter 8, and
\cite{Lindp}, Chapter 4. See also Section 11.2 for more comments on
the fact that passing to $V(Du)$ allows for a gain in
differentiability. Observe that the main essence here is that the
product between the differentiability and the integrability indexes
of the fractional spaces involved for $Du$ and $V(Du)$,
respectively, is invariant \eqn{propro}
$$
 \frac{2}{p}\cdot p =1\cdot 2\;.
$$
This phenomenon extends to measure data problems as well:
\begin{theorem}[Non-linear Calder\'on-Zygmund estimates]\label{noncz}
Under the assumptions \trif{asp} and \trif{lipi} with $p \in (2,n]$,
let $u \in W^{1,q}_0(\Omega)$ be the solution to \trif{Dir1} found
in Theorem \ref{main1}. Then \eqn{nonfra1}
$$V(Du) \in
W^{\frac{p}{2(p-1)}-\ep,\frac{2(p-1)}{p}}_{\loc}(\Omega,\er^n),
\qquad \mbox{for every}\ \ep \in (0,1)\;.$$ Moreover, for any open
subset $\Omega'\CC \Omega$, we have \eqn{nonapest}
$$
 \int_{\Omega'} \int_{\Omega'} \!
\frac{|V(Du(x)) - V(Du(y)) |^{\frac{2(p-1)}{p}}}{|x-y|^{n+1-\ep}} \
dx \ dy \leq \tilde{c} |\mu|(\Omega)+\tilde{c}s^{p-1}|\Omega|\;,
$$
where $\tilde{c} \equiv \tilde{c} (n,p,\ratio,\ep, \dist(\Omega',
\partial \Omega),\Omega)$.
\end{theorem}
When the problem is non-degenerate, that is $s>0$ in \rif{datadata},
something more can be proved: $W^{1,p}$-solutions to \rif{omogenea}
belong to $W^{2,2}_{\loc}$. The following corollary of Theorem
\ref{noncz} contains the analogue in the measure data case.
\begin{cor}[The non-degenerate case]\label{nondegenere}
Under the assumptions \trif{asp} and \trif{lipi} with $p \in (2,n]$,
let $u \in W^{1,q}_0(\Omega)$ be the solution to \trif{Dir1} found
in Theorem \ref{main1}, and assume $s>0$. Then \eqn{secondeu}
$$Du \in
W^{\frac{p}{2(p-1)}-\ep,\frac{2(p-1)}{p}}_{\loc}(\Omega,\er^n),
\qquad \mbox{for every}\ \ep \in (0,1)\;.$$ Moreover estimate
\trif{nonapest} holds with $Du$ replacing $V(Du)$ provided the
constant $\tilde{c}$ is replaced by the new one:
$s^{(2-p)(p-1)/p}c(n,p)\tilde{c} $.
\end{cor}
\begin{remark} In Theorems \ref{main1es} and \ref{noncz}, and Corollary
\ref{nondegenere}, the constants $c, \tilde{c}$ depending on
$q,\ep,\sigma$ blow-up as $q\nearrow b$, $\ep \searrow 0$, $\sigma
\nearrow \sigma(q)$. Also observe that formally letting $p=2$ in the
three previous statements we obtain \rif{duedue}.
\end{remark} Combining inclusions \rif{duedue} and \rif{nonfra1} with Proposition \ref{redu} below we gain
\begin{cor}[BV-type behavior]\label{bv}
Let $\Sigma_u$ denote the set of non-Lebesgue points of the solution
found in Theorem \ref{main1}, in the sense of \eqn{singularset}
$$
\begin{array}{c}
 \displaystyle \Sigma_u:= \left\{ x \in \Omega :
\liminf_{\rho\searrow 0}
 \mean{B(x,\rho)} |Du(y)-(Du)_{x,\rho}| \ dy
 >0\right. \\ \displaystyle\left. \hspace{60mm} \mbox{or}\qquad \limsup_{\rho\searrow 0}
 |(Du)_{x,\rho}| =\infty \right\}\;.
\end{array}
$$
Then its Hausdorff dimension $\dim(\Sigma_u)$ satisfies \eqn{stimah}
$$
\dim(\Sigma_u)\leq n-1\;.
$$
The same result holds replacing $Du$ by $V(Du)$ in
\trif{singularset}.
\end{cor}
Therefore solutions behave as BV functions \cite{AFP}. For $p=2$ one
can guess this, with some brave heuristics, by looking at $\triangle
u =\mu$, ``replacing" $\triangle u$ by $D^2u$.

Before going on let us observe that the above results are only
local, while we are dealing with a Dirichlet problem; this is a
precise, simplifying choice of ours. Indeed the techniques presented
here are suitable to be carried out up to the boundary under
additional regularity assumptions on $\partial \Omega$, say $C^2$
for instance, or Lipschitz in some cases, but since they are already
delicate and involved, at this stage we prefer to confine ourselves
to the local versions of the results, in order to highlight the main
new ideas. For the case $p< 2$ see also \cite{minmis}; here the
results change.

\subsection{\bf Diffusive measures} The sharpness of \rif{Marbase}
and \rif{pipi}-\rif{fra1} stems from considering counterexamples
involving Dirac measures, see again Section 11.1. It is therefore
natural to wonder whether things change when considering measures
diffusing on sets with higher Hausdorff dimension. A natural way to
quantify this, also suggested by a classical result of Frostman, see
\cite{AdHe}, Theorem 5.1.12, is to consider the following density
condition: \eqn{morrey}
$$
|\mu|(B_R) \leq M R^{n-\theta} \qquad \qquad 0 \leq \theta \leq n,
\qquad M\geq 0,
$$
satisfied for any ball $B_R\subset \er^n$ of radius $R$. Assuming
\rif{morrey} does not allow $\mu$ to concentrate on sets with
Hausdorff dimension less than $n-\theta$, and indeed higher
regularity of solutions can be obtained. We have anyway to
distinguish two cases.

\subsection{\bf The super-capacitary case} This is when $\theta \geq p$,
making sense only if $p\leq n$. We see that in all the above results
the role of the dimension $n$ is actually played by $\theta$ in
\rif{morrey}; in particular, the critical exponent $b$ appearing in
Theorem \ref{main1} and in \rif{boccardo} is replaced by the larger
one \eqn{mingione}
$$m:=\frac{\theta(p-1)}{\theta-1}\;.$$
\indent The first improvement is in the integrability properties of
$Du$, detectable in two different scales: Marcinkiewicz and Morrey
ones, see \rif{mama} and \rif{inclusionimorrey}.
\begin{theorem} [Marcinkiewicz-Morrey regularity] \label{main5}
Under the assumptions \trif{asp} with $p\leq n$, and \trif{morrey}
with $\theta \geq p$, there exists a solution $u \in
W^{1,1}_0(\Omega)$ to the problem \trif{Dir1} such that \eqn{lolo}
$$
Du \in \MM^{m, \theta}_{\loc}(\Omega,\er^n)\subseteq
\MM^{m}_{\loc}(\Omega,\er^n)\;,
$$
where $m$ is in \trif{mingione}. Moreover, for any open subset
$\Omega'\CC \Omega$ we have \eqn{apestlo}
$$
\|Du\|_{\MM^{m, \theta}(\Omega')} \leq cM^{\frac{1}{p-1}}
+cs|\Omega|^{\frac{1}{m}}\;,
$$
where $c \equiv c(n,p,\ratio, \Omega',\Omega)$, and $M$ appears in
\trif{morrey}. In particular, in the limit case $\theta=p$ we have
\eqn{lolo2}
$$
Du \in \MM^{p,p}_{\loc}(\Omega,\er^n)\subseteq
\MM^{p}_{\loc}(\Omega,\er^n)\;.
$$
\end{theorem}
The exponent $m$ is expected to be the best possible in \rif{lolo}
for every $p\geq 2$, and it actually is when $p=2$, see Section
11.3: Theorem \ref{main5} may be also regarded as the non-linear
version of a classical result of Adams \cite{adamsduke}. As
explained below, when $\theta<p$, the solution $u$ is uniquely found
in $W^{1,p}_0(\Omega)$, so that \trif{morrey} provides the natural
scale that allows to pass from \rif{Marbase}, when $\theta=n$, to
\rif{lolo2}, when $\theta=p$; in this last case the
$W^{1,p}$-regularity of the solution is missed just by a natural
Marcinkiewicz-scale factor. A warning for the reader: {\em in
$\MM^{q,\theta}$ the second exponent does not ``tune" the first
one}; these are not like Lorentz spaces: indeed $\MM^{q,0}\equiv
L^{\infty}$ and $\MM^{q,n}\equiv \MM^q$, for every $q\geq 1$.
Finally, note that \rif{lolo} does not require \rif{lipi} since we
are not dealing with higher derivatives of the gradient, and we do
not need to ``differentiate" equation \rif{Dir1}$_1$, that obviously
needs \rif{lipi}.

The second effect of condition \rif{morrey} is an expansion of the
range \rif{ssii}. The fractional derivatives are themselves in a
Morrey space, see the definition in \rif{defrac2}. This leads to the
final and central stage of our regularity program:
\begin{theorem} [Sobolev-Morrey regularity] \label{main3}
Under the assumptions \trif{asp} and \trif{lipi} with $p\leq n$, and
\trif{morrey} with $\theta \geq p$, let $u \in W^{1,1}_0(\Omega)$ be
the solution found in Theorem \ref{main5}. Then\eqn{fra13mo}
$$
Du \in W^{\frac{\sigma(q,\theta)-\ep}{q},q,
\theta}_{\loc}(\Omega,\er^n)\;,
$$
for every $\ep \in (0,\sigma(q,\theta))$, where $m$ is in
\trif{mingione}, and \eqn{ssii2}
$$
p-1 \leq q <
    \frac{\theta(p-1)}{\theta-1}=m,  \hspace{6mm}  \ \sigma(q,\theta):=
    \theta- \frac{q(\theta-1)}{p-1}=\theta(1- q/m) \;.
$$
In particular \eqn{dueduemo}
 $$Du \in W^{\frac{1-\ep}{p-1},p-1,\theta}_{\loc}(\Omega,\er^n), \quad \mbox{and} \quad
  Du \in W^{1-\ep,1,\theta}_{\loc}(\Omega,\er^n) \quad  \mbox{when}\quad p=2\;.$$
Moreover, for any open subset $\Omega'\CC \Omega$ and $\sigma \in
(0,\sigma(q,\theta))$, we have \eqn{apestmo}
$$
\|Du\|_{W^{\sigma/q,q, \theta}(\Omega')} \leq
cM^{\frac{1}{p-1}}+cs|\Omega|^{\frac{1}{q}}\;,
$$
where $c \equiv c(n,p,\ratio,q,\sigma, \dist(\Omega', \partial
\Omega),\Omega)$, and $M$ is in \trif{morrey}.
\end{theorem}
Originally introduced in \cite{Ca63, Ca64}, Sobolev-Morrey spaces
$W^{\alpha,q,\theta}$ appear in various forms in several pde issues
as they provide the natural scaling properties of solutions
\cite{Kozono, Mazzucato, Taylor}. Estimate \rif{apestmo} extends to
the case of non-linear equations with measure data the classical
Morrey space results for linear elliptic equations \cite{Ca70, G,
Caffmi, Difazio93, Difazio99, Lieb03}; see the definition in
\rif{dede1mo} below. The standard result for the model case
$\triangle v =f$ is that $Dv \in W^{1,q,\theta}$ when $f \in
L^{q,\theta}$ for $q>1$, that is
$$\int_{B_R} |D^2v|^q\, dx \leq c R^{n-\theta}\;.$$
Inclusion \rif{dueduemo} sharply extends this to the case $q=1$,
that is previous inequality is replaced by the following analogue:
$$
\int_{B_{R}} \int_{B_{R}} \! \frac{|Du(x) - Du(y)
|}{|x-y|^{n+1-\ep}} \ dx \ dy \leq c R^{n-\theta}\;,
$$
which is valid for every $\ep \in (0,1)$ and every ball $B_R \CC
\Omega$ of radius $R$, where $c$ depends on $\ep$ and on the
distance between $B_R$ and $\partial \Omega$. Finally, in light of
\rif{lolo} we can interpret \rif{fra13mo} and therefore also
\rif{fra1} as a scale of regularity for $Du$ leading, as $q \nearrow
m$, from the maximal differentiability \rif{dueduemo} towards the
maximal integrability \rif{lolo}.

\subsection{\bf The capacitary case} This is
when $\theta <p $; this case is simpler and we will be shorter. Here
$\mu \in W^{-1, p'}$, that is, the dual space of $W^{1,p}$, and
moreover $\mu$ cannot charge null $p$-capacity sets. When $p>n$ this
follows from Sobolev's embedding theorem; one-point sets have
positive $p$-capacity. When $\theta<p\leq n$ a basic theorem of
D.~R.~Adams \cite{adams, AdHe} still ensures that $\mu \in W^{-1,
p'}$; here note that \rif{morrey} implies $|\mu|(B_R)\preceq
R^{p-\theta}$Cap$_p(B_R)$. At the end \rif{Dir1} can be solved by
monotonicity methods \cite{lions}, and the existence of a unique
solution in the natural space $W^{1,p}_0(\Omega)$ follows.
\begin{theorem}[Capacitary Calder\'on-Zygmund estimates]\label{main4}
Assume \trif{asp}, \trif{lipi}, and \trif{morrey} with $\theta < p$.
Then the unique solution $u \in W^{1,p}_0(\Omega)$ to the problem
\trif{Dir1} satisfies \eqn{fra123}
$$
Du \in W^{\frac{\sigma(p)-\ep}{p},p}_{\loc}(\Omega,\er^n),\qquad
\qquad \qquad \sigma(p):= \frac{p-\theta}{p-1}\;,
$$
for every $\ep \in (0,\sigma(p))$. Moreover, for any open subset
$\Omega'\CC \Omega$ we have \eqn{apestsub}
$$
\int_{\Omega} \! |Du|^{p}\, dx + \int_{\Omega'} \int_{\Omega'} \!
\frac{|Du(x) - Du(y) |^{p}}{|x-y|^{n+\sigma}} \ dx \ dy \leq
cM^{\frac{1}{p-1}}|\mu|(\Omega) +cs^p|\Omega|\;,
$$
where $c \equiv c(n,p,\ratio,\sigma, \theta, \dist(\Omega', \partial
\Omega),\Omega)$, and $M$ appears in \trif{morrey}.
\end{theorem}
As a corollary of \rif{fra123} and of the fractional Sobolev
embedding Theorem \ref{fraemb}, we also have the following higher
integrability result: \eqn{highr}
$$
Du \in L^{t}_{\loc}(\Omega,\er^n) \qquad \qquad \mbox{for every} \ t
< \frac{np}{n-\sigma(p)}\;.
$$
We point out the analogy between \rif{fra123} and the results in
\cite{KZ} for the case $\theta < p \leq n$, stating that solutions
to \rif{pmeas} are $C^{0,\alpha}$-regular with $\alpha =\sigma(p)$;
see also \cite{Kilp, KiMaly, liebe}. Theorem \ref{main4} extends to
general elliptic systems, see Section 11.6.

\subsection{\bf Additional results} For the proof of the above theorems we
shall need the following intermediate result, which may have its own
interest; see also \cite{cirmi, DL} for a particular case.
\begin{theorem} [Morrey space regularity] \label{mainmorrey}
Under the assumptions \trif{asp} with $p\leq n$, and \trif{morrey}
with $\theta \geq p$, let $u \in W^{1,1}_0(\Omega)$ be the solution
found in Theorem \ref{main5}. Then $Du$ belongs to the Morrey space
$ L^{q, \delta(q)}_{\loc}(\Omega,\er^n)$, for every $q$ and
$\delta(q)$ such as \eqn{ildelta}
$$1\leq q < \frac{\theta(p-1)}{\theta-1}=m, \qquad \qquad \delta(q):= \frac{q(\theta-1)}{p-1}\;.$$
For every $\Omega' \CC \Omega$, there exists $c\equiv
c(n,p,\ratio,q,\dist(\Omega', \partial \Omega),\Omega)$ such that
\eqn{stimamorrey}
$$
\|Du\|_{L^{q,\delta(q)}(\Omega')} \leq c
M^{\frac{1}{p-1}}+cs|\Omega|^{\frac{1}{q}}\;.
$$
\end{theorem}
 For the sake of
completeness we also include a corollary that in different forms,
but not in the following one, already appears in the literature
\cite{ DHM2, FF, zhong}.
\begin{theorem} [$\BMO/\VMO$ regularity] \label{main6}
Under the assumptions \trif{asp} with $p\leq n$, and \trif{morrey}
with $\theta =p$, the solution $u \in W^{1,1}_0(\Omega)$ found in
Theorem \ref{main5} belongs to  \BMO$_{\loc}(\Omega)$. Moreover, if
\eqn{locuni}
$$
\lim_{R \searrow 0} \frac{|\mu|(B_R)}{R^{n-p}}=0\;,
$$
locally uniformly in $\Omega,$ then $u \in $\VMO$_{\loc}(\Omega)$.
Finally, for every open subset $\Omega' \CC \Omega$ there exists a
constant $c \equiv c(n,p,\ratio,\dist(\Omega',\partial \Omega),
\Omega)$ such that \eqn{apbmo}
$$
[u]_{\BMO(\Omega')}\leq cM^{\frac{1}{p-1}}+ cs|\Omega|\;.
$$
\end{theorem}
For the exact meaning of ``locally uniformly" in \rif{locuni} see
Definition \ref{lloocc} below; see also Remark \ref{abbs}. Observe
that also in this case the result complements the ones in the
literature: as soon as $\theta<p$ solutions are H\"older continuous
\cite{Kilp, liebe}.

\begin{remark} In Theorems \ref{main3}-\ref{main6} the constants $c$ depending on
$q,\ep$ blow-up as $q\nearrow m$, $\ep \searrow 0$, $\sigma\nearrow
\sigma(q,\theta)$; $\sigma \nearrow\sigma(p)$ in case of Theorem
\ref{main4}.
\end{remark}

Finally, a road-map to the paper. Some of the results presented are
obtained via a delicate interaction between various types of
regularity scales. For instance, as for Theorems \ref{main5} and
Theorem \ref{main3}, we have \eqn{catena}
$$ Du \in L^{q,
\delta(q)}_{\loc}(\Omega,\er^n), \ q<b \Longrightarrow Du \in
\MM^{m, \theta}_{\loc}\Longrightarrow Du \in
W^{\frac{\sigma(q,\theta)-\ep}{q},q, \theta}_{\loc}\;.$$ In Section
2 we collect a miscellanea of preliminary material and notations.
Section 3 includes some results for elliptic problems, that in the
form presented are not explicitly contained in the literature. In
Section 4 we collect a few preparatory lemmas of comparison type,
while in Section 5 we fix the basic approximation procedure. Section
6 contains the proofs of Theorems \ref{main1}-\ref{noncz} and
Corollary \ref{nondegenere}, while Section 7 contains the one of
Theorem \ref{main4}. Section 8 contains the proof of Theorems
\ref{mainmorrey} and \ref{main6}. Section 9 contains the proof of
Theorem \ref{main5}, while Section 10 contains the one of Theorem
\ref{main3}. Finally, in Section 11 we discuss the sharpness of some
of the results obtained.

Part of the results obtained in this paper have been announced in
the preliminary Comptes Rendus note \cite{mincras}.

{\em Acknowledgements.} Part of this paper was carried out while the
author was visiting FIM at ETH Z\"urich in January 2006, and the
Institut f\"ur Mathematik at Humboldt University Berlin in April
2006, institutions that provided a nice and stimulating environment;
thanks go to Tristan Rivi\`ere and Barbara Niethammer, respectively,
for their invitations. Thanks also go to Nicola Fusco, Jan
Kristensen, Francesco Leonetti, Juan J.~Manfredi for stimulating
conversations on problems with measure data; special thanks finally
go to Anna F\"oglein and Tero Kilpel\"ainen for a careful reading of
an earlier version of the paper.
\section{Preliminaries, function spaces}
\subsection{General notation} In this paper we shall adopt the usual, but somehow arguable
convention to denote by $c$ a general constant, that may vary from
line to line; peculiar dependence on parameters will be properly
emphasized in parentheses when needed, while special occurences will
be denoted by $c_*, c_1, c_2$ or the like. With $x_0\equiv
(x_{0,1},\dots ,x_{0,n}) \in \er^n$, we denote $$B_R(x_0)\equiv
B(x_0,R) :=\{x \in \er^n :  |x-x_0|<R\}\;,$$ and $$Q_R(x_0)\equiv
Q(x_0,R):=\{x \in \er^n  :  \sup |x_i-x_{0,i}|<R, \ 1\leq i \leq
n\}\;,$$ the open ball and cube, respectively, with center $x_0$ and
``radius" $R$. We shall often use the short hand notation $B_R
\equiv B(x_0,R)$ and $Q_R\equiv Q(x_0,R)$, when no ambiguity will
arise. Moreover, with $B,Q$ being balls and cubes, respectively, by
$\gamma B, \gamma Q$ we shall denote the concentric balls and cubes,
with radius magnified by a factor $\gamma.$ If $g\colon A \to \er^k$
is an integrable map with respect to the Borel measure $\mu$, and
$0<\mu(A)<\infty$, we write
$$
    (g)_{\mu,A}  := \aint_{A} \! g(x) \, d\mu :=
    \frac{1}{\mu(A)}\int_{A} \! g(x) \, d\mu \,.
$$
When $\mu$ is the Lebesgue measure and $A \equiv B(x_0,R)$, we may
also use the short hand notation $(g)_{\mu,A} \equiv (g)_A \equiv
(g)_{B_R} \equiv (g)_B$.

{\em Permanent conventions.}  In the estimates the constants will in
general depend on the parameters $n,p,\nu,L$. The dependence on
$\nu,L$ is actually via the {\em ellipticity ratio} $\ratio$, and
will be given directly in this way. This can be seen by passing to
rescaled vector fields $a/\nu$. When considering a function space
$X(\Omega,\er^k)$ of possibly vector valued measurable maps defined
on an open set $\Omega \subset \er^n$, with $k \in \en$, e.g.:
$L^p(\Omega, \er^k), W^{\alpha,p}(\Omega, \er^k)$, we shall define
in a canonic way the local variant $X_{\loc}(\Omega, \er^k)$ as that
space of maps $f:\Omega \to \er^k$ such that $f \in X(\Omega',
\er^k)$, for every $\Omega' \CC \Omega.$ Moreover, also in the case
$f$ is vector valued, that is $k >1$, we shall also use the short
hand notation $X(\Omega, \er^k) \equiv X(\Omega)$, or even
$X(\Omega, \er^k)\equiv X$.

\subsection{\bf The map $V(z)$, and the monotonicity of $a(x,z)$} With $s \geq
0$, we define \eqn{Vfun}
$$
V(z) = V_{s}(z) := (s^2+|z|^2)^{\frac{p-2}{4}}z\;, \qquad \qquad z
\in \er^n\;,
$$
which is easily seen to be a locally bi-Lipschitz bijection of
$\er^n$. A basic property of $V$, whose proof can be found in
\cite{Ham}, Lemma 2.1, is the following: For any $z_1, z_2 \in
\er^n$, and any $s\geq 0$, it holds \eqn{V}
$$
c^{-1}\Bigl( s^2+|z_1|^2+|z_2|^2 \Bigr)^{\frac{p-2}{2}}\leq
\frac{|V(z_2)-V(z_1)|^2}{|z_2-z_1|^2} \leq c\Bigl(
s^2+|z_1|^2+|z_2|^2 \Bigr)^{\frac{p-2}{2}}\;,
$$
where $c\equiv c(n,p)$, is independent of $s$. We also notice that
\eqn{elemV}
$$
|z|^p\leq |V(z)|^2\leq 2(s^p+|z|^p)\;.
$$
Indeed when $p=2$ this is trivial, otherwise when $p>2$ we just use
Young's inequality with conjugate exponents $(p, p/(p-2))$; in what
follows we shall also need another elementary property of $V$:
\eqn{elemV2}
$$
V_{s/A}(z/A)=A^{-p/2}V_{s}(z), \qquad \mbox{for every} \ s \geq 0,\
\mbox{and} \ A>0\;.
$$
The strict monotonicity properties of the vector field $a$ implied
by the left hand side in \rif{asp}$_1$ can be recast using the map
$V$. Indeed combining \rif{asp}$_1$ and \rif{V} yields, for $c\equiv
c(n,p,\nu)>0$, and whenever $z_1,z_2 \in \er^n$ \eqn{mon3}
$$   c^{-1}
|V(z_2)-V(z_1)|^2 \leq  \langle
a(x,z_2)-a(x,z_1),z_2-z_1\rangle\;.$$ Moreover, since $p\geq 2$,
assumption \rif{asp}$_1$ also implies \eqn{mon2}
$$   c^{-1}
|z_2-z_1|^p \leq  \langle a(x,z_2)-a(x,z_1),z_2-z_1\rangle\;.$$
Finally, inequality \rif{asp}$_1$, together with  \rif{asp}$_3$ and
a standard use of Young's inequality, yield for every $z \in
\er^n$\eqn{veramon}
$$
c^{-1}(s^2+|z|^2)^{\frac{p-2}{2}}|z|^2-cs^p \leq \langle
a(x,z),z\rangle, \qquad \qquad c \equiv c(n,p,\ratio)\;,
$$
while \rif{asp}$_2$ and again \rif{asp}$_3$ give, again via Young's
inequality \eqn{veracre}
$$
|a(x,z)|\leq c (s^2+|z|^2)^{\frac{p-1}{2}}\;.
$$
In the following and for the rest of the paper, unless otherwise
stated, we shall denote $V \equiv V_s$ with $s$ fixed at the
beginning, in \rif{datadata}.

\subsection{\bf Fractional Sobolev/Nikolski spaces, and difference operators}
We recall some basic facts about fractional order Sobolev spaces,
using the standard notation from \cite{Ad}, adapted to the
situations we are going to deal with. For a bounded open set $A
\subset \er^n$ and $k \in \en$, parameters $\alpha \in (0,1)$ and $q
\in [1,\infty )$, we write $w \in W^{\alpha ,q}(A,\er^k )$ provided
the following Gagliardo-type norm is finite:
\begin{eqnarray}
\nonumber \| w \|_{W^{\alpha,q}(A )} & := & \left(\int_{A} \!
|w(x)|^{q}\, dx \right)^{\frac{1}{q}} + \left(\int_{A} \int_{A} \!
\frac{|w(x)
- w(y) |^{q}}{|x-y|^{n+\alpha q}} \ dx \, dy \right)^{\frac{1}{q}}\\
&   =:&  \|w \|_{L^q(A)} +[ w ]_{\alpha ,q;A} <
\infty\;.\label{defrazionari}
\end{eqnarray}

For a possibly vector valued function $w \colon \Omega \to \er^k$,
and a real number $h \in \er$, we define the finite difference
operator $\tau_{i,h}$ for $i \in \{1,\ldots, n\}$ as\eqn{detau}
$$
\tau_{i,h}w(x)\equiv \tau_{i,h}(w)(x):=w(x+he_i)-w(x)\;.
$$
where  $\{e_i\}_{1\leq i \leq n}$ denotes the standard basis of
$\er^n$. This makes sense whenever $x, x +he_i \in A$, an assumption
that will be satisfied whenever we use $ \tau_{i,h}$ in the
following. In particular, we shall very often take $x \in A$ where
$A \CC \Omega$ is an open subset of $\Omega$, and where $|h|\leq
\dist(A,\partial \Omega)$. Accordingly, the Nikolski space
$\N^{\alpha,q}(A)$, with $A \CC \Omega$ is hereby defined by saying
that $u \in \N^{\alpha,q}(A)$ if and only if
$$
\| w \|_{\N^{\alpha,q}(A )} :=\| w\|_{L^{q}(A)} +
\sum_{i=1}^{n}\sup_{h}|h|^{-\alpha}\|\tau_{i,h} w\|_{L^{q}(A)} \;,
$$
is finite, where $0< |h|\leq \dist(A,\partial \Omega)$. In the
following we shall also let $W^{0,q}\equiv \N^{0,q}\equiv L^q$. The
strict inclusions $$ W^{\alpha,q}(A)\subset \N^{\alpha,q}(A)\subset
W^{\alpha-\ep,q}(A)\,, \qquad \forall \ \ep \in (0,\alpha)\,,$$ are
well known, and the next lemma somehow quantifies the last one.
\begin{lemma} \label{diff2}
Let $w\in L^q (\Omega ,\er^k)$, $q>1$, and assume that for
$\bar{\alpha} \in (0,1]$, $S \geq 0$ and an open set $\tilde{\Omega}
\CC\Omega$ we have \eqn{crucial}
$$\|\tau_{i,h}w\|_{
L^q(\tilde{\Omega})} \leq S |h|^{\bar{\alpha}}\;, $$ for every
$1\leq i\leq n$ and every $h\in \er$ satisfying $0<|h|\leq d$, where
$0< d \leq \min \{ 1,\dist (\tilde{\Omega},
\partial \Omega )\}$. Then $w\in W^{\alpha
,q}_{\loc}(\tilde{\Omega},\er^k )$ for every $\alpha \in
(0,\bar{\alpha})$, and for each open set $A \CC \tilde{\Omega}$
there exists a constant $c\equiv c(d, \bar{\alpha}-\alpha,\dist
(A,\partial \tilde{\Omega}))$, independent of $S$ and $w$, such that
\eqn{eess}
$$
\| w\|_{W^{\alpha ,q}(A, \er^k)}\leq c\left(S+ \| w\|_{L^{q}(A,
\er^k )}\right).
$$
\end{lemma}
Basic references for the last result are \cite{Ad}, Chapter 7 or
\cite{KM}; see also \cite{ELM3}, Lemma 3, from which the previous
lemma follows via a covering argument. The following result is
nothing but Sobolev's embedding theorem in the case of fractional
spaces; see again \cite{Ad} and also \cite{ELM3}, Lemma 3.
\begin{theorem}[Fractional Sobolev embedding] \label{fraemb}
Let $w\in W^{\alpha,q} (A ,\er^k)$,with $q \geq 1$  and $\alpha \in
(0,1]$, such that $\alpha  q<n$, and let $A\subset \er^n$ be a
Lipschitz domain. Then $w \in L^{nq/(n-\alpha q)}(A ,\er^k)$, and
there exists a constant $c \equiv c(n,\alpha,q, [\partial A]_{0,1})$
such that
$$
\|w\|_{L^{\frac{nq}{n-\alpha q }(A)}}\leq c \|w\|_{W^{\alpha, q
}(A)}\;.
$$
\end{theorem}
For the following fact see for instance \cite{Min1} and related
references.
\begin{prop}[Fractional Poincar\'e inequality]\label{frapoes}
If $B \equiv B(x_0,R) \subset \er^n$ is a ball and $w \in W^{\alpha,
q}(B,\er^k)$, then \eqn{frapo}
$$
\int_{B} \! |w -w_B|^{q} \, dx \leq c(n)R^{\alpha q}\int_{B}
\int_{B} \! \frac{|w(x) - w(y) |^{q}}{|x-y|^{n+\alpha q}} \, dx\,dy
= c(n)R^{\alpha q}[w]_{\alpha,q;B}^q\;.
$$
\end{prop}
The following result is classical in potential theory \cite {AdHe};
see again \cite{Min1} for an elementary proof that avoids potential
theory.
\begin{prop} \label{redu}
Let $w \in W^{\alpha,q}_{\loc}(\Omega, \er^k)$, where $0 < \alpha
<1$, $q \geq  1$ are such that $\alpha q <n$. Let $\Sigma_w$ denote
the set of non-Lebesgue points of $w$ in the sense of
$$
\begin{array}{c}
\displaystyle \Sigma_w:= \left\{ x \in \Omega :
\liminf_{\rho\searrow 0}
 \mean{B(x,\rho)} |w(y)-(w)_{x,\rho}|^q \ dy
 >0\ \displaystyle \mbox{or}\  \limsup_{\rho\searrow 0}
 |(w)_{x,\rho}| =\infty \right\}.
\end{array}
$$
Then its Hausdorff dimension $\dim(\Sigma_w)$ satisfies
$\dim(\Sigma_w)\leq n-\alpha q.$
\end{prop}
\subsection{\bf Morrey spaces, BMO, VMO} We shall adopt a slightly modified
definition of Morrey spaces, or more precisely: there are several
possible, essentially equivalent definitions in the literature; we
choose one. With $A\subset \er^n$ being an open subset, we define
the Morrey space $L^{q,\theta}(A)$, with $q\geq 1$ and $\theta \in
[0,n]$ as that of those measurable maps $w\in L^{q}(A)$ such that
the following quantity is finite: \eqn{dede1mo}
$$
\|w\|_{L^{q,\theta}(A)}^q:= \sup_{B_R \subset A, R\leq 1}
R^{\theta-n} \int_{B_R} |w|^q\, dx\;.
$$
 In the following, when considering the space $\mathcal M
(A)$ of Borel measures with finite mass on $A \subset \er^n$, we
shall automatically consider them extended on the whole $\er^n$ in
the trivial way: $|\mu|(\er^n\setminus A):=0$. When considering
$L^{1,\theta}(A)$, as in \cite{adamsduke}, we include measures $\mu
\in \mathcal M (A)$ defining in this case
$$
\|\mu\|_{L^{1,\theta}(A)}:= \sup_{B_R \subset A, R \leq 1}
R^{\theta-n}|\mu|(B_R) < \infty\;,
$$
and actually $L^{1,\theta}(A)$ will be considered as a subspace of
$\mathcal M (A)$. Trivially, if $\mu$ satisfies \rif{morrey} then
$\mu \in L^{1,\theta}(A)$ for every open subset $A \subset \er^n$
and $\|\mu\|_{L^{1,\theta}(A)}\leq M$. Information on Morrey spaces
are in \cite{adams, G}. Our definition differs from the usual one in
that we consider only balls contained in $A$ when stating
\rif{dede1mo}, and with radius not larger than one, because we shall
treat interior regularity, and information near the boundary
$\partial \Omega$ will play no role. Such a modification is truly
inessential, and will simplify the already heavy technical treatment
in the following pages; observe that our definition is anyway
consistent with the one in \cite{Stamp65}, Definition 1.1.

The following lemma is elementary and can be obtained via a standard
scaling argument; the simple proof is left to the reader.
\begin{lemma}\label{scalamorrey} Let $g \in L^{q,\theta}(B(x_0,r))$ and define
$\tilde{g}(y):=g(x_0+ry)$ for $y \in B(0,1)\equiv B_1$. Then
$\tilde{g} \in L^{q,\theta}(B_1)$ and
$\|\tilde{g}\|_{L^{q,\theta}(B_1)}=r^{-\theta/q}\|g\|_{L^{q,\theta}(B(x_0,r))}$.
\end{lemma}
We now pass to recall the definition of $\BMO$ and $\VMO$ spaces,
introduced in \cite{JoNi, sara} respectively. As already in the case
of Morrey spaces, we shall also modify a bit the definition in order
to adapt it to the local statement we are giving in the following.
The space $\BMO(A)$ is that of those measurable maps $w:A \to \er^n$
such that the semi-norm
$$[w]_{\BMO(A)} := \sup_{B_R \subset A} \mean{B_R} |w - (w)_{B_R}|\, dx $$
is finite.  Further information can be found for instance in
\cite{G}, and its references. Finally the space
$\VMO_{\loc}(\Omega)$. Let $\Omega' \CC \Omega$ be an open subset,
and define
$$\omega_{w}(R, \Omega'):= \sup_{B_r, r\leq R} \mean{B_r} |w - (w)_{B_r}|\, dx\;,$$
where $B_r \CC \Omega$ is centered in $\Omega'$.  Then a map $w :
\Omega \to \er^n$ belongs to $\VMO_{\loc}(\Omega)$ if and only if $
\lim_{R\searrow 0}\omega_{w}(R, \Omega')=0 $ for every open subset
$\Omega' \CC \Omega$. In connection to VMO spaces we shall need the
following:
\begin{definition}\label{lloocc} A Borel measure
$\mu \in \MM (\Omega)$ is said to satisfy the condition in
\rif{locuni} locally uniformly in $\Omega$ iff for every $\Omega'\CC
\Omega$ and every $\ep >0$ there exists $\bar R>0$, depending on
$\ep $ and $\dist(\Omega', \partial \Omega)$, such that $
|\mu|(B_R)\leq \ep R^{n-p},$ whenever $B_R \subset \Omega'$ and $R
\leq \bar R$.
\end{definition}
\begin{remark}\label{abbs} When $p=n$ it is a simple exercise in basic measure theory to check that the measure $\mu$ fulfills Definition \ref{lloocc} iff  $\mu$ has no atoms, i.e.:
$\mu(\{x_0\})=0$ for every $x_0 \in \Omega$. This allows to view the
local $\VMO$ regularity results of \cite{FF} as a particular case of
Theorem \ref{main6}.
\end{remark}
\subsection{\bf Sobolev-Morrey spaces} Beside that of Morrey spaces, we recall
the definition of Sobolev-Morrey spaces of fractional order; also in
this case we propose an inessential modification of the usual
definition to simplify the treatment in the following. Let $A\subset
\er^n$ be an open subset; we say that a map $w \in W^{\alpha, q}(A,
\er^k)$, belongs to $W^{\alpha, q,\theta}(A, \er^k)$, with $\alpha
\in (0,1]$, $q \in [1,\infty)$, $\theta \in [0,n]$, iff $w \in
W^{\alpha, q}(A,\er^k)$, and moreover \eqn{defrac}
$$
[w]_{\alpha, q, \theta;A}^q:=\left\{\begin{array}{ccc} \displaystyle
\sup_{B_R \subset A, R\leq 1} R^{\theta-n} [w]_{\alpha,q;B_R}^q<
\infty & \
\  \mbox{if}\ \ & \alpha < 1  \\ \\
\|Dw\|_{L^{q,\theta}(A)}^q & \ \  \mbox{if}\ \ & \alpha = 1
\end{array}\right. < \infty\;.
$$
In any case we let \eqn{defrac2}
$$
\| w \|_{W^{\alpha,q,\theta}(A )}:= \| w \|_{W^{\alpha,q}(A )}+
[w]_{\alpha, q, \theta;A}\;.
$$
For such spaces and generalizations, see the original papers
\cite{Ca63, Ca64} and \cite{Mazzucato, Taylor}.

\subsection{\bf Marcinkiewicz spaces} Finally, Marcinkiewicz spaces $\mathcal
M^{t}(A,\er^k)$, $t\geq 1$, also called Lorentz-Marcinkiewicz spaces
and denoted by $L^{t,\infty}(A)$, or by $L_{w}^t(A)$, when they are
called ``weak-$L^t$" spaces. A measurable map $w:A \to \er^k$
belongs to $\MM^{t}(A,\er^k)$ iff \eqn{dema}
$$
\sup_{\lambda \geq 0} \lambda^t|\{x \in A \ : \ |w|> \lambda\}|=:
\|w\|_{\MM^{t}(A)}^t< \infty\;.
$$
Yet, we recall the definition of Marcinkiewicz-Morrey spaces
\cite{adamsduke}. A map $w \in \MM^{t}(A,\er^k)$ belongs to the
space $\MM^{t, \theta}(A,\er^k)$ with $\theta \in [0,n]$ iff
$$
 \sup_{B_R \subset A, R\leq 1}
R^{\theta-n} \|w\|_{\MM^{t}(B_R)}^t< \infty\;.
$$
Accordingly, we let \eqn{mama}
$$
\|w\|_{\MM^{t,\theta}(A)}:=\|w\|_{\MM^{t}(A)}+ \left[\sup_{B_R
\subset A, R\leq 1} R^{\theta-n}
\|w\|_{\MM^{t}(B_R)}^t\right]^{1/t}\;.
$$
Obviously $$\|w\|_{\MM^{t,n}(A)} \approx \|w\|_{\MM^{t}(A)}, \qquad
\qquad  \MM^{t,n}(A)\equiv \MM^{t}(A)\;,$$ and, according to the
definition in \rif{dede1mo} \eqn{inclusionimorrey}
$$\MM^{t,\theta_1}(A)\subset  \MM^{t, \theta_2}(A)
\Longleftrightarrow \theta_1 < \theta_2, \quad \mbox{with}\quad
\|w\|_{\MM^{t,\theta_2}(A)}\leq  \|w\|_{\MM^{t,\theta_1}(A)}\;.$$
\begin{lemma}[H\"older's inequality in $\MM^t$]\label{hollo} Let $f \in
\MM^t(A)$ with $t >1$. Then $f  \in L^q(A)$ with $1\leq q < t$ and
it holds \eqn{hollore}
$$
\|f\|_{L^q(A)}\leq \left( \frac{t}{t-q} \right)^{\frac{1}{q}} |A|
^{\frac{1}{q}-\frac{1}{t}}\|f\|_{\MM^t(A)}\;.$$
\end{lemma}
The proof of the latter result is standard \cite{DHM2, zhong}. Next,
a so called ``trace type inequality"  \cite{adams,MZ}.
\begin{theorem}\label{adamo} Let $\lambda$ be a non-negative
Radon measure on $\er^n$ such that $\lambda(B_R)\leq MR^{n-\theta}$,
for every ball $B_R \subset \er^n,$ where $0\leq \theta < p\leq n$
and $M>0$. Then when $p<n$ for every $w \in W^{1,p}_0(B_R)$ we have
\eqn{ada}
$$
\int_{B_R} |w|^{p}\, d \lambda \leq c M R^{p-\theta}\int_{B_R}
|Dw|^{p}\, dx\;,
$$
where $c\equiv c (n,p,\theta)$. In the limit case $p=n$ inequality
\trif{ada} holds replacing $c M R^{n-\theta}$ by $c M R^{\sigma}$,
for any $\sigma <n-\theta$, where $c\equiv c (n,\theta, \sigma)$.
\end{theorem}
\begin{proof}
We did not find any direct reference for this result, therefore we
sketch the proof for the reader's convenience, based on the results
in \cite{adams}. Firstly the case $p<n$. Letting \eqn{pitheta}
$$ p_{\theta}:=\frac{(n-\theta)p}{n-p}, \qquad
\mbox{and} \qquad \tilde{\lambda}:=\lambda/M\;,$$ we obviously have
$\tilde{\lambda} (B_R) \leq R^{n-\theta}$, and then it holds
$$
\left(\int_{B_R} |w|^{p_{\theta}}\, d \tilde{\lambda}
\right)^{\frac{1}{p_{\theta}}}\leq c(n,p) \left(\int_{B_R}
|Dw|^{p}\, dx \right)^{\frac{1}{p}}\;.
$$
This is Adams' inequality, see \cite{MZ}, Corollary 1.93; see also
\cite{AdHe}, comments at Chapter 7 to see the earlier contributions
of Mazy'a, and the original paper of Adams \cite{adams}. Using
H\"older's inequality, as $p_{\theta}\geq p$, we have
$$
\left(\mean{B_R} |w|^{p}\, d \tilde{\lambda}
\right)^{\frac{1}{p}}\leq \left(\mean{B_R} |w|^{p_{\theta}}\, d
\tilde{\lambda} \right)^{\frac{1}{p_{\theta}}}\leq c
\tilde{\lambda}(B_R)^{-\frac{1}{p_{\theta}}} \left(\int_{B_R}
|Dw|^{p}\, dx \right)^{\frac{1}{p}}\;,
$$
therefore, using again that $\tilde{\lambda} (B_R) \leq
R^{n-\theta}$ and \rif{pitheta} we have
$$\int_{B_R} |w|^{p}\, d \tilde{\lambda} \leq c  \tilde{\lambda}(B_R)^{\frac{p-\theta}{n-\theta}}
\int_{B_R} |Dw|^{p}\, dx\leq c R^{p-\theta} \int_{B_R} |Dw|^{p}\, dx
\;, $$ and \rif{ada} follows scaling back to $\lambda$. Now we treat
the case $p=n$. In this case observe that $ \tilde{\lambda}(B_R)
\leq c\log^{q(1-n)/n}(R^{-1})$ when $R\leq 1/2$, and
$\tilde{\lambda}(B_R) \leq 2^{\sigma-q}R^q$ when $R> 1/2$, where
$q>n$, and $c\equiv(n,p,\theta,q)$ is a suitable constant. Therefore
we may apply Theorem 1.94 from \cite{MZ}, see also \cite{AdHe},
Theorem 7.2.2, to have
$$
\left(\int_{B_R} |w|^{q}\, d \lambda \right)^{\frac{1}{q}}\leq
c(n,p,\theta,q) \left(\int_{B_R} |Dw|^{n}\, dx
\right)^{\frac{1}{n}}\;.
$$
Applying H\"older's inequality and $\tilde{\lambda} (B_R) \leq
R^{n-\theta}$ yields, with $c \equiv c(n,p,\theta,q)$
$$
\int_{B_R} |w|^{n}\, d \lambda\leq c(n,q,\sigma)
\tilde{\lambda}(B_R)^{\frac{q-n}{q}}\int_{B_R} |Dw|^{n}\, dx\leq
cR^{(n-\theta)\left(1-\frac{n}{q}\right)} \int_{B_R} |Dw|^{n}\,
dx\;.
$$
The assertion follows taking $q \equiv q(\sigma)$ large enough, and
scaling back to $\lambda$.
\end{proof}
\begin{remark} In a similar way, if $w \in W^{1,p}_0(\Omega)$ with $c\equiv c (n,p,\theta,\Omega)$ we have\eqn{ada2}
$$
\int_{\Omega} |w|^{p}\, d \lambda \leq c M \int_{\Omega} |Dw|^{p}\,
dx\;.
$$
\end{remark}
\subsection{Technical lemmata}
The following is a simple variant of a well known iteration result.
See for instance \cite{G}, Chapter 7, or \cite{zhong}, last section.
\begin{lemma} \label{lemmaiterazione}Let $\phi: [0,\bar{R}]\to [0,\infty)$ be a non-decreasing function such that
$$
\phi(\varrho) \leq c_0
\left(\frac{\varrho}{R}\right)^{\delta_0}\varphi(R) + \B
R^{\gamma}\;, \qquad \mbox{for every}\ \ \varrho < R \leq
\bar{R},\quad \B \geq 0\;,
$$
where $ \gamma \in (0,\delta_0)$. Then if $\delta_1 \in [\gamma,
\delta_0)$, there exists $c_1 \equiv c_1(c_0, \delta_1, \gamma)$
such that
$$
\phi(\varrho) \leq c_1
\left(\frac{\varrho}{R}\right)^{\delta_1}\varphi(R) + c_1\B
\varrho^{\gamma}\;, \qquad \mbox{for every}\  \ \varrho \leq R \leq
\bar{R}\;.
$$
\end{lemma}
Then, Giaquinta \& Giusti's ``simple but fundamental lemma",
\cite{G}, Chapter 6.
\begin{lemma} \label{simpfun}Let $\varphi: [R_0, 2R_0]\to [0,\infty)$ be a function such that
$$
\varphi(t) \leq (1/2) \varphi(\varrho) + \B\dro^{-\beta} +K\;,
\qquad \mbox{for every}\ \ R_0 < t < \varrho < 2R_0\;,
$$
where $\B, K\geq 0$ and $\beta >0$. Then $ \varphi(R_0) \leq
c(\beta)\B R_0^{-\beta}+cK$.
\end{lemma}
Finally, a standard fact.
\begin{lemma} \label{stanom}Let $\Omega' \CC \Omega \subset \er^n$ be bounded
domains. There exists another open subset $\Omega''$ such that
$\Omega' \CC \Omega'' \CC \Omega$ and \eqn{distanze}
$$
\dist(\Omega',\partial \Omega'')=\dist(\Omega'',\partial
\Omega)=\dist(\Omega',\partial \Omega)/2\;.
$$
\end{lemma}
\begin{proof}
Just let $ \Omega'' :=\{x \in \Omega \ : \ \dist(x,
\overline{\Omega'})< 1/2\ \}. $
\end{proof}
\section{Regularity for homogeneous problems}
In this section we recall some results on the regularity of
solutions to homogeneous elliptic systems and equations with
$p$-growth; some of them are well-known; some others, much less if
not at all, especially in the explicit form needed in this paper. In
such cases we shall give - sometimes sketchy - proofs; anyway a good
general reference for the whole section is \cite{G}, Chapters 6 and
7.

Let us start with a simple but rather rarely used lemma on reverse
H\"older inequalities. For the proof it suffices to follow Remark
6.12, page 205 in \cite{G}; see also \cite{BI} for this kind of
result.
\begin{lemma}\label{revq}
Let $g : A \to \er^k$ be a measurable map, and $\chi_0>1$, $c,s \geq
0$, such that
$$
\left(\mean{B_R} |g|^{\chi_0} \, dx \right)^{\frac{1}{\chi_0}}\leq c
\mean{B_{2R}} (|g|+s) \, dx\;,
$$
whenever $B_{2R} \CC A$, where $A \subset \er^n$ is an open subset.
Then, for every $t \in (0,1]$ there exists a constant $c_0 \equiv
c_0(n,c,t)$ such that, for every $B_{2R} \subseteq A$
$$
\left(\mean{B_R} |g|^{\chi_0} \, dx \right)^{\frac{1}{\chi_0}}\leq
c_0 \left(\mean{B_{2R}} (|g|^t+s^t) \, dx \right)^{\frac{1}{t}}\;.
$$
\end{lemma}
The next two lemmata will be of fundamental importance in the
following in that they provide {\em estimates below the natural
growth exponent} $p$. For reasons that will become clear in Section
11 the first one is stated directly for systems, that is when $u$ is
a vector valued map and therefore $N\geq 1$.
\begin{lemma}\label{revbase}
Let $v_0\in W^{1,p}(A, \er^{N})$ be a weak solution to the system
$$\divo \ a_0(Dv_0)=0 \qquad \qquad \mbox{in} \quad  A\;.$$ Here $a_0 : \Ma \to
\Ma$ satisfies the assumptions \trif{asp} obviously recast to fit
the vectorial case with no $x$-dependence, and $A \subset \er^n$ is
an open subset. Then $V(Dv_0) \in W^{1,2}_{\loc}(A, \Ma)$, and there
exists $c \equiv c(n,N,p,\ratio)$ such that for every $z_0 \in \Ma$
and every ball $B_R \subseteq A$, we have \eqn{cacc2}
$$\int_{B_{R/2}} \! |D(
V(Dv_0))|^{2} \, dx \leq \frac{c}{R^2} \int_{B_R} \!
|V(Dv_0)-V(z_0)|^{2} \, dx\;.
$$
Moreover, for every $t \in (0,1]$ there exists  $c \equiv
c(n,N,p,\ratio,t)$ such that \eqn{revca}
$$\left(\mean{B_{R/2}} \! |
V(Dv_0)-V(z_0)|^{2} \, dx \right)^{\frac{1}{2}}\leq c
\left(\mean{B_R} \! |V(Dv_0)-V(z_0)|^{2t} \,
dx\right)^{\frac{1}{2t}}\;.
$$
All the constants named $c$ involved in \trif{cacc2}-\trif{revca}
are independent of the choice of $z_0 \in \Ma$.
\end{lemma}
\begin{proof} {\em Step 1: Regularization}. We first regularize the problem
following a few smoothing arguments similar to those in \cite{ELM2}.
We consider a standard, symmetric and non-negative mollifier $\phi :
\Ma \to \er$, such that $\phi \in C^{\infty}_0(B_1)$, and
$\|\phi\|_{L^1(\Ma)}=1$.
 Moreover, in order to apply the technique of \cite{ELM2} we also need a technical assumption, namely we need to take
 $\phi$
 such that
\eqn{stranamol}
$$
\int_{B_1\setminus B_{1/2}} \phi(z)\, dz \geq 1/1000\;.
$$
For every $k \in \en$ set $\phi_k(z):=k^{\Ma}\phi(kz)$, and then
define the smooth vector field $a_{k}(z)$ via convolution as
follows:
$$a_{k}(z):=(a_0*\phi_k)(z):=\int_{B(0,1)}a_0(z+k^{-1} y)\phi(y)\,
dy\;.$$ Assumptions \rif{asp} and a few convolution estimates also
using \rif{stranamol}, similar to those of \cite{ELM2}, Lemma 3.1,
imply that each $a_k$ satisfies
 \eqn{asp2}
$$
\left\{
    \begin{array}{c}
    |a_k(z)|+|\partial_{z}a_k(z)|(s_k^2+|z|^2)^{\frac{1}{2}} \leq c(s_k^2+|z|^2)^{\frac{p-1}{2}} \\ \\
    c^{-1}(s_k^2+|z|^2)^{\frac{p-2}{2}}|\lambda|^{2} \leq \langle\partial_{z}a_k(z)\lambda, \lambda
    \rangle\\ \\
    |a_0(z)-a_k(z)|\leq
    ck^{-1}(s_k^2+|z|^2)^{\frac{p-2}{2}}\;,
    \end{array}
    \right.
$$ whenever $z, \lambda \in \Ma$, where $s_k:=s+k^{-1}$, $k \in \en$, and $c\equiv c(n,N,p,\ratio)$.
Note that the new ellipticity/growth constant $c$ is actually
independent of $k \in \en$. Moreover each $a_k$ satisfies the
assumptions \rif{asp} with $s$ replaced by $s_k$, for different
constants $\nu, L$, but still depending on the original ones. This
fact and standard monotonicity methods \cite{lions} allow to define,
with $B_R \CC \Omega$ as in the statement, $v_k \in
v_0+W^{1,p}_0(B_R)$ as the unique solution to \eqn{Dirckappa}
$$
\left\{
    \begin{array}{cc}
    -\divo \ a_k(Dv_k)=0 & \qquad \mbox{in} \ B_R\\
        v_k= v_0&\qquad \mbox{on $\partial B_R$.}
\end{array}\right.
$$
\indent {\em Step 2: Estimates}. Under assumptions
\rif{asp2}$_{1,2}$ the proof of Caccioppoli's type inequality
\rif{cacc2} with $c\equiv c(n,N,p,\ratio)$, $V(Dv_0) \equiv
V_{s_k}(Dv_k)$, and any ball $B_r \subseteq B_R$, can be inferred
from \cite{Ca}, Theorem 1.1, with minor variants, see also
\cite{Ham, KM, ELM2}. As for the proof of \rif{revca}, set $$ \chi_0
:=\left\{
\begin{array}{ccc}\frac{n}{n-2} & \mbox{if}&  n>2\\ \\
2 & \mbox{if}&  n=2 \end{array}>1\;. \right.$$ Using a simple
scaling argument, and applying Sobolev embedding theorem to the map
$V_{s_k}(Dv_k)-V_{s_k}(z_0)$, we get that there exists a constant
$c\equiv c(n)$ such that for any ball $B_r \subseteq B_R$
\begin{eqnarray*}
\left( \mean{B_{r/2}} \! |V_{s_k}(Dv_k)-V_{s_k}(z_0)|^{2\chi_0} \,
dx\right)^{\frac{1}{\chi_0}} &\leq & c  \mean{B_{r/2}} \!
|V_{s_k}(Dv_k)-V_{s_k}(z_0)|^{2} \, dx \\ && \qquad +
cr^2\mean{B_{r/2}} \! |D( V_{s_k}(Dv_k))|^{2} \, dx\;.
\end{eqnarray*}
We now use \rif{cacc2} with $V(Dv_0)\equiv V_{s_k}(Dv_k)$ for the
last integral, thereby getting
$$
\left( \mean{B_{r/2}} \! |V_{s_k}(Dv_k)-V_{s_k}(z_0)|^{2\chi_0} \,
dx\right)^{\frac{1}{\chi_0}}\leq c \mean{B_{r}} \!
|V_{s_k}(Dv_k)-V_{s_k}(z_0)|^{2} \, dx\;.
$$
In the last two inequalities it is $c\equiv c(n,N,p,\ratio)$.
Inequality \rif{asp2} for $V(Dv_0) \equiv V_{s_k}(Dv_k)$ now follows
from Lemma \ref{revq}, and then H\"older's inequality again.

{\em Step 3: Approximation.}  Assumptions \rif{asp2}$_{1,2}$ imply
in a rather standard way that $a_k(z)$ satisfy the growth and
monotonicity assumptions \rif{veramon}-\rif{veracre} with $s$
replaced by $s_k$, uniformly with respect to $k \in \en$. In turn
this implies that $v_k$ is a so-called $Q$-minimum of the functional
\eqn{funcase}
$$w \mapsto
\int_{B_R}(|Dw|^p+s^p +k^{-p})\,dx$$ with $Q \equiv
Q(n,N,p,\ratio)\geq 1$ being independent of $k \in \en$; for such a
conclusion see Theorem 6.1 from \cite{G} applied to the functional
in \rif{funcase}, when, with the notation of \cite{G}, it is
$a_1(x)\equiv [a_2(x)]^{p/(p-1)}\equiv s_k^p $. The $Q$-minimality
of $v_k$ now easily yields \eqn{approssimano2}
$$
\|Dv_k\|_{L^p(B_R)}\leq c(n,N,p,\ratio) \||Dv_0|+s+1\|_{L^p(B_R)}\;.
$$
Therefore, up to a non-relabeled subsequence we may assume that
$\{v_k\}_k$ weakly converges to some map in $W^{1,p}_0(B_R)$;
actually we may assume that $v_k \to v_0$ strongly in $W^{1,p}(B_R)$
too. Indeed, using that both $v_0$ and $v_k$ are solutions, and that
$v_0 \equiv v_k$ on $\partial B_R$, and making also use of
\rif{mon2}, we have
\begin{eqnarray*}
\int_{B_R} |Dv_k-Dv_0|^p \, dx &\leq &
c\int_{B_R}\langle a_0(Dv_k)-a_0(Dv_0),Dv_k-Dv_0\rangle  \, dx\\
&= &  c\int_{B_R}\langle  a_0(Dv_k)-a_k(Dv_k),Dv_k-Dv_0\rangle  \, dx\\
&\leq  & \frac{1}{2}\int_{B_R}|Dv_k-Dv_0|^p \, dx \\ && \qquad
\qquad + c\int_{B_R} |a_0(Dv_k)-a_k(Dv_k)|^{\frac{p}{p-1}}  \, dx\;.
\end{eqnarray*}
The last integral tends to zero as $k \nearrow \infty$ by
\rif{asp2}$_3$ and \rif{approssimano2}. Therefore $Dv_k \to Dv_0$
strongly in $L^p(B_R)$, and since all the $v_k,v_0$ share the same
boundary datum it follows $v_k \to v_0$ strongly in $W^{1,p}(B_R)$.
In turn this and \rif{elemV} imply that up to another subsequence
$V_{s_k}(Dv_k) \to V(Dv_0)$ strongly in $L^2(B_R)$. Fix $z_0 \in \Ma
$ as in the statement; applying estimate \rif{cacc2} to
$V_{s_k}(Dv_k)$ instead of $V(Dv_0)$, which is allowed by the
previous step, we infer that $\{V_{s_k}(Dv_k)\}_k$ is bounded in
$W^{1,2}(B_{R/2})$ and therefore it also holds $V_{s_k}(Dv_k)
\rightharpoonup V(Dv_0)$ weakly in $W^{1,2}(B_{R/2})$ up to yet
another subsequence. We are ready to conclude: writing estimate
\rif{cacc2} with $V(Dv_0) \equiv V_{s_k}(Dv_k)$ and letting $k
\nearrow \infty$ we find the final form of \rif{cacc2} for the
original $V(Dv_0)$ using strong convergence for the right hand side,
and lower semicontinuity for the left hand one. From this last fact
the inclusion $V(Dv_0) \in W^{1,2}_{\loc}(A)$ in turn follows via a
covering argument. In the same way \rif{revca} follows from the
similar inequality for the $V_{s_k}(Dv_k)$ given in Step 2.
\end{proof}
Finally, a few basic consequences of De Giorgi's regularity theory
for elliptic equations, and Gehring's lemma for elliptic problems
and variational integrals; see for instance \cite{G}, Chapters 6-7
for a reasonable overview of the subject.
\begin{lemma}\label{superrev}
Let $v \in W^{1,p}(A)$ with $p \in (1,n]$ be a weak solution to the
equation \eqn{roueq}
$$\divo \ a(x,Dv)=0 \qquad \qquad \mbox{in} \quad  A\,,$$ under the assumptions
\eqn{monass}
$$|a(x,z)|\leq c(s^2+|z|^2)^{\frac{p-1}{2}}, \qquad \qquad
c^{-1}|z|^p -cs^p \leq \langle a(x,z),z \rangle \;,$$ for every $x
\in \Omega$ and $z \in \er^n$, where $c \equiv c(\ratio)$ and $\nu,
L$ are the numbers given in \trif{asp}. There exists $\beta \equiv
\beta(n,p,\ratio) \in (0,1]$, such that for every $q \in (0,p]$
there exists  $c \equiv c(n,p,\ratio,q)$ such that, whenever $B_R
\subseteq A$, and $0<\varrho \leq R$ it holds  \eqn{decmo}
$$
\int_{B_\varrho} ( |Dv|^q +s^q ) \, dx \leq c
\left(\frac{\varrho}{R}\right)^{n-q+\beta q} \int_{B_R} ( |Dv|^q
+s^q )\, dx\;.
$$
Moreover, there exists $\chi \equiv \chi(n,p,\ratio)>1$, such that
$Dv \in L^{p\chi}_{\loc}(A,\er^n)$ and\eqn{zz2}
$$
\left(\mean{B_{R/2}} |Dv|^{p\chi} \, dx\right)^{\frac{1}{p\chi}}\leq
c \left(\mean{B_R} (|Dv|^q +s^q)\, dx \right)^{\frac{1}{q}}\;,
$$
where again $c \equiv c(n,p,\ratio,q)$.
\end{lemma}
\begin{proof}
First observe that by \rif{monass} we may apply Theorem 6.1 from
\cite{G} with the choice $a_1(x)\equiv [a_2(x)]^{p/(p-1)}\equiv s^p
$, concluding that the solution $v$ is a $Q$-minimum of the
functional \eqn{partfun}
$$w \mapsto \int_{A}(|Dw|^p+s^p)\, dx\,$$
with $Q\equiv Q(n,p,\ratio)\geq 1$. This in turn allows to use
Theorem 6.7 from \cite{G} that in the special case of the functional
in \rif{partfun} applies with the choice $\theta(x,u)\equiv s$, and
ultimately yields the existence of higher integrability exponent
$\chi \equiv \chi (n,p,\ratio)>1$ such that $Dv \in
L^{p\chi}_{\loc}(A,\er^n)$; moreover \rif{zz2} holds for $q=p$. In
turn \rif{zz2} for every $q \in (0,p]$ follows applying Lemma
\ref{revq} with the choice $\chi_0 \equiv \chi$. In particular
H\"older's inequality gives \eqn{zz2palle}
$$
\left(\mean{B_{R/2}} |Dv|^{p} \, dx\right)^{\frac{q}{p}}\leq c
\mean{B_R} (|Dv|^q +s^q)\, dx \;.
$$
 It remans to
establish \rif{decmo}; this inequality is standard when $q=p$, that
is \eqn{zz1}
$$\mean{B_\varrho} |Dv|^p \, dx \leq c \left(\frac{\varrho}{R}\right)^{-p+\beta p}\mean{B_R} |Dv|^p \, dx
+c \left(\frac{\varrho}{R}\right)^{-p+\beta p}s^p\;,$$ where $\beta
>0$ is as specified in the statement. For the proof of \rif{zz1} see Remark \ref{esore} below.
Therefore, when $\varrho \in (0,R/2]$, using H\"older's inequality
yields
\begin{eqnarray*}
\mean{B_{\varrho}} |Dv|^q \, dx &\leq & \left(\mean{B_{\varrho}}
|Dv|^p \, dx\right)^{\frac{q}{p}}\\ & \stackrel{\rif{zz1}}{\leq}& c
\left(\frac{\varrho}{R}\right)^{-q+\beta q}\left(\mean{B_{R/2}}
|Dv|^p \, dx\right)^{\frac{q}{p}} +c
\left(\frac{\varrho}{R}\right)^{-q+\beta q}s^q\\&
\stackrel{\rif{zz2palle}}{\leq} &c
\left(\frac{\varrho}{R}\right)^{-q+\beta q}\mean{B_{R}} (|Dv|^q
+s^q)\, dx\;.
\end{eqnarray*}
Summing $s^q$ to both sides of the previous inequality, taking into
account that $\varrho \leq R$ and $q-\beta q\geq 0$, and finally
getting rid of the averages
 gives \rif{decmo} when $\varrho \in
(0,R/2]$; the case $\varrho \in (R/2,R]$ trivially follows and the
lemma is completely proved.
\end{proof}
\begin{remark}[An esoteric detail]\label{esore} By carefully tracing
the dependence of the constants back in De Giorgi's theory - see in
particular \cite{G}, Paragraphs 7.1-7.3 - we are giving here a
justification of inequality \rif{zz1}. Using Theorem 7.7 from
\cite{G}, applied to the particular functional in \rif{partfun} when
$a(x)\equiv s^p$, and taking into account Remark 7.7 again from
\cite{G}, we have that \rif{zz1} holds in the preliminary form
\eqn{decmoeso}
$$\mean{B_\varrho} |Dv|^p \, dx \leq c
\left(\frac{\varrho}{R}\right)^{-p+\beta p}\mean{B_R} |Dv|^p \, dx
+c \|s^p\|_{L^t(B_R)}\varrho^{-p+\beta p}\;,$$ where $\beta :=
\min\{\tilde{\beta}(n,p,\ratio), n\ep/p\}$ and
$\tilde{\beta}(n,p,\ratio)>0$, $1/t:=p/n-\ep$; here $\ep \in
(0,p/n)$ can be picked arbitrarily small. In fact, choose $\ep\equiv
\ep(n,p,\ratio)$ small enough in order to have that $\beta =
n\ep/p$; then
$$
\|s^p\|_{L^t(B_R)}\varrho^{-p+\beta p}\leq
c(n,p)s^pR^{p-n\ep}\varrho^{-p+\beta
p}=cs^p\left(\frac{\varrho}{R}\right)^{-p+\beta p}\;.
$$
Merging the latter inequality with \rif{decmoeso} yields \rif{zz1}.
\end{remark}

\section{Comparison estimates}
Let us first introduce some notation that we shall keep for the rest
of the paper; accordingly to \rif{ssii} and \rif{ssii2}, in the case
$\theta \in [p,n]$ in \rif{morrey} we define \eqn{ssiidef}
$$
\left\{\begin{array}{c}
 \sigma(q,\theta):=
    \theta(1- q/m) = \theta - \frac{q(\theta-1)}{p-1} \qquad \mbox{when} \qquad p-1 \leq q < m\\
    \\
    \sigma(q):=n - \frac{q(n-1)}{p-1} \qquad \mbox{when} \qquad 1 \leq q < b\;.
\end{array}\right.
$$
     Here, as in the rest of the paper, $b$ will denote the
number defined in \rif{boccardoex}, and $m$ the one in
\rif{mingione}. For the rest of the section we fix a ball $$B_R
\equiv B(x_0,R) \CC \Omega, \qquad \mbox{with}\qquad R\leq 1\;.$$
The first two lemmas are dealing with weak solutions to more regular
problems i.e.~$u \in W^{1,p}_0(\Omega)$ will be the unique solution
to \eqn{Dirapp2}
$$
\left\{
    \begin{array}{cc}
    -\divo \ a(x,Du)=f \in L^{\infty}(\Omega) & \qquad \mbox{in $\Omega$}\\
        u= 0&\qquad \mbox{on $\partial\Omega$,}
\end{array}\right.
$$
for a fixed $f$ to be eventually chosen; such a solution exists via
standard monotonicity methods \cite{lions} as $f$ belongs to the
dual of $W^{1,p}$. By the same argument we introduce $v \in u +
W^{1,p}_0(B_R)$, defined as the unique solution to \eqn{Dirc1}
$$
\left\{
    \begin{array}{cc}
    -\divo \ a(x,Dv)=0 & \qquad \mbox{in} \ B_R\\
        v= u&\qquad \mbox{on $\partial B_R$.}
\end{array}\right.
$$
\begin{lemma} \label{coco1} Under the assumptions \trif{asp} with $p\leq n$, with $u \in W^{1,p}(B_R)$ as in \trif{Dirapp2},
and $ v \in u +W^{1,p}_0(B_R)$
 as in \trif{Dirc1}, we have
for any $1 \leq q <
   b$ that \eqn{comp1}
$$
\int_{B_R} \VVq+|Du-Dv|^q\, dx \leq c \left( \int_{B_R} |f|\, dx
\right)^{\frac{q}{p-1}} R^{\sigma(q)}\;,
$$
where $\sigma(q)$ is in \trif{ssiidef}, and $c \equiv c(n,p,\nu,q)$.
\end{lemma}
\begin{proof} {\em Step 1}: Here we observe that we can
assume $B(x_0,R) \equiv B_1$ by a scaling argument. Indeed, changing
variables, we let for $y \in B_1$ \eqn{scalaraggio}
$$\left\{
\begin{array}{ccc}
\tilde{u}(y) := R^{-1}u(x_0+Ry), &\tilde{v}(y) := R^{-1}v(x_0+Ry),
&\\ \\
\tilde{a}(y,z) := a(x_0+Ry,z),& \tilde{f}(y) := Rf(x_0+Ry), & B_R
\equiv B(x_0,R),\\ \\-\divo\ \tilde{a}(y,D\tilde{u}) =
\tilde{f},&-\divo\ \tilde{a}(y,D\tilde{v}) =0 \;. &
\end{array}\right.
$$
Obviously $\tilde{u} \equiv \tilde{v}$ on $\partial B_1$. It is
readily verified that the new vector field $\tilde{a}$ still
satisfies \rif{asp}. Now \rif{comp1} follows by writing its
corresponding version for $R=1$, and scaling back to $B_R$.

{\em Step 2}: Here we prove the following implication: \eqn{unob}
$$\int_{B_1} |f|\, dx\leq 1\Longrightarrow \int_{B_1} \VVq+ |Du-Dv|^q\, dx \leq
c_2\;, $$with $c_2 \equiv c_2(n,p,\nu,q)$. Notice that the following
computations remain valid also for $q \in [1,p-1)$. In order to
prove \rif{unob} we shall revisit the technique of \cite{BG1},
reporting the necessary modifications in some detail for the sake of
clarity. For $k >  0$, let us define the following truncation
operators, classical after \cite{BG1}: \eqn{troncamenti}
$$
T_k(s):= \max \{-k, \min\{k,s\}\}, \quad \quad
\Phi_k(s):=T_1(s-T_k(s)), \quad s \in \er\;.
$$
Since both $u$ and $v$ are solutions, we test the weak formulation
\eqn{weak}
$$
\int_{B_1} \langle a(x,Du) -a(x,Dv),D \varphi \rangle\, dx =
\int_{B_1} f  \varphi\, dx \;,
$$
with $\varphi\equiv T_k(u-v)$; this function is admissible as it is
in $L^\infty(B_R)\cap W^{1, p}_0(B_R)$, and we have \rif{veracre}.
Using the monotonicity inequalities \rif{mon3}-\rif{mon2}, and the
bound in \rif{unob}, we easily obtain with $c\equiv c(n,p,\nu)$
\eqn{boc01}
$$
\int_{D_k} |V(Du)-V(Dv)|^2+|Du-Dv|^p\, dx \leq ck\int_{B_1} |f|\, dx
\leq ck \;.
$$
Here we have set \eqn{dkappa}
$$D_k:=\{x \in B_1  :  |u(x)-v(x)|\leq k\}\;.$$
Moreover, testing again \rif{weak} with $\varphi \equiv
\Phi_k(u-v)$, and again using \rif{mon3}-\rif{mon2} and the bound in
\rif{unob}, we obtain \eqn{boc1}
$$
\int_{C_k}\VVd+ |Du-Dv|^p\, dx \leq c\int_{B_1} |f|\, dx
\stackrel{\rif{unob}}{\leq }c\;,
$$
where this time \eqn{ckappa}
$$C_k:=\{x \in B_1  : k  < |u(x)-v(x)|
\leq k+1\}\;,$$ and $c\equiv c(n,p,\nu)$. By H\"older's inequality,
and the very definition of $C_k$, we find
\begin{eqnarray}
&&\nonumber \int_{C_k} \VVq+ |Du-Dv|^q\, dx \\&  & \qquad \qquad
\leq c|C_k|^{1-\frac{q}{p}} \left(\int_{C_k} \VVd+ |Du-Dv|^p\,
dx\right)^{\frac{q}{p}} \notag  \\ & & \qquad
\qquad\stackrel{\rif{boc1}}{\leq }c|C_k|^{1-\frac{q}{p}}\leq
\frac{c}{k^{q^*\left(1-\frac{q}{p}\right)}}\left(\int_{C_k}
|u-v|^{q^*}\, dx\right)^{1-\frac{q}{p}}\;.\label{gal}
\end{eqnarray}
With $0\not=k_0 \in \en $ to be fixed later, we have, using the
previous inequality and H\"older's inequality for sequences, and
finally Sobolev's embedding theorem, as $q<n$ under the present
assumptions:
\begin{eqnarray}
\nonumber && \int_{B_1} \VVq +|Du-Dv|^q\, dx\\\nonumber  &  & \qquad
= \int_{D_{k_0}} \VVq+ |Du-Dv|^q\, dx \\&& \nonumber \qquad \qquad +
\sum _{k=k_0}^{\infty}\int_{C_k} \VVq+ |Du-Dv|^q\, dx \nonumber
\\&
  &\qquad \qquad \qquad  \stackrel{\rif{boc01}}{\leq }ck_0 + c\left[\sum
_{k=k_0}^{\infty}\frac{1}{k^{q^*\left(\frac{p-q}{q}\right)}}\right]^{\frac{q}{p}}\left(\int_{B_1}
|u-v|^{q^*}\, dx\right)^{1-\frac{q}{p}} \nonumber \\& &  \qquad
\qquad \qquad \leq ck_0 + cH(k_0)\left(\int_{B_1} |Du-Dv|^{q}\,
dx\right)^{\frac{q^*}{q}\left(1-\frac{q}{p}\right)}\;,\label{semi}
\end{eqnarray}
with $$ H(k_0):=\left[\sum
_{k=k_0}^{\infty}\frac{1}{k^{q^*\left(\frac{p-q}{q}\right)}}\right]^{\frac{q}{p}},
\qquad \mbox{and}\quad c\equiv c(n,p,\nu,q)\;.$$  Here $ H(k_0)$ is
finite since $q   < b $ implies that $ q^*(p/q-1)>1$. We finally
distinguish two cases. If $p<n$ then we take $k_0=1$ in \rif{semi},
and observe that $\gamma:= (q^*/q)(1-q/p)<1$. Therefore, applying
Young's inequality in \rif{semi} with conjugate exponents $1/\gamma$
and $1/(1-\gamma)$ we find \rif{unob}. In the case $p=n$ we have
that $\gamma=1$ and the previous argument does not work; instead, we
choose $k_0$ large enough in order to have $cH(k_0)=1/2$ in
\rif{semi}, and \rif{unob} follows again. Observe that this
determines $k_0\equiv k_0(n,p,\nu,q)$ possibly large, and this
finally reflects in the constant $c$ appearing in \rif{unob}.

{\em Step 3.} We are ready to conclude the whole proof, again by
mean of a scaling argument. We shall prove the validity of the
estimate for $B_R \equiv B_1$, and then we shall conclude using Step
1. Without loss of generality we assume that $A :=
\|f\|_{L^1(B_1)}^{1/(p-1)}>0$, otherwise $u\equiv v$ and the
assertion is trivially verified. We define the new solutions
$\tilde{u}:=A^{-1}u$, $\tilde{v}:=A^{-1}v,$ the new datum
$\tilde{f}:=A^{1-p}f$, and the new vector field
$\tilde{a}(x,z):=A^{1-p}a(x,Az).$ Therefore we have that
$\tilde{u}\equiv \tilde{v}$ on $\partial B_1$, and moreover $\divo \
\tilde{a}(x,D\tilde{u})=\tilde{f}$, $\divo\
\tilde{a}(x,D\tilde{v})=0$, in the weak sense. We make sure that we
can apply the result in Step 2. Trivially
$\|\tilde{f}\|_{L^1(B_1)}=1$ and moreover it is easy to see that the
vector field $\tilde{a}(x,z)$ satisfies \rif{asp} with $s$ replaced
by $s/A\geq 0$. Therefore the inequality in \rif{unob} holds in the
form
$$
\int_{B_1} |V_{s/A}(D\tilde{u})-V_{s/A}(D\tilde{v})|^{2q/p}+
|D\tilde{u}-D\tilde{v}|^q \, dx \leq c_2, \qquad  c_2 \equiv
c_2(n,p,\nu,q)\;.
$$
 Re-scaling back from
$\tilde{u}-\tilde{v}$ to $u-v$ and using \rif{elemV2}, we find
$$
\int_{B_1} \VVq+ |Du-Dv|^q \, dx  \leq c_2 \left( \int_{B_1} |f|\,
dx \right)^{\frac{q}{p-1}} \;,
$$
and the proof is concluded via Step 1.
\end{proof}
\begin{remark} Although the previous lemma has been stated for $q\geq
1$ we shall use it only for the case $q\geq p-1$.
\end{remark}
\begin{lemma} \label{coco2} Under the assumptions \trif{asp} with $p\leq n$, assume $p-1 \leq q <
   b,$ and $f \in L^{1,\theta}(B_R)$.
 With $u \in W^{1,p}_0(B_R)$ as in \trif{Dirapp2},
and $ v \in u +W^{1,p}_0(B_R)$
 as in \trif{Dirc1}, we have
for any $R\leq 1$ \eqn{comp1m}
$$
\int_{B_R} |V(Du)-V(Dv)|^{\frac{2q}{p}}+|Du-Dv|^q\,dx \leq c
\|f\|_{L^{1,\theta}(B_R)}^{\frac{q-p+1}{p-1}} \int_{B_R} |f|\, dx \,
R^{\sigma(q,\theta)}\;,
$$
where $\sigma(q,\theta)$ is in \trif{ssiidef}, and $c \equiv
c(n,p,\nu,q)$.
\end{lemma}
\begin{proof}
First observe that the definition in \rif{ssiidef} implies
\eqn{reprep}
$$
(n-\theta)
\left(\frac{q}{p-1}-1\right)+\sigma(q,n)=\sigma(q,\theta)\;.
$$ Now, since $p-1\leq q$ we may estimate
\eqn{rep2}
$$
\left( \int_{B_R} |f|\, dx \right)^{\frac{q}{p-1}}  \leq
R^{(n-\theta) \left(\frac{q}{p-1}-1\right)}
\|f\|_{L^{1,\theta}(B_R)}^{\frac{q-p+1}{p-1}} \int_{B_R} |f|\, dx\;,
$$
and then we conclude by merging \rif{rep2} with \rif{comp1}, taking
\rif{reprep} into account.
\end{proof} The next twin lemmata are about the capacitary case
$\theta <p$. In the following $u$ will be the solution to
\rif{Dir1}, and $\mu$ the Radon measure in \rif{Dir1}$_1$. We have $
u\in W^{1,p}_0(\Omega)$, and $u$ is the unique solution, since under
the assumptions considered in the next two lemmata it is $\mu \in
W^{-1,p'}(\Omega)$ by a theorem of D.~R.~Adams \cite{adams, AdHe}.
\begin{lemma} \label{coco5} Under the assumptions \trif{asp} with $p >
n$, and  with $u, v \in W^{1,p}(B_R)$ as in \trif{Dir1} and
\trif{Dirc1} respectively, if \trif{morrey} holds then \eqn{baba}
$$
\int_{B_R} \VVd\, dx\leq c
M^{\frac{1}{p-1}}|\mu|(B_R)R^{\sigma(p)}\, \;,
$$
where $\sigma(p)=(p-\theta)/(p-1)$ is as in \trif{fra123}, and
$c\equiv c(n,p,\nu)$.
\end{lemma}
\begin{proof}
Notice that here it can be also $\theta=n$. We test the weak
formulation \eqn{weak2}
$$
\int_{B_R} \langle a(x,Du) -a(x,Dv),D \varphi \rangle\, dx =
\int_{B_R} \varphi \, d\mu \;,
$$
with $\varphi\equiv u-v$, which is admissible as $p>n$ and therefore
both $u$ and $v$ are H\"older continuous. Moreover, using
Morrey-Sobolev's embedding theorem, and the fact that $u\equiv v$ on
$\partial B_R$, we estimate \begin{eqnarray*} \left| \int_{B_R}
(u-v) \, d\mu \right| & \leq & ||u-v||_{L^\infty(B_R)} |\mu|(B_R)
\\ &\leq & cR^{\frac{p-n}{p}} \|Du-Dv\|_{L^p(B_R)} |\mu|(B_R)\;.
\end{eqnarray*} Combining the last inequality with \rif{weak2} and
using \rif{mon3}-\rif{mon2} we gain
$$
\|V(Du)-V(Dv)\|_{L^2(B_R)}^2+\|Du-Dv\|_{L^p(B_R)}^p\leq
cR^{\frac{p-n}{p}} \|Du-Dv\|_{L^p(B_R)} |\mu|(B_R),
$$
thereby, applying Young's inequality and then using \rif{morrey} we
conclude
\begin{eqnarray*}
&& \|V(Du)-V(Dv)\|_{L^2(B_R)}^2+\|Du-Dv\|_{L^p(B_R)}^p\\ &&
\hspace{25mm}\leq cR^{\frac{p-n}{p-1}}
 [|\mu|(B_R)]^{\frac{p}{p-1}}\leq  cR^{\frac{p-\theta}{p-1}}M^{\frac{1}{p-1}}
 |\mu|(B_R)\;.
\end{eqnarray*}
\end{proof}
\begin{lemma} \label{coco52} Under the assumptions \trif{asp} and \trif{morrey} with $0\leq \theta < p<n$,
and with $u, v \in W^{1,p}(B_R)$ as in \trif{Dir1} and \trif{Dirc1}
respectively, we have that $\trif{baba}$ holds, with $\sigma(p)$ as
in \trif{fra123} and $c\equiv c(n,p,\nu,\theta)$. In the case $0\leq
\theta<p=n$ estimate \trif{baba} remains valid modulo replacing the
right hand side by $cM^{1/(n-1)}|\mu|(B_R)R^{\sigma}$, for any
choice $\sigma < \sigma(n)=(n-\theta)/(n-1)$, where $c\equiv
c(n,\nu,\theta,\sigma)$.
\end{lemma}
\begin{proof} Firstly we deal with the case $p<n$. We test \rif{weak2} with $\varphi\equiv u-v$, which is
again admissible since $\theta <p$ implies that $\mu \in
W^{-1,p'}(\Omega)$. Therefore using again monotonicity \rif{mon2} as
for the previous lemma, using H\"older's inequality and applying
Theorem \ref{adamo} with the measure $\lambda \equiv |\mu|$, we have
\begin{eqnarray}
\nonumber&&\|V(Du)-V(Dv)\|_{L^2(B_R)}^2+\|Du-Dv\|_{L^p(B_R)}^p \leq
c \left| \int_{B_R} (u-v) \, d \mu\right|\\ & &\nonumber
\hspace{25mm}\leq c [|\mu|(B_R)]^{1-\frac{1}{p}}\left(\int_{B_R}
|u-v|^p \, d |\mu|\right)^{\frac{1}{p}} \\ & &\hspace{25mm}\leq
cM^{\frac{1}{p}}R^{\frac{p-\theta}{p}}[|\mu|(B_R)]^{1-\frac{1}{p}}\left(\int_{B_R}
|Du-Dv|^p \, d x\right)^{\frac{1}{p}}\label{adamin}\;.\end{eqnarray}
Using again Young's inequality yields \rif{baba}. In order to treat
the case $p=n$ it suffices to use Theorem \ref{adamo} again, and
applying it in \rif{adamin} as for the case $p<n$.
\end{proof} Finally, a comparison lemma of a different type.
With $B_{\bar{R}}\equiv B(x_0,{\bar{R}}) \subseteq B(x_0,R)$ and $v$
as in \rif{Dirc1}, let us define $v_0 \in v +
W^{1,p}_0(B_{\bar{R}})$ as the unique weak solution to \eqn{Dirc2}
$$
\left\{
    \begin{array}{cc}
    -\divo \ a(x_0,Dv_0)=0 & \qquad \mbox{in} \ B_{\bar{R}}\\
        v_0= v&\qquad \mbox{on $\partial B_{\bar{R}}$.}
\end{array}\right.
$$
\begin{lemma} \label{coco6} Under the assumptions \trif{asp} and \trif{lipi},
with $v$ as in \trif{Dirc1} and $v_0$ as in \trif{Dirc2}, with
$c\equiv c(n,p,\ratio)$ we have
 \eqn{baba3}
$$
\int_{B_{\bar{R}}} |V(Dv_0)-V(Dv)|^2 \, dx\leq c
\bar{R}^{2}\int_{B_{\bar{R}}}(|Dv|^p+s^p)\, dx \;.
$$
\end{lemma}\begin{proof} Using \rif{veramon}-\rif{veracre} it
follows that $v_0$ is a $Q$-minimum of the functional $w \mapsto
\int_{B_R}(|Dw|^p+s^p )dx$, with $Q \equiv Q(n,p,\ratio)$, see
Theorem 6.1 from \cite{G} that implies \eqn{ellsta}
$$
\int_{B_{\bar{R}}} |Dv_0|^p\, dx \leq
c(n,p,\ratio)\int_{B_{\bar{R}}} (|Dv|^p+s^p)\, dx\;.
$$
In turn, using \rif{asp}$_1$, and the fact that both $v$ and $v_0$
are solutions, we have
\begin{eqnarray*}
&& \int_{B_{\bar{R}}}
(s^2+|Dv_0|^2+|Dv|^2)^{\frac{p-2}{2}}|Dv-Dv_0|^2\,
dx\\
& &\qquad \leq c\int_{B_{\bar{R}}} \langle a(x_0,Dv)-a(x_0,Dv_0),
Dv-Dv_0\rangle \, dx\\& &\qquad =c\int_{B_{\bar{R}}} \langle
a(x_0,Dv)-a(x,Dv), Dv-Dv_0\rangle \, dx\\& &\qquad
\stackrel{\rif{lipi}}{\leq} c
\bar{R}\int_{B_{\bar{R}}}(s^2+|Dv_0|^2+|Dv|^2)^\frac{p-1}{2}|Dv-Dv_0|\,
dx \\ && \qquad \stackrel{\textnormal{Young}}{\leq}
\frac{1}{2}\int_{B_{\bar{R}}}
(s^2+|Dv_0|^2+|Dv|^2)^{\frac{p-2}{2}}|Dv-Dv_0|^2\, dx\\
&& \qquad \qquad \qquad + c\bar{R}^{2}\int_{B_{\bar{R}}}
(s^2+|Dv_0|^2+|Dv|^2)^{\frac{p}{2}}\, dx\;.\end{eqnarray*} Using
 \rif{V} for the left hand side, we get
$$
\int_{B_{\bar{R}}} |V(Dv_0)-V(Dv)|^2 \, dx\leq c
\bar{R}^{2}\int_{B_{\bar{R}}}(|Dv|^p+|Dv_0|^p+s^p)\, dx \;,
$$
and \rif{baba3} follows by merging the latter inequality with
\rif{ellsta}.
\end{proof}
\begin{remark}\label{scalabene}(Global estimates by scaling).
 We
consider \rif{Dirapp2}, and we find a global a priori estimate for
$u$, making explicit the ones in \cite{BG1}. Let us go back to Lemma
\ref{coco1}, Step 2, and let's modify a bit the estimates given
there. Assume not only that $\|f\|_{L^1(\Omega)}\leq 1$ as in
\rif{unob}, but now also that $s^q|\Omega|\leq 1$. Just forget about
$v$, and test \rif{Dirapp2}$_1$ on the whole $\Omega$ with $T_k(u)$
and $\Phi_k(u)$, as $u \in W^{1,p}_0(\Omega)$; then use
\rif{veramon} to get \eqn{apglo1}
$$\left\{\begin{array}{c}\displaystyle \int_{D_k} |Du|^p\, dx \leq
ck\|f\|_{L^1(\Omega)}+cs^p|D_k|\, \\ \\ \displaystyle
 \int_{C_k} |Du|^p\, dx \leq c\|f\|_{L^1(\Omega)}+cs^p|C_k|\;,\end{array}\right.$$
 where this time $D_k:=\{x \in \Omega  :  |u(x)|\leq k\}$ and
$C_k:=\{x \in \Omega  : k \ < |u(x)| \leq k+1\}$, consistently with
\rif{dkappa} and \rif{ckappa} respectively, and the constant $c$
just depends on $n,p,\ratio$. Then proceed as in \rif{gal}, but
using \rif{apglo1}, and we get $$ \|Du\|_{L^q(C_k)}^q \leq
c|C_k|^{1-q/p}+cs^q|C_k|\;.$$ Summing up these inequalities as for
\rif{semi}, the terms $|C_k|^{1-q/p}$ are treated as in \rif{gal}
and subsequent estimates, while, obviously,
$s^q|D_{k_0}|+\sum_{k\geq k_0} s^q|C_k|= s^q|\Omega|$. Therefore,
when $p<n$, it follows that $\|Du\|_{L^q(\Omega)}^q\leq
c(1+s^q|\Omega|)\leq \tilde{c}$, where $\tilde{c}$ is universal in
the sense it only depends on $n,p,\ratio,q$, and on $\Omega$. In the
case $p=n$, which already requires a different treatment in Lemma
\ref{coco1}, $\tilde{c}$ must be replaced by
$\tilde{c}(|\Omega|^{1/q-1/p}+1)$; indeed we need to use also
$$ \|Du\|_{L^q(D_{k_0})}^q \leq
c(k_0^{q/p}\|f\|_{L^1(\Omega)}^{q/p}|\Omega|^{1-q/p}+s^q|\Omega|)\leq
\tilde{c}(|\Omega|^{1-q/p}+1)\;,$$ that comes from \rif{apglo1}
exactly as in Lemma \ref{coco1}, Step 2, case $p=n$. Now we use a
scaling argument to treat the general case. Define, $A :=
\|f\|_{L^1(\Omega)}^{1/(p-1)}+s|\Omega|^{1/q}>0$, and accordingly,
$\tilde{u}:=A^{-1}u$, $\bar{f}:=A^{1-p}f$,
$\tilde{a}(x,z):=A^{1-p}a(x,Az)$, so that the vector field
$\tilde{a}$ satisfies \rif{asp} with $s$ replaced by $s/A$. Moreover
$\tilde{u}$ satisfies div $\tilde{a}(x,\tilde{u})=\bar{f}$ in the
weak sense and obviously $\tilde{u} \in W^{1,p}_0(\Omega)$. By the
definition of $A$ we have that $\|\tilde{f}\|_{L^1(\Omega)}\leq 1$
and $(s/A)^q|\Omega|\leq 1$, therefore we get the universal bounds
$\|D\tilde{u}\|_{L^q(\Omega)}\leq \tilde{c}$ when $p<n$, and
$\|D\tilde{u}\|_{L^q(\Omega)}\leq \tilde{c}(|\Omega|^{1/q-1/n}+1)$
when $p=n$. Taking into account the definitions of $\tilde{a}$ and
$A$ the latter inequalities readily give \eqn{apscala}
$$
\| Du\|_{L^{q}(\Omega)} \leq c\| f \|_{L^1(\Omega)}^{1/(p-1)}+c
s|\Omega|^{1/q}\;,
$$
that is the estimate we were looking for; the constant $c$ in
\rif{apscala} will depend on $n,p,q,\ratio$, and $\Omega$. The
dependence on $\Omega$ is on $(|\Omega|^{1/q-1/n}+1)$ in the case
$p=n$. As for Step 2 from Lemma \ref{coco1}, here everything works
for $q \in [1,p-1)$ too.
\end{remark}
\begin{remark} An a priori estimate can be derived for the super-capacitary case of
Theorem \ref{main4} too. Testing \rif{desol} with $u$, and this is
possible since $\theta <p$ implies $\mu \in W^{-1,p'}$, using
\rif{veramon}, and proceeding as in Lemmata \ref{coco5} and
\ref{coco52}, but using \rif{ada2} instead of \rif{ada}, we have,
with $c$ as in Theorem \ref{main4} \eqn{apestsub2222}
$$
\|Du\|_{L^{p}(\Omega)}^p\leq cM^{1/(p-1)}|\mu|(\Omega)
+cs^p|\Omega|\;.
$$
\end{remark}
\section{Basic approximation}
In order to establish the existence and regularity results for the
problem \rif{Dir1} for a general measure $\mu$, a standard device
\cite{BG1,DHM1} is to consider solutions to suitable approximate
problems, and then to prove a priori estimates; the main feature of
such solutions is to be in the natural space $W^{1,p}_0(\Omega)$.
Then the final assertion follows by a suitable passage-to-the-limit
argument. We remark that this procedure is not necessary when
considering the assumptions of Theorem \ref{main4}, that is when
$\theta < p$ ($p$-capacitary measures). In this section we set up
the approximation scheme, considered in the rest of the paper for
the case $\theta \in [p,n]$. In fact, as already remarked in the
Introduction, in the case $\theta < p$ the measure $\mu$ belongs to
the dual space $W^{-1,p'}$ (see \cite{adams}), and at this point the
standard monotone operator theory provides a unique solution to
\rif{Dir1} in the natural energy space $W^{1,p}_0(\Omega)$,
therefore no approximation with $W^{1,p}$-solutions is obviously
needed.

We consider a standard, symmetric and non-negative mollifier $\phi
\in C^{\infty}_0(B_1)$ such that $\|\phi\|_{L^1(\er^n)}=1$, and then
define, for $k \in \en$, $\phi_k(x):=k^n\phi(kx)$. Finally the
functions $f_k:\er^n \to \er$ are defined via convolution, $
f_k(x):=(\mu*\phi_k)(x).$ Since in particular $f_k \in
L^\infty(\Omega)$, applying standard monotonicity methods
\cite{lions} we can find a unique $u_k \in W^{1,p}_0(\Omega)$ such
that \eqn{Dirapp}
$$
\left\{
    \begin{array}{cc}
    -\divo \ a(x,Du_k)=f_k & \qquad \mbox{in $\Omega$}\\
        u_k= 0&\qquad \mbox{on $\partial\Omega$.}
\end{array}\right.
$$
From now on and for the rest of the paper the sequence
$\{u_k\}_k\subset W^{1,p}_0(\Omega)$ will be the one fixed by
\rif{Dirapp}. Let us collect some basic facts now. Up to extracting
a non-relabeled subsequence we can assume \eqn{convergenzabase}
$$ f_k \rightharpoonup \mu
\hspace{1cm} \mbox{weakly in the sense of measures}.$$ Moreover,
looking at \cite{my}, Proposition 2.7, we have\eqn{total}
$$ \|f_k\|_{L^1(\Omega)}\leq |\mu|(\Omega), \qquad \|f_k\|_{L^{1,\theta}(\Omega)}
\leq \|\mu\|_{L^{1,\theta}(\Omega)}\stackrel{\rif{morrey}}{\leq} M\;
,$$ and \eqn{total2}
$$
 \|f_k\|_{L^1(B_R)}\leq |\mu|(B_{R+1/k}), \qquad \|f_k\|_{L^{1,\theta}(B_R)}\leq \|\mu\|_{L^{1,\theta}(B_{R+1/k})}\;.
$$
Applying Remark \ref{scalabene} and in particular estimate
\rif{apscala} to $u_k$, and eventually using \rif{elemV}, we get
\eqn{apbase}
$$
\int_{\Omega} \! |V(Du_k)|^{2q/p}\, dx + \int_{\Omega} \!
|Du_k|^{q}\, dx \leq c[|\mu|(\Omega)]^{\frac{q}{p-1}}+c
s^q|\Omega|,\qquad  \forall \ q < b\;,
$$
where $c \equiv c(n,p,\ratio,q,\Omega)$, and is independent of $k
\in \en$. Therefore we immediately obtain that up to a non-relabeled
subsequence \eqn{convergenzeu}
$$u_k \rightharpoonup u \ \mbox{ weakly in} \
W^{1,q}(\Omega)\qquad \mbox{ and } \qquad u_k \to u\ \mbox{ strongly
in} \ L^{q}(\Omega)\;.$$ The function $u$ is eventually shown to be
a solution of \rif{Dir1}. The proof of this fact usually involves
certain truncation arguments \cite{BG1, BG2} to prove the strong
convergence of the gradients. Thanks to the stronger a priori
estimates we derive here we shall give a very short proof of such
convergence; see Theorem \ref{main1}. The function $u$ is the
solution the Theorems and results of Section 1 in the
super-capacitary case $\theta \geq p$.
\section{General measures}
This section is mainly devoted to the proof of Theorems
\ref{main1}-\ref{noncz} and Corollary \ref{nondegenere}. The
ingredients will be: the lemmata of Section 4, the key estimate
below the growth exponent \rif{revca}, and a variant of a fractional
regularity technique recently introduced in \cite{KM} in order to
obtain singular sets estimates for variational problems.

{\bf Warning.} In the rest of the paper we shall very often deal
with a solution $u$ to problem \rif{Dirapp2}, for a fixed, but a
priori un-specified $L^{\infty}$ function $f$. Eventually we shall
take $f \equiv f_k$ and $u \equiv u_k $, where $u_k,f_k$ appear in
\rif{Dirapp}.

Keeping \rif{ssiidef} in mind, let us define \eqn{ildelta0}
$$
\delta:= \frac{p\sigma(q,\theta)}{2q}\;;\qquad
\gamma(t):=\frac{\delta}{\delta+1-t},\qquad \mbox{for every} \ t \in
[0,\delta+1)\;.
$$
\begin{remark}\label{deltauno}We have $\delta \leq 1$. Indeed, by
\rif{ssiidef}, when $p-1\leq q$ then $\sigma(q,\theta)\leq 1$,
therefore
$$
\delta \leq \frac{p}{2q}\leq \frac{p}{2(p-1)}\leq 1\;,
$$
which holds since $p\geq 2$. For the same reason we have $2q/p \geq
1$ for $q \geq p-1$.
\end{remark}
We shall start deriving a priori estimates for $W^{1,p}$-solutions
to \rif{Dirapp2}. We set \eqn{tilf}
$$ \bar{f}:=\|f\|_{L^{1,\theta}(\Omega'')}^{\frac{q-p+1}{p-1}}f \qquad \mbox{if} \qquad p-1\leq q
\;,
$$
where $\Omega'' \subseteq \Omega$ will be clarified in Lemmas
\ref{2it} and \ref{222it} below. When $\theta=n$ it follows directly
from definition \rif{dede1mo} that \eqn{triviale}
$$\|\bar{f}\|_{L^1(\Omega)}\leq
\|f\|_{L^1(\Omega)}^{\frac{q}{p-1}}\;.$$
\begin{lemma}\label{2it}
Let $u \in \Wp_0(\Omega)$ be the unique solution to \trif{Dirapp2},
under the assumptions \trif{asp} and \trif{lipi} with $p \leq n$,
and let $q$ be such that $p-1 \leq q < b.$ Assume that \eqn{initial}
$$ V(Du) \in
W^{t,2q/p}_{\rm loc} (\Omega ,\er^n)\,, \qquad \mbox{for some}\ \ t
\in [0,\delta)\;,$$ where $\delta$ is as in \trif{ildelta0}, and that
for every couple of open subsets $\Omega^{\prime} \CC \Omega^{\prime
\prime} \CC \Omega$ there exists $c_1\equiv
c_1(\dist(\Omega',\partial\Omega''))$, such that \eqn{eq1it}
$$
[V(Du)]_{t ,2q/p;\Omega^{\prime}}^{2q/p} \leq c_1
\int_{\Omega^{\prime \prime}} \! ( |Du|^{q}+s^q+|\bar{f}| ) \, dx\;.
$$
Then \eqn{final}
$$ V(Du) \in W^{\tilde{t}, 2q/p}_{\rm loc} (\Omega
,\er^n )\,,\qquad \mbox{for every}\  \ \tilde{t} \in
[0,\gamma(t))\;,$$ where $\gamma(\cdot)$ is in \trif{ildelta0}, and
for every couple of open subsets $\Omega^{\prime} \CC \Omega^{\prime
\prime} \CC \Omega$ there exists a new constant $c$ depending only
on $n,p,\ratio,q, \distc,\tilde{t},c_1$, such that \eqn{bettere}
$$
[V(Du)]_{\tilde{t},2q/p;\Omega^{\prime}}^{2q/p} \leq c
\int_{\Omega^{\prime \prime}} \! ( |Du|^{q}+s^q+|\bar{f}| ) \, dx\;.
$$
Moreover, for every $i \in \{1,\ldots, n\}$ and with $0 < |h| <
\dist(\Omega',\partial\Omega'')$ \eqn{bettereh1}
$$
\sup_{h} \int_{\Omega'}
\frac{|\tau_{i,h}V(Du(x))|^{2q/p}}{|h|^{\gamma(t)2q/p}}\, dx\leq c
\int_{\Omega^{\prime \prime}} \! ( |Du|^{q}+s^q+|\bar{f}| ) \, dx\;.
$$
\end{lemma}
\begin{proof} We fix a notation that we shall keep for the rest of the
paper. Let us take $B\CC \Omega$, a ball of radius $R$; we shall
denote by $\Qin \equiv \Qin(B)$ and $\Qou \equiv \Qou(B)$ the
largest and the smallest cubes, concentric to $B$ and with sides
parallel to the coordinate axes, contained in $B$ and containing
$B$, respectively; clearly $|B|\approx |\Qin|\approx|\Qou|\approx
R^n$. The cubes $\Qin(B)$ and $\Qou(B)$ will be called the {\em
inner} and the {\em outer} cubes of $B$, respectively. We also
denote the enlarged ball as $\hat{B} \equiv 16B$. Consistently with
such a notation we put $ \Qin \equiv \Qin(B)$ and
$\hat{Q}_{\textnormal{out}} \equiv \Qou(\hat{B})$, and therefore we
have the following chain of inclusions: \eqn{chain}
$$ \Qin \subset B \CC 2B \CC 4B \subset \Qin(\hat{B}) \subset \hat{B}
\subset \hat{Q}_{\textnormal{out}}\;.$$ Now we fix arbitrary open
subsets $\Omega^{\prime} \CC \ompp \CC \Omega$, and then take $\beta
\in (0,1)$ to be chosen later, and let $h \in \er$ be a real number
satisfying \eqn{ha}
$$
0<|h|\leq \min \left\{\left(
\frac{\textnormal{dist}(\Omega^{\prime},
\partial \Omega^{\prime \prime})}{10000\sqrt{n}}\right)^{\frac{1}{\beta}},
\left(\frac{1}{10000}\right)^{\frac{1}{1-\beta}}\right\}=:d<\textnormal{dist}(\Omega^{\prime},
\partial \Omega^{\prime \prime})\;.
$$
We take $x_0 \in \omp$, and fix a ball of radius $|h|^{\beta}$
\eqn{pallah}
$$B
\equiv B(h) = B(x_{0},|h|^{\beta})\;.$$ By \rif{ha} we have
$\hat{Q}_{\textnormal{out}} \subset \Omega^{\prime \prime}$. Let us
first define $v \in u+W^{1,p}_{0}(\hat B)$, and then $v_0 \in
v+W^{1,p}_{0}(8 B)$, as the unique solutions to the following
Dirichlet problems: \eqn{cdir2}
$$
\left\{
    \begin{array}{cc}
    -\divo \ a(x,Dv)=0 & \ \  \mbox{in} \ \hat B\\
        v= u&\ \  \mbox{on $\partial \hat B$,}
\end{array}\right.
$$
and
\eqn{cdir2bis}
$$ \left\{
    \begin{array}{cc}
    -\divo \ a(x_0,Dv_0)=0 & \ \ \mbox{in} \ 8 B\\
        v_0= v&\ \  \mbox{on $\partial 8B$,}
\end{array}\right.
$$
respectively. Now we fix $i \in \{1,\ldots, n\}$ and write, using
that $|h|\leq d$ from \rif{ha}
\begin{eqnarray}
\nonumber \int_{B} \! |\tau_{i,h}V(Du)|^{2q/p} \, dx &\leq&
c\int_{B} \! |\tau_{i,h}V(Dv_0)|^{2q/p} \, dx \\
 && +c\int_{B} \! |V(Du(x+he_{i}))-V(Dv(x+he_{i}))|^{2q/p} \, dx \nonumber\\
&& +c\int_{B} \! |V(Dv(x+he_{i}))-V(Dv_0(x+he_{i}))|^{2q/p} \, dx \nonumber\\
&& \nonumber +c\int_{B} \! |V(Du)-V(Dv)|^{2q/p} \, dx
\\&& +c\int_{B} \! |V(Dv)-V(Dv_0)|^{2q/p} \, dx \nonumber \\ &\leq&
\hspace{-1mm}c\int_{B} \! |\tau_{i,h}V(Dv_0)|^{2q/p} \, dx
+c\int_{\hat B} \! |V(Du)-V(Dv)|^{2q/p}\, dx \nonumber\\ &&
\nonumber \hspace{3cm}+c\int_{ 2B} \! |V(Dv)-V(Dv_0)|^{2q/p} \,
dx\nonumber
\\& =: & I + II + III\;. \label{lage}
\end{eqnarray}
In order to estimate $II$ we shall use Lemmas \ref{coco1} and
\ref{coco2}, this last one when $q \geq p-1$ and $\theta<n$; by the
definition of $\sigma(q,\theta)$ in \rif{ssiidef}, we have
\eqn{comp111}
$$ \int_{\hat B}
|V(Du)-V(Dv)|^{2q/p} \, dx \leq c  \left(\int_{\hat B} |\bar{f}|\,
dx \right) \, |h|^{\beta \sigma(q,\theta)} \;,
$$
where we used \rif{tilf} too. To estimate $III$ we first appeal to
Lemma \ref{coco6} that gives \eqn{dopodopo}
$$
\int_{8B} |V(Dv_0)-V(Dv)|^2 \, dx\leq c
\left(\int_{8B}(s^2+|Dv|^2)^\frac{p}{2}\, dx\right) |h|^{\beta 2}\;,
$$
and then apply Lemma \ref{superrev} to $v$ in \rif{cdir2}; with
$\chi \equiv \chi (n,p,\ratio)>1$ being the exponent determined in
in Lemma \ref{superrev} we have
\begin{eqnarray*}
\int_{8B} |V(Dv_0)-V(Dv)|^{2q/p}\, dx & \leq & c
|h|^{n\beta\left(1-\frac{q}{p} \right) } \left(\int_{8B}
|V(Dv_0)-V(Dv)|^2 \, dx \right)^{\frac{q}{p}}\\
& \stackrel{\rif{dopodopo}}{\leq} & c |h|^{\frac{\beta 2q }{p} +n
\beta } \left(\mean{8B}
(s^2+|Dv|^2)^\frac{p}{2} \, dx \right)^{\frac{q}{p}} \\
& \leq & c |h|^{\frac{\beta 2q }{p} +n \beta } \left(\mean{8B}
(s^2+|Dv|^2)^\frac{p\chi}{2} \, dx \right)^{\frac{q}{p\chi}} \\
& \stackrel{\rif{zz2}}{\leq} & c |h|^{\frac{\beta 2q }{p}}
\int_{\hat B} (s^2+|Dv|^2)^\frac{q}{2} \, dx \\ &
\stackrel{\rif{elemV}}{\leq} & c |h|^{\frac{\beta 2q}{p}} \int_{\hat
B} (s^2+|V(Dv)|^{\frac{4}{p}})^{\frac{q}{2}} \, dx\\&
\stackrel{\rif{elemV}}{\leq} & c |h|^{\frac{\beta 2q}{p}} \int_{\hat
B} (s^q+|Du|^q\\&& \qquad \qquad  + |V(Du)-V(Dv)|^{2q/p})\,
dx\\
& \stackrel{\rif{comp111}}{\leq} & c |h|^{\frac{\beta 2q }{p}}
\int_{\hat B} (s^q+|Du|^q+ |\bar{f}|) \, dx\;.
\end{eqnarray*}
We recall that $16 B = \hat B$. Summarizing the latter estimate and
\rif{comp111} yields
$$
II + III \leq c\left[  |h|^{\beta \sigma(q,\theta)}+|h|^{\beta
2q/p}\right] \int_{\hat B} (s^q+|Du|^q+ |\bar{f}|) \, dx\;,
$$
where $c \equiv c(n,p,\ratio,q)$ is independent of any of the balls
considered. Recalling \rif{ildelta0} and Remark \ref{deltauno} that
gives $\delta \leq 1$, we estimate $|h|^{\beta 2q/p}\leq |h|^{\beta
\sigma(q,\theta)}=|h|^{\beta \delta 2q/p}$ as $|h|\leq 1$, therefore
\eqn{IIIII}
$$ II + III \leq c |h|^{\beta \delta 2q/p} \int_{\hat B} (s^q+|Du|^q+
|\bar{f}|) \, dx\;.
$$
Implicit in the previous inequality is\eqn{late}
$$ \int_{8B} |V(Du)-V(Dv_0)|^{2q/p} \, dx \leq c |h|^{\beta \delta 2 q/p} \int_{\hat B} (s^q+|Du|^q+
|\bar{f}|) \, dx\;.
$$
Now we turn to $I$. Applying Lemma \ref{revbase} to $v_0$ taking
$a_0(z)\equiv a(x_0,z)$, \rif{cacc2} gives  \eqn{cacc}
$$
\aint_{2B} \! |D ( V(Dv_0) ) |^{2} \, dx \leq
c|h|^{-2\beta}\aint_{4B} \! |V(Dv_0)-V(z_{0})|^{2} \, dx,
$$
for every $z_{0} \in \er^n$, while using \rif{revca} with $t= q/p$,
we also have \eqn{revcac}
$$
\left(\mean{4B} \! |V(Dv_0)-V(z_{0})|^{2} \,
dx\right)^{\frac{q}{p}}\leq c\mean{8 B} \! |V(Dv_0)-V(z_{0})|^{2q/p}
\, dx\;.
$$
Now, again using H\"older's inequality yields \eqn{Iuno}
$$ \int_{B}
\! |\tau_{i,h}V(Dv_0)|^{2q/p} \, dx\leq  c
|h|^{n\beta\left(1-q/p\right)} \left(\int_{B} \!
|\tau_{i,h}V(Dv_0)|^{2} \, dx \right)^{\frac{q}{p}}\;.
$$
Using the definition of the operator $\tau_{i,h}$ in \rif{detau},
elementary properties of Sobolev functions, and again the
restriction on $|h|$ imposed in \rif{ha} that in this case serves to
ensure that $B(x_0, |h|^\beta) + B(0,|h|)\subset B(x_0,
2|h|^\beta)$, we have
\begin{eqnarray}
\int_{B} \! |\tau_{i,h}V(Dv_0)|^{2} \, dx  &\leq & c|h|^2 \nonumber
\int_{2B} \! |DV(Dv_0)|^{2} \, dx\\
&\stackrel{\rif{cacc}}{\leq} & c|h|^{2-2\beta} \int_{4B} \!
|V(Dv_0)-V(z_0)|^{2} \, dx\;.\label{Idue}
\end{eqnarray}
Combining \rif{Iuno} and \rif{Idue} gives
$$
\int_{B} \! |\tau_{i,h}V(Dv_0)|^{2q/p} \, dx \leq c
|h|^{(1-\beta)2q/p+n\beta \left(1-q/p\right)} \left(\int_{4B} \!
|V(Dv_0)-V(z_0)|^{2} \, dx\right)^{\frac{q}{p}}\;.
$$
Using now \rif{revcac} gives with $c\equiv c(n,p,\ratio,q)$
\eqn{cala}
$$
I = c \int_{B} \! |\tau_{i,h}V(Dv_0)|^{2q/p} \, dx \leq c
|h|^{(1-\beta)2q/p} \int_{8B} \! |V(Dv_0)-V(z_0)|^{2q/p} \, dx\;,
$$
and we estimate the last integral; recall that in the latter
estimate $z_0 \in \er^n$ is still to be chosen. We shall distinguish
two cases now.

{\em Case } $t=0$. In this case we take $z_0=0$ in \rif{cala}; then
\rif{late} and \rif{elemV} yield
\begin{eqnarray} && \int_{8B} \! |V(Dv_0)-V(z_0)|^{2q/p} \,
dx \leq c \int_{8B}
(s^q+|Dv_0|^q)  \, dx\nonumber \\
& & \qquad \leq c \int_{8B} (s^q+|Du|^q)  \, dx +c\int_{8B}
|V(Du)-V(Dv_0)|^{2q/p}  \, dx\nonumber \\
&&  \qquad \qquad \qquad \leq c \int_{\hat B} (s^q+|Du|^q+|\bar{f}|)
\, dx \;.\label{coplno}
\end{eqnarray}

{\em Case } $t>0$. In this case we choose $z_0$ as the following
``average":\eqn{copl2}
$$
z_0:= V^{-1}\left((V(Du))_{8B}\right)\;;
$$
observe that such a choice is possible since the map $V$ is
bijective. Now, first
\begin{eqnarray} \int_{8B} \! |V(Dv_0)-V(z_0)|^{2q/p} \,
dx &\leq & c \int_{8B} \! |V(Dv_0)-V(Du)|^{2q/p} \,
dx\nonumber \\
&& \qquad +c \int_{8B} \! |V(Du)-V(z_0)|^{2q/p}\, dx\;.\label{prim2}
\end{eqnarray}
Then by \rif{initial} and Proposition \ref{frapoes} with
\rif{copl2}, we have \eqn{copl3}
$$ \int_{8B} \!
|V(Du)-V(z_0)|^{2q/p} \, dx  \leq  c|h|^{\beta t2q/p
}[V(Du)]_{t,2q/p;8B}^{2q/p} \;.
$$
Combining \rif{copl3} and \rif{late} with \rif{prim2} we have
\begin{eqnarray}
&& \nonumber \int_{8B} \! |V(Dv_0)-V(z_0)|^{2q/p} \, dx\\ && \qquad
\qquad \leq c |h|^{ \beta t2q/p }\left\{\int_{\hat B} (s^q+|Du|^q+
|\bar{f}|) \, dx + [V(Du)]_{t,2q/p;\hat
B}^{2q/p}\right\}.\label{copl10}
\end{eqnarray}
Observe that we have used $ t < \delta $ to estimate $|h|^{\beta
\delta 2q/ p}\leq |h|^{\beta t2q /p}$ as $|h|\leq 1$.

Now let us define for any measurable set $A \CC \Omega$ the
following set function: \eqn{setta}
$$
\lambda(A):=\int_{A} (s^q+|Du|^q+ |\bar{f}|) \, dx +
\chi(t)[V(Du)]_{t,2q/p;A}^{2q/p}\;,$$ where $\chi(t)=0$ if $t=0$,
and $\chi(t)=1$ if $t>0$. Summarizing \rif{cala}, \rif{coplno} and
\rif{copl10} we have
$$
I = c \int_{B} \! |\tau_{i,h}V(Dv_0)|^{2q/p} \, dx \leq c
|h|^{[(1-\beta)+ t\beta]2q/p} \lambda(\hat B)\;.
$$
Combining the latter estimate with \rif{IIIII}, and in turn with
\rif{lage}, we find
$$
\int_{B} \! |\tau_{i,h}V(Du)|^{2q/p} \, dx \leq c
\left[|h|^{[(1-\beta)+ t\beta]2q/p}+ |h|^{\beta \delta2q/p}
\right]\lambda(\hat B)\;.
$$
Since by \rif{chain} $\Qin(B)\equiv \Qin\subset B$ and $\hat B
\subset \hat{Q}_{\textnormal{out}} \equiv \Qou(\hat B)$, we finally
obtain \eqn{dasomma}
$$
\int_{\Qin} \! |\tau_{i,h}V(Du)|^{2q/p} \, dx \leq \tilde{c}
\left[|h|^{[(1-\beta)+ t\beta]2q/p}+ |h|^{\beta \delta2q/p}
\right]\lambda(\hat{Q}_{\textnormal{out}})\;.
$$
Now we conclude with a covering argument.~Preliminary, observe that
the set function $\lambda(\cdot)$ in \rif{setta} is not a measure
due to the presence of $[V(Du)]_{t,2q/p;A}$ in its definition, but
it is nevertheless countably super-additive, that is \eqn{cou}
$$
\sum \lambda(A_j) \leq \lambda\left(\cup A_j\right)\;,
$$
whenever $\{A_j\}_j$ is a countable family of mutually disjoint
subsets. The covering argument goes now as follows: first recall
that all the cubes here have sides parallel to the coordinate axes;
then for each $h \in \er\setminus\{0\}$ satisfying \rif{ha} we can
find balls $B_1 \equiv B(x_1,|h|^\beta)$, ... , $B_J \equiv
B(x_J,|h|^\beta)$, $J \equiv J(h) \in \en$ of the type considered in
\rif{pallah} such that the corresponding inner cubes
$\Qin(B_1)$,\ldots, $\Qin(B_J)$ are disjoint and cover $\Omega'$ up
to a negligible set \eqn{ricoprimento}
$$|\Omega' \setminus \bigcup \Qin(B_j)|=0, \qquad \qquad \Qin(B_i) \cap \Qin(B_j) \not=\emptyset \Longleftrightarrow i =j\;.$$
Actually we are proceeding as follows: we first take a lattice of
cubes $\{Q_j\}$ with equal side length, comparable to $|h|^\beta$,
and sides parallel to the coordinate axes, in order to obtain
\rif{ricoprimento}. They must be centered in $\Omega'$. Then we view
them as the inner cubes of the balls $\{B(x_j, |h|^\beta)\}$,
according to \rif{chain}. Now we sum up inequalities \rif{dasomma}
for $j\leq J$ and get
\begin{eqnarray} \nonumber &&\sum \int_{\Qin(B_j)} \!
|\tau_{i,h}V(Du)|^{2q/p} \, dx \\ && \qquad \qquad \qquad \leq
\tilde{c} \left[|h|^{[(1-\beta)+ t\beta]2q/p}+ |h|^{\beta \delta
2q/p} \right]\sum\lambda(\Qou(\hat{B_j}))\;.
\label{somma}\end{eqnarray} By construction, and in particular by
\rif{ha}, we have $\Qou(\hat{B}_{j}) \subset \Omega^{\prime
\prime}$, for every $j \leq J$. Moreover by \rif{ricoprimento} each
of the dilated outer cubes $\Qou(\hat{B}_{j})$ intersects the
similar ones $\Qou(\hat{B}_{k})$ less than $(32\sqrt{n})^{n}$ times.
Using all these facts and \rif{cou}, in turns out that
\rif{ricoprimento}-\rif{somma} imply \eqn{quasi22}
$$
\int_{\Omega'} \! |\tau_{i,h}V(Du)|^{2q/p} \, dx \leq
2^{8n}\tilde{c} \left[|h|^{[(1-\beta)+ t\beta]2q/p}+ |h|^{\beta
\delta2q/p} \right]\lambda(\Omega'')\;.
$$
Now we determine $\beta$ in order to minimize the right-hand side
with respect to $|h|$; this yields $ [(1-\beta)+ t\beta]=\beta
\delta$, that is $ \beta=\gamma(t)/\delta$, see \rif{ildelta0}.
Observe that we are requiring everywhere that $\beta<1$, see
\rif{ha}, and the choice $\beta=\gamma(t)/\delta$ is admissible
since $t< \delta$ implies $\gamma(t)/\delta <1$. Accordingly, for
any $h$ as in \rif{ha}, \rif{quasi22} becomes \eqn{analog}
$$
\int_{\Omega'} \! |\tau_{i,h}V(Du)|^{2q/p} \, dx \leq c_0
|h|^{\gamma(t)2q/p}\lambda(\Omega'')\;,
$$
for $c_0 \equiv c_0(n,p,\ratio,q)$. Therefore, since $ i \in
\{1,\ldots,n\}$ is arbitrary, the crucial inequality \rif{crucial}
of Lemma \ref{diff2} is satisfied with $d$ as in \rif{ha}, $q$
replaced by $2q/p$, $\bar{\alpha}\equiv \gamma(t)$, and finally
$S\equiv [c_0\lambda(\Omega'')]^{p/2q}$. Up to changing the subsets
according to Lemma \ref{stanom}, that is passing to inner and outer
subsets to $\Omega''$ and $\Omega'$ respectively, we may apply Lemma
\ref{diff2} that now gives $V(Du) \in
W^{\tilde{t},2q/p}_{\loc}(\Omega',\er^n)$, for every $\tilde{t} <
\gamma(t)$; as $\Omega'$ is arbitrary, this proves the first part of
the assertion. Changing again the subsets, since $\Omega' \CC
\Omega''$ are themselves arbitrary, using estimate \rif{eess}, and
finally \rif{elemV}, we have that for every couple of open subsets
$\Omega' \CC \Omega''$ there exists a constant $c \equiv c
(n,p,\ratio,q, \dist(\Omega',\partial \Omega''))$ such that
\eqn{stimaccia}
$$
[V(Du)]_{\tilde{t},2q/p;\Omega^{\prime}}^{2q/p} \leq c
\int_{\Omega^{\prime \prime}} \! ( s^q+|Du|^{q} +|\bar{f}| ) \,
dx+c[V(Du)]_{t,2q/p;\Omega''}^{2q/p}\;.
$$
We have used \rif{elemV} to estimate the integral of $V$ arising
when applying \rif{eess}: \eqn{veratriv}
$$
\int_{\Omega^{\prime \prime}} \! |V(Du)|^{2q/p}\, dx\leq
c\int_{\Omega^{\prime \prime}} \! ( s^q+|Du|^{q} ) \, dx\;.
$$
Using \rif{stimaccia} in combination with \rif{eq1it}, and again
changing the subsets via Lemma \ref{stanom}, we finally obtain
\rif{bettere} with the specified dependence of $c$. In a completely
similar way using \rif{analog} it follows \rif{bettereh1} with
$|h|\leq d$ as in \rif{ha}. The full case $0 < |h| <
\dist(\Omega',\partial \Omega'')$ follows by increasing the constant
$c$ in \rif{analog} by a number depending on $n,p,q$ and $
\dist(\Omega',\partial \Omega'')$; indeed when $ |h| \in
(d,\dist(\Omega',
\partial \Omega''))$
\begin{eqnarray*}
&&\sup_{h} \int_{\Omega'}
\frac{|\tau_{i,h}V(Du(x))|^{2q/p}}{|h|^{\gamma(t)2q/p}}\, dx\\
&& \qquad \leq \frac{c}{d^{\gamma(t)2q/p}} \int_{\Omega'}
|V(Du(x+he_i))|^{2q/p}+|V(Du(x))|^{2q/p}\, dx \\ && \qquad \leq
\frac{c}{d^{\gamma(t)2q/p}} \int_{\Omega''} |V(Du)|^{2q/p}\,
dx\stackrel{\rif{elemV}}{\leq}\frac{c}{d^{\gamma(t)2q/p}}
\int_{\Omega''} (s^q+|Du|^q)\, dx \;.
\end{eqnarray*}
The proof, also of \rif{final}, is complete as the open subsets
considered are arbitrary.
\end{proof}
\begin{lemma}\label{222it}
Let $u \in \Wp_{0}(\Omega ,\er^N)$ be the unique solution to
\trif{Dirapp2}, under the assumptions \trif{asp} and \trif{lipi}
with $p\leq n$, and let $q$ be such that $p-1 \leq q < b.$ Then
\eqn{start22}
$$
V(Du) \in W^{t,2q/p}_{\rm loc} (\Omega ,\er^n), \quad Du \in
W^{2t/p,q}_{\rm loc} (\Omega ,\er^n), \ \mbox{for every} \ t \in
[0,\delta) \;,
$$
where $\delta$ is in \trif{ildelta0}. Moreover, for every couple of
open subsets $\Omega^{\prime} \CC \Omega'' \CC \Omega$ there exists
a constant $c \equiv c(n,p,\ratio,q,t, \dist(\Omega',
\partial \Omega''))$ such that \eqn{bettere2}
$$
[V(Du)]_{t,2q/p;\Omega^{\prime}}^{2q/p}+[ Du ]_{2t/p
,q;\Omega^{\prime}}^q \leq c \int_{\Omega''} \! ( |Du|^{q}+s^q
+|\bar{f}| ) \, dx\;
$$
and \eqn{bettereh2}
$$
\sup_{h} \int_{\Omega'} \frac{|\tau_{i,h}Du(x)|^{q}}{|h|^{t2q/p}}\,
dx\leq c \int_{\Omega^{\prime \prime}} \! ( |Du|^{q}+s^q+|\bar{f}| )
\, dx\;,
$$ for every $i \in \{1,\ldots, n\}$, where $0 < |h| <
\dist(\Omega',\partial\Omega'')$
\end{lemma}
\begin{proof} The proof follows from Lemma \ref{2it} via iteration.
 We first prove the assertion about $V(Du)$. The function
$\gamma(\cdot)$ in \trif{ildelta0} is seen to be increasing and it
satisfies
 \eqn{gapo1}
$$t \in (0,\delta) \Longrightarrow \gamma(t) \in (t,\delta)\qquad  \mbox{and }\qquad \gamma (\delta)= \delta\;.$$
Now, let us inductively define the two sequences $\{t_k\}_{k \geq
1}$ and $\{s_k\}_{k \geq 1}$ as \eqn{succ}
$$
s_1:=\frac{\delta}{4(\delta+1)},\quad t_1=2s_1, \quad
s_{k+1}:=\gamma(s_k), \quad
t_{k+1}:=\frac{\gamma(s_k)+\gamma(t_{k})}{2}\;.
$$
From \rif{gapo1} it follows that $s_k \nearrow \delta$, moreover,
since $\gamma(\cdot )$ is increasing we have that $s_k < t_k <
\delta$, so that also $t_k \nearrow \delta$ holds. We prove by
induction that $V(Du) \in W^{t_k,2q/p}_{\rm loc} (\Omega ,\er^n)$,
for every $k \in \en$; this will prove the first assertion in
\rif{start22}. Applying Lemma \ref{2it} with $t=0$ we immediately
get $V(Du) \in W^{t_1,2q/p}_{\rm loc} (\Omega ,\er^n)$, with a
corresponding estimate of the type \rif{eq1it}. Now assuming that
$V(Du) \in W^{t_k,2q/p}_{\rm loc} (\Omega ,\er^n)$, we may apply
again Lemma \ref{2it} with $t=t_k$, to get that $V(Du) \in
W^{t,2q/p}_{\rm loc} (\Omega ,\er^n)$ for every $t< \gamma(t_k)$.
Now observe, that since $\gamma(\cdot)$ is increasing and $s_k<t_k$,
we have that $t_{k+1}< \gamma(t_k)$, and therefore $V(Du) \in
W^{t_{k+1},2q/p}_{\rm loc} (\Omega ,\er^n)$, with corresponding
estimates of the type \rif{bettere} and \rif{bettereh1}. Taking into
account the fact that the open subsets $\Omega' \CC \Omega'' \CC
\Omega$ in Lemma \ref{2it} are arbitrary, and the estimates
\rif{eq1it} and \rif{bettere}, the part of \rif{bettere2} regarding
$V(Du)$ also follows by induction. In the same way, by induction on
\rif{bettereh1}, for every $i \in \{1,\ldots, n\}$ and considering
$0 < |h| < \dist(\Omega',\partial \Omega'')$, we have
\eqn{bettereh3}
$$
\sup_{h} \int_{\Omega'}
\frac{|\tau_{i,h}V(Du(x))|^{2q/p}}{|h|^{t2q/p}}\, dx\leq c
\int_{\Omega^{\prime \prime}} \! ( |Du|^{q}+s^q+|\bar{f}| ) \, dx,
\quad \forall \ t < \delta\;.
$$
The assertions concerning $Du$ instead follows using \rif{V} and the
fact that $p\geq 2$:
\begin{eqnarray}
\nonumber [ Du ]_{\frac{2t}{p} ,q;\Omega^{\prime}}^q&=&
\int_{\Omega^{\prime}} \int_{\Omega^{\prime}}
\frac{|Du(x)-Du(y)|^q}{|x-y|^{n+2tq/p}}\ dx  dy\\ \nonumber & \leq&
\int_{\Omega^{\prime}} \int_{\Omega^{\prime}}
\frac{\left[(s+|Du(x)|+|Du(y)|)^{p-2}|Du(x)-Du(y)|^2\right]^{q/p}}{|x-y|^{n+2tq/p}}\
dx  dy\\ \nonumber &\leq&  c\int_{\Omega^{\prime}}
\int_{\Omega^{\prime}}
\frac{|V(Du(x))-V(Du(y))|^{2q/p}}{|x-y|^{n+2tq/p}}\ dx  dy\\
&=& c[V(Du)]_{t,2q/p;\Omega^{\prime}}^{2q/p}\;,\label{converti}
\end{eqnarray}
for any $\Omega' \CC \Omega$, where $c \equiv c(n,p)$; this gives
\rif{bettere2}. A completely similar argument allows to get
\rif{bettereh2} from \rif{bettereh3}, and the proof is complete.
\end{proof}
\begin{proof}[Proof of Theorems \ref{main1} and \ref{main1es}] Firstly, observe that since $p \geq 2$, then $q\geq p-1$ implies
$2q/p\geq 1$, and therefore Lemma \ref{222it} can be used in the
full range \rif{ssii}. We consider the approximation sequence
$\{u_k\}_k$ built in Section 5. Applying to each $u_k$ the result of
Lemma \ref{222it}, and keeping in mind \rif{total}-\rif{apbase}, we
have \eqn{cisiamo1}
$$
\| Du_k \|_{{W^{\sigma/q,q}}(\Omega')}^q \leq c \int_{\Omega} \! (
s^q+|Du_k|^{q} +|\bar{f}_k| ) \, dx \leq
c[|\mu|(\Omega)]^{\frac{q}{p-1}} +cs^q |\Omega| \;,
$$
with the obvious definition of
$\bar{f}_k:=\|f_k\|_{L^1(\Omega)}^{\frac{q-p+1}{p-1}}f_k$: look at
\rif{tilf}-\rif{triviale} and recall that here it is $\theta =n$.
The constant $c$ depends as in the statement of Theorem \ref{main1},
while $q \in [p-1 ,b)$, and $\sigma \in (0,\sigma(q))$. Now estimate
\rif{apest} follows from \rif{apbase},\rif{convergenzeu} and
\rif{cisiamo1}, together with a standard lower semicontinuity
argument to handle the left hand sides of \rif{apbase},
\rif{cisiamo1}. We conclude showing that $u$ solves \rif{Dir1} in
the sense of \rif{desol}. The a priori estimate \rif{cisiamo1}
allows for a quick derivation of this fact. Indeed, thanks to
Rellich's compactness theorem in the case of fractional Sobolev
spaces \cite{Ad}, we have that, up to extracting a diagonal
subsequence, $Du_k$ strongly converges to $Du$ in
$L^{t}_{\loc}(\Omega,\er^n)$ for every $t < nq/(n-\sigma(q))$, and
on the other hand note that $nq/(n-\sigma(q))=n(p-1)/(n-1)>p-1$.
Taking into account the growth condition \rif{veracre}, and that
$f_k \rightharpoonup \mu$ by \rif{convergenzabase}, we can pass to
the limit in \rif{Dirapp}$_1$ using \rif{veracre} and a well known
variant of Lebesgue's dominated convergence theorem, getting that
$u$ finally satisfies \rif{desol}. The proof of Theorem \ref{main1}
is now complete, and estimate \rif{apest} is also proved. It remains
to prove \rif{stimaloc1}, to this aim we use a scaling argument.
Take $B_R \subset \Omega$, let $u \in W^{1,p}(\Omega)$ be the
solution to \trif{Dirapp2} with a fixed $f$, and scale it back as in
\rif{scalaraggio} in order to obtain $\tu(y)$, a solution in $B_1$.
Now observe that we may apply Lemma \ref{222it} to $\tu$ since the
whole argument of the lemma is local, and makes no use of the
boundary information on the solution considered. Therefore estimate
\rif{bettere2} applied to $\tu$ with $\Omega' \equiv B_{1/2}$ gives
$$
[ D\tu ]_{\sigma/q ,q;B_{1/2}}^q \leq c \||D\tu|+s\|_{L^{q}(B_1)}^q
+c\|\tilde{f}\|_{L^1(B_1)}^{q/(p-1)} \;,
$$
for every $\sigma < \sigma(q)$; here we also used \rif{triviale}
while $c \equiv c(n,p,\ratio,\sigma,  q).$ Scaling back to $B_R$,
observing that $[ D\tu ]_{\sigma/q ,q;B_{1/2}}^q = R^{\sigma-n}[ Du
]_{\sigma/q ,q;B_{R/2}}^q$ we have
$$
[Du]_{\sigma/q, q;B_{R/2}}^q\leq c
R^{-\sigma}\||Du|+s\|_{L^{q}(B_R)}^q+c
R^{\sigma(q)-\sigma}\|f\|_{L^1(B_R)}^{q/(p-1)}\;.
$$
We used that $n-\sigma(q)=q(n-1)/(p-1)$ by \rif{ssiidef}. Writing
the latter estimate for $u \equiv u_k$, and using the approximation
scheme of Section 5 and in particular \rif{convergenzabase} and
\rif{total2}, we finally obtain estimate \rif{stimaloc1}.
\end{proof}
\begin{remark}\label{integerex} The crucial case in the proof Theorem
\ref{main1} is actually \rif{pipi}. The case \rif{fra1} can be
obtained by embedding from \rif{pipi} \cite{rusi}, 2.2.3. Indeed,
for a space $W^{\alpha,q}$ the number $\alpha-n/q$ is called {\em
integer dimension}; all the spaces in \rif{fra1} share the same
integer dimension if $\ep=0$, and this allows for using a suitable
embedding. We gave here a self-contained proof, which is on the
other hand even shorter than the one using abstract embedding
theorems for Besov spaces.
\end{remark}
\begin{proof}[Proof of Theorem \ref{noncz}.] The proof goes along
the lines of the one for Theorem \ref{main1}. Take $q=p-1$ in Lemma
\ref{222it}, in such a way that now \rif{ildelta0} gives
\eqn{deltalele}
$$\delta = \frac{p}{2(p-1)}\;.$$
Now we proceed as for the proof of Theorems
\ref{main1}-\ref{main1es}, again applying Lemma \ref{222it} first to
the approximating solutions $u_k$ defined in Section 5, and then
passing to the limit $k \nearrow \infty$ the resulting a priori
estimates. The equality in \rif{deltalele} together with
\rif{start22} finally leads to
$$
V(Du) \in W^{t,\frac{2(p-1)}{p}}_{\loc}(\Omega, \er^n) \qquad
\mbox{for every} \qquad t < \frac{p}{2(p-1)}\,,
$$
which establishes \rif{nonfra1} in Theorem \ref{noncz}. In order to
get \rif{nonapest} and therefore completing the proof we just use
the a priori estimate (\ref{bettere2}) for the approximate solutions
$u_k$, and then we let $k\nearrow  \infty$ as for the proof of
Theorem \ref{main1es}.
\end{proof}
\begin{remark}
As for \rif{ssii} we can prove, using Lemma \ref{222it}, that the
solution $u$ found in Theorem \ref{noncz} satisfies
$$V(Du) \in
W^{\frac{p\sigma(q)}{2q}-\ep,\frac{2q}{p}}_{\loc}(\Omega,\er^n),
\qquad \mbox{for every}\ \ \ep >0\;,
$$
for the values of $q$ described in \rif{ssii}.
\end{remark}
\begin{proof}[Proof of Corollary \ref{nondegenere}.] This is based
on inequality \rif{V}. Set $q_0= 2(p-1)/p$; since $p \geq 2$ we have
\begin{eqnarray*}
&& s^{\frac{(p-2)q_0}{2}}\int_{\Omega^{\prime}}
\int_{\Omega^{\prime}} \frac{|Du(x)-Du(y)|^{q_0}}{|x-y|^{n+1-\ep}}\
dx\, dy\\ \nonumber & &\qquad  \leq \int_{\Omega^{\prime}}
\int_{\Omega^{\prime}}
\frac{\left[(s+|Du(x)|+|Du(y)|)^{p-2}|Du(x)-Du(y)|^2\right]^{q_0/2}}{|x-y|^{n+1-\ep}}\
dx\,  dy\\ \nonumber && \qquad  \leq c(n,p)\int_{\Omega^{\prime}}
\int_{\Omega^{\prime}}
\frac{|V(Du(x))-V(Du(y))|^{q_0}}{|x-y|^{n+1-\ep}}\ dx\,  dy\;,
\end{eqnarray*}
and the proof is concluded using estimate \rif{nonapest}.
\end{proof}
\section{The capacitary case}
Here we give the proof of Theorem \ref{main4}, that will be along
the lines of the one for Theorem \ref{main1}; therefore we shall
confine to report the necessary modifications. The main point here
is that we do not need estimates below the growth exponent like
\rif{revca} and \rif{zz2}, as the solution $u$ to \rif{Dir1} is
uniquely determined in $W^{1,p}_0(\Omega)$; for the same reason no
approximation scheme as in Section 5 is needed.

As for \rif{ildelta0} we first we need to define \eqn{ildelta2}
$$
\delta:= \frac{\sigma(p)}{2},\quad
\gamma(t):=\frac{\delta}{\delta+1-t} \quad t \in [0,\delta+1), \quad
\bar{\mu}:=M^{\frac{1}{p-1}}|\mu|\;,
$$
where $M$ appears in \rif{morrey}, and $\sigma(p)$ is defined in
\rif{fra123}. Next lemma is the analog of Lemma \ref{2it}.
\begin{lemma}\label{3it}
Let $u \in \Wp_0(\Omega ,\er^N)$ be the unique solution to
\trif{Dir1}, under the assumptions \trif{asp}, \trif{lipi} and
\trif{morrey} for $\theta < p$. Assume that \eqn{initial2}
$$ V(Du)
\in W^{t,2}_{\rm loc} (\Omega ,\er^n)\,,\qquad \mbox{for some} \ \ t
\in [0,\delta)\;,$$ where $\delta$ is as in \trif{ildelta2}, and
that for every couple of open subsets $\Omega^{\prime} \CC
\Omega^{\prime \prime} \CC \Omega$ there exists a constant $c_1
\equiv c_1(\distc)$ such that \eqn{eq1it2}
$$
[V(Du)]_{t ,2;\Omega^{\prime}}^{2} \leq c_1 \int_{\Omega^{\prime
\prime}} \! ( |Du|^{p}+s^p)\, dx +c_1\bar{\mu}(\Omega'')\;.
$$
 Then \eqn{final2}
 $$ V(Du) \in W^{\tilde{t}, 2}_{\rm loc} (\Omega ,\er^n)\,, \qquad
\mbox{for every}\ \ \tilde{t} \in [0,\gamma(t))\;,$$ where
$\gamma(\cdot)$ is in \trif{ildelta2}. Moreover, for every couple of
open subsets $\Omega^{\prime} \CC \Omega^{\prime \prime} \CC \Omega$
there exists a constant
 $c \equiv c(n,p,\ratio, \distc,\tilde{t},c_1)$
such that \eqn{bettere2222}
$$
[V(Du)]_{\tilde{t},2;\Omega^{\prime}}^{2} \leq c
\int_{\Omega^{\prime \prime}} \! ( |Du|^{p}+s^p  ) \, dx+
c\bar{\mu}(\Omega'')\;.
$$
\end{lemma}
\begin{proof} The proof follows the one of Lemma
\ref{2it}, therefore we shall keep the notation introduced there,
giving the suitable modifications. Let us firstly treat the case
$p\not=n$. Once again $h$, $v,v_0$ are as in \rif{ha} and
\rif{cdir2}-\rif{cdir2bis}, respectively. As for \rif{lage},
\begin{eqnarray}
\nonumber \int_{B} \! |\tau_{i,h}V(Du)|^{2} \, dx &\leq& c\int_{B}
\! |\tau_{i,h}V(Dv_0)|^{2} \, dx +c\int_{\hat B} \!
|V(Du)-V(Dv)|^{2}\, dx \nonumber\\ &&\qquad  +c\int_{ 2B} \!
|V(Dv)-V(Dv_0)|^{2} \, dx=: I + II + III\;. \label{lage222}
\end{eqnarray}
The term $III$ is estimated via \rif{dopodopo}, while for $II$ we
use Lemmata \ref{coco5}-\ref{coco52}:
$$
\int_{\hat B} |V(Du)-V(Dv)|^{2} \, dx\leq c \bar{\mu}(\hat
B)|h|^{\beta\sigma(p)}\, \;.
$$
Therefore, as $\sigma(p)\leq 2$ when $p\geq 2$, we have
\eqn{IIIII222} $$ \int_{\hat B} |V(Du)-V(Dv_0)|^{2} \, dx + II + III
\leq
 c |h|^{\beta2 \delta} \left\{\int_{\hat B} (s^p+|Du|^p)\, dx
+\bar{\mu}(\hat B)\right\}\;.$$ As for $I$ we shall simply use
estimate \rif{Idue}. We again distinguish two cases:

{\em Case } $t=0$. Taking $z_0=0$ we have, using \rif{IIIII222} and
\rif{elemV}
\begin{eqnarray} \int_{8B} \! |V(Dv_0)-V(z_0)|^{2} \,
dx &\leq & c \int_{\hat B} (s^p+|Du|^p)  \, dx +\int_{8B}
|V(Du)-V(Dv_0)|^{2}  \, dx\nonumber \\
&\leq & c \int_{\hat B} (s^p+|Du|^p)\, dx +c\bar{\mu}(\hat B)
\;.\label{coplno222}
\end{eqnarray}

{\em Case } $t>0$. In this case we choose $z_0$ as in \rif{copl2}.
Again we estimate
\begin{eqnarray} \int_{8B} \! |V(Dv_0)-V(z_0)|^{2} \,
dx &\leq & c \int_{8B} \! |V(Dv_0)-V(Du)|^{2} \,
dx\nonumber \\
&& \qquad +c \int_{8B} \! |V(Du)-V(z_0)|^{2}\, dx\;.\label{prim2222}
\end{eqnarray}
Using Proposition \ref{frapoes}, together with  \rif{initial2} and
the choice \rif{copl2}, gives \eqn{copl3222}
$$ \int_{8B} \!
|V(Du)-V(z_0)|^{2} \, dx  \leq  c|h|^{ \beta 2t}[V(Du)]_{t,2;8
B}^{2}\;.
$$
Combining \rif{copl3222} and \rif{IIIII222} with \rif{prim2222} we
have, as $t< \delta$ \eqn{copl10222}
$$
\int_{8B} \! |V(Dv_0)-V(z_0)|^{2} \, dx \leq c |h|^{\beta
2t}\left[\int_{\hat B} (s^p+|Du|^p)\, dx +\bar{\mu}(\hat B)+
[V(Du)]_{t,2;\hat B}^{2}\right].
$$
Now let us set for any measurable set $A \CC \Omega$
$$
\lambda(A):=\int_{A} (s^p+|Du|^p)\, dx +\bar{\mu}(A) +
\chi(t)[V(Du)]_{t,2;A}^{2}\;,$$ where again $\chi(t)=0$ if $t=0$,
and $\chi(t)=1$ if $t>0$. Summarizing \rif{lage222}, \rif{coplno222}
and \rif{copl10222} we have
$$
I \leq c \int_{B} \! |\tau_{i,h}V(Dv_0)|^{2} \, dx \leq c
|h|^{2[(1-\beta)+ t\beta]} \lambda(\hat B)\;.
$$
Combining this last estimate with \rif{IIIII222} and \rif{lage222}
we finally find
$$
\int_{B} \! |\tau_{i,h}V(Du)|^{2} \, dx \leq c
\left[|h|^{2[(1-\beta)+ t\beta]}+ |h|^{\beta 2\delta}
\right]\lambda(\hat B)\;.
$$
From now on we can proceed with the covering argument adopted in the
proof of Lemma \ref{2it}, up to formula \rif{quasi22}, arriving at
\eqn{quasi}
$$
\int_{\Omega'} \! |\tau_{i,h}V(Du)|^{2} \, dx \leq c
\left[|h|^{2[(1-\beta)+ t\beta]}+ |h|^{\beta 2\delta}
\right]\lambda(\Omega'')\;.
$$
Taking $\beta=\gamma(t)/\delta \in (0,1)$ now yields
$$
\int_{\Omega'} \! |\tau_{i,h}V(Du)|^{2} \, dx \leq c
|h|^{2\gamma(t)}\lambda(\Omega'')\;,
$$
that is the analog of \rif{analog}. From this point on the proof
proceeds as for Lemma \ref{2it}, and the case $p\not=n$ is complete.
As for $p=n$, Lemma \ref{coco52} allows to re-do the whole proof
where this time $\delta:=\sigma'/2$, for any $\sigma' \in
(0,\sigma(p))$; therefore we obtain $ V(Du) \in W^{\tilde{t},
2}_{\rm loc} (\Omega ,\er^n)$ for every $\tilde{t} <
(\sigma'/2)/[(\sigma'/2)+1-t]$. Since $\sigma'$ can be chosen
arbitrarily close to $\sigma(p)$ the statement follows again, and
the proof is complete. In particular \rif{final2} follows from the
fact that the open subsets considered are arbitrary.
\end{proof}
\begin{proof}[Proof of Theorem \ref{main4}] The proof goes as the
one for Lemma \ref{222it}, but directly for the solution $u$ to
\rif{Dir1}. Applying repeatedly Lemma \ref{3it} with $t \equiv t_k$
as in Lemma \ref{222it}, and $\{t_k\}$ is the sequence defined in
\rif{succ} with $\delta= \sigma(p)/2$, we get that $ V(Du) \in
W^{t_k,2}_{\rm loc} (\Omega ,\er^n)$ for every $k \in \en$, with a
corresponding estimate of the type \rif{eq1it2}. The assertion
finally follows observing that this time $t_k\nearrow \sigma(p)/2$,
passing to from $V(Du)$ to $Du$ as in \rif{converti}, and using
\rif{apestsub2222} to get the global bound in \rif{apestsub}.
\end{proof}
\section{Morrey estimates}
In this section we give the proofs of Theorems \ref{mainmorrey} and
\ref{main6}. We shall actually argue as follows: we first prove
Theorem \ref{mainmorrey} in the special case $q < b$, at least as a
priori estimate. This will allows us to prove Theorem \ref{main6}
immediately, and also Theorem \ref{main5} in the next section. In
turn Theorem \ref{main5} will finally imply Theorem \ref{mainmorrey}
for the full range $q< m$; compare with \rif{catena}. Therefore we
shall start with
\begin{lemma} \label{mainmorrey2}
Let $u \in W^{1,p}_0(\Omega)$ be the solution to \trif{Dirapp2} for
a fixed $f \in L^{\infty}(\Omega)$, under the assumptions \trif{asp}
with $p \leq n$. Then with $$p-1 \leq q < \frac{n(p-1)}{n-1}=b,
\quad \mbox{and} \quad \delta(q):= \frac{q(\theta-1)}{p-1}\;,$$ as
in \trif{ildelta}, for every couple of open subsets $\Omega' \CC
\Omega'' \CC \Omega$ there exists $c \equiv
c(n,p,\ratio,q,\dist(\Omega',\partial \Omega''))$ such that whenever
$\theta \in [p,n]$ \eqn{sstimamo}
$$
\||Du|+s\|_{L^{q,\delta(q)}(\Omega')} \leq
c\||Du|+s\|_{L^{q}(\Omega'')}+c\|f\|_{L^{1,\theta}(\Omega'')}^{1/(p-1)}\;.
$$
Moreover there exists $c \equiv
c(n,p,\ratio,q,\dist(\Omega',\partial \Omega),\Omega)$ such that
\eqn{stimamorreypre}
$$
\||Du|+s\|_{L^{q,\delta(q)}(\Omega')} \leq c
\|f\|_{L^{1}(\Omega)}^{1/(p-1)}+c\|f\|_{L^{1,\theta}(\Omega)}^{1/(p-1)}
+cs|\Omega|^{1/q}\;.
$$
\end{lemma}
\begin{proof} We shall apply a standard comparison technique to get
Morrey estimates. Let us take $B_R \CC \Omega''$ with $R \leq 1$,
and define $v \in u + W_0^{1,p}(B_R)$ as the unique solution to
\rif{Dirc1}. Using Lemma \ref{superrev}, estimate \rif{decmo}, for
any $\varrho \in (0,R)$ \eqn{primodec}
$$
\int_{B_\varrho} (|Dv|^q+s^q) \, dx \leq c
\left(\frac{\varrho}{R}\right)^{n-q+\beta q} \int_{B_{R}} (|Dv|^q
+s^q)\, dx\;,
$$
where $c\equiv c(n,p,\ratio,q)$, and $ \beta \equiv
\beta(n,p,\ratio)\in (0,1]$. Now we compare $u$ and $v$ in $B_R$,
that is, using the latter estimate
\begin{eqnarray}
&& \int_{B_\varrho} (|Du|^q+s^q) \, dx \leq  c \int_{B_\varrho}
(|Dv|^q+s^q) \, dx+ c \int_{B_\varrho} |Dv-Du|^q \, dx \nonumber \\
&&\qquad \qquad \leq   c \left(\frac{\varrho}{R}\right)^{n-q+\beta
q} \int_{B_{R}} (|Dv|^q
+s^q)\, dx + c \int_{B_R} |Dv-Du|^q \, dx\nonumber \\
&&\qquad \qquad  \leq c \left(\frac{\varrho}{R}\right)^{n-q+\beta q}
\int_{B_{R}} (|Du|^q +s^q)\, dx + c \int_{B_R} |Dv-Du|^q \,
dx\;.\label{primo2}
\end{eqnarray}
Using Lemma \ref{coco2}, with $c \equiv c(n,p,\ratio,q)$, and $q \in
[p-1,b)$, we get \eqn{primo3}
$$
\int_{B_R} |Du-Dv|^q \, dx\leq c
\|f\|_{L^{1,\theta}(B_R)}^{\frac{q-p+1}{p-1}} \int_{B_R} |f|\, dx \,
R^{\sigma(q,\theta)}\leq
c\|f\|_{L^{1,\theta}(B_R)}^{\frac{q}{p-1}}R^{n-\delta(q)}\;.
$$
Observe that $f \in L^{\infty}$, therefore
$\|f\|_{L^{1,\theta}(\Omega)}< \infty$. Combining \rif{primo2} and
\rif{primo3} yields \eqn{stringata}
$$
\int_{B_\varrho} (|Du|^q+s^q) \, dx\leq c
\left(\frac{\varrho}{R}\right)^{n-q+\beta q} \int_{B_R} (|Du|^q
+s^q)\, dx +c\|f\|_{L^{1,\theta}(B_R)}^{\frac{q}{p-1}}
R^{n-\delta(q)}\;,
$$
where $c \equiv c(n,p,\ratio)$. Observe now that $\theta \geq p$
implies $n-q\geq n- \delta(q)$, therefore we can apply Lemma
\ref{lemmaiterazione} with the choice $$\varphi(t):= \int_{B_t}
(|Du|^q+s^q) \, dx, \qquad \qquad \B:=
\|f\|_{L^{1,\theta}(B_R)}^{\frac{q}{p-1}}\;,$$ and $
\delta_0:=n-q+\beta q> \delta_1:= n-q+\beta q/2>\gamma\equiv
n-\delta(q)$ in order to have, after an elementary manipulation
\begin{eqnarray} && \int_{B_\varrho} (|Du|^q+s^q) \, dx\nonumber  \\ && \leq
c_1\left\{c_*(R)\left(\frac{\varrho}{R}\right)^{\beta q
/2}\int_{B_R} (|Du|^q+s^q) \, dx +
\|f\|_{L^{1,\theta}(B_R)}^{\frac{q}{p-1}}\right\}\varrho^{n-\delta(q)}\;,\label{servelittle}
\end{eqnarray}
for every $\varrho \leq R$, where $c_1 \equiv c_1(n,p,\ratio,q)$,
and $c_*(R)=R^{\delta(q)-n}$. Now take $\bar R:=\dist(\Omega',
\partial \Omega'')/4$, then use \rif{servelittle} on the generic ball of radius $\bar{R}$ centered in
$\Omega'$; of course such a ball is contained in $\Omega''$. Also
observe that such a choice of $\bar{R}$ determines $c_* \equiv
c_*(n,\dist(\Omega',
\partial \Omega''))$ in \rif{servelittle}. All in all such choices
give \eqn{morreydefi}
$$ \int_{B_\varrho} (|Du|^q+s^q) \, dx \leq
c\left[\||Du|+s\|_{L^{q}(\Omega'')}^q +
\|f\|_{L^{1,\theta}(\Omega'')}^{\frac{q}{p-1}}\right]\varrho^{n-\delta(q)}\;,
$$ with $c \equiv
c(n,p,\ratio, q,\dist(\Omega',
\partial \Omega''))$, for any $\varrho \leq \bar{R}$.
This procedure, and an elementary estimation involving the
definition in \rif{dede1mo}, yield \rif{sstimamo} with the specified
dependence of the constant. More precisely, \rif{morreydefi} is
satisfied for $\varrho \leq \bar{R}$, but then is satisfied also for
any ball $B_{\varrho} \subset \Omega'$, with $\varrho \leq 1$,
modulo increasing the constant $c$ of the factor
$\bar{R}^{\delta(q)-n}$ in the case $\bar{R}< 1$; recall that $\bar
R:=\dist(\Omega',
\partial \Omega'')/4$. Finally, in order to get \rif{stimamorreypre}, fix $\Omega' \CC \Omega$, and determine
$\Omega''$ according to Lemma \ref{stanom}; at this point
\rif{stimamorreypre} follows using \rif{apscala} in \rif{sstimamo},
since $\dist(\Omega',
\partial \Omega'')=\dist(\Omega',
\partial \Omega)/2$.
\end{proof}
\begin{proof}[Proof of Theorem \ref{main6}] Take $\Omega' \CC \Omega$ as in the statement of the Theorem, and determine
$\Omega''$ according to Lemma \ref{stanom}. We go back to the proof
of Lemma \ref{mainmorrey2}, and apply the arguments  to $u_k$, that
is the solution to \rif{Dirapp}, with such a choice of $\Omega',
\Omega''$. We recall that everywhere both $\dist(\Omega',
\partial \Omega'')$ and $\dist(\Omega'',
\partial \Omega)$ depend on $\dist(\Omega',
\partial \Omega)$ via \rif{distanze}. We start from
\rif{servelittle}; as by \rif{ildelta} $\delta(q)=q$ when
$\theta=p$, we use Poincar\'e's inequality in order to estimate the
left hand side of \rif{servelittle} from below. With $c_1$ being the
one in \rif{servelittle} up to multiplicative constant $c(n,q)$, we
have
\begin{eqnarray}
&&\nonumber  \mean{B_\varrho} |u_k - (u_k)_{B_\varrho}|^q \, dx \\
&& \qquad \leq
c_1\left\{R^{q-n}\left(\frac{\varrho}{R}\right)^{\beta q
/2}\int_{B_R} (|Du_k|^q+s^q) \, dx +
\|f_k\|_{L^{1,\theta}(B_R)}^{\frac{q}{p-1}}\right\}\;.\label{servelittle2}
\end{eqnarray}
Now, fix $\Omega' \CC \Omega'' \CC \Omega$ as in the proof of Lemma
\ref{mainmorrey2}, and using the same argument used to prove Morrey
regularity in the previous proof we find $$ [u_k]_{BMO(\Omega')}
\leq
c\|f_k\|_{L^{1}(\Omega)}^{1/(p-1)}+c\|f_k\|_{L^{1,\theta}(\Omega)}^{1/(p-1)}+cs|\Omega|^{1/q}\;,$$
with $c \equiv c(n,p,\ratio, q,\dist(\Omega',
\partial \Omega))$. Letting $k \nearrow \infty$, and using of \rif{total}, we finally obtain
$$ [u]_{BMO(\Omega')} \leq
c[|\mu|(\Omega)]^{1/(p-1)}+cM^{1/(p-1)}+cs|\Omega|^{1/q}\;.$$ Now
\rif{apbmo} is finally proved combining the last estimate with the
following trivial consequence of \rif{morrey}: \eqn{trimo}
$$
|\mu|(\Omega)\leq  [\textnormal{diam}(\Omega)]^{n-\theta}M\;,
$$
and taking into account that $c$ may depend on $\Omega$ too.

In order to prove the local VMO regularity we assume that $\mu$
satisfies \rif{locuni} locally uniformly in the sense of Definition
2 in Section 2.4. In order to conclude it suffices to prove that:
For every $\Omega' \CC \Omega$ and every $\ep >0$ there exist
$\bar{k}\equiv \bar{k}(\ep,\dist(\Omega',
\partial \Omega))\in \en $ and
$\bar{\varrho}\equiv \bar{\varrho}(\ep,\dist(\Omega',
\partial \Omega))\in \en $, possibly
also depending on $n,p,\ratio,q,s,\Omega$, such that \eqn{finalebmo}
$$
\mean{B_{\varrho}} |u_k - (u_k)_{B_{\varrho}}|^q \, dx\leq \ep,
\qquad \qquad k \geq \bar{k}, \ \ \ \ \varrho \leq \bar{\varrho}\;,
$$
whenever $B_{\varrho} \CC \Omega''$ is ball centered in $\Omega'$.
This with \rif{convergenzeu} will finally prove the whole theorem as
$\Omega' \CC \Omega$ is arbitrary. Using \rif{apbase} and
\rif{total2} with \rif{servelittle2} we have
\begin{eqnarray}
\nonumber && \mean{B_\varrho} |u_k - (u_k)_{B_{\varrho}}|^q \, dx\\
&& \qquad \leq
c_1\left\{R^{q-n}\left(\frac{\varrho}{R}\right)^{\beta q
/2}\left[[|\mu|(\Omega)]^{\frac{q}{p-1}} +s^q|\Omega|\right]+
\|\mu\|_{L^{1,\theta}(B_{R+1/k})}^{\frac{q}{p-1}}\right\}\;.\label{servelittle3}
\end{eqnarray}
Determine a positive radius $\bar{R}\leq \dist(\Omega',\partial
\Omega'')/4$, depending on $\ep, \dist(\Omega',\partial \Omega)$ and
on $n,p,\ratio,q$, such that $|\mu|(B_r)\leq (2c_1)^{-1}\ep r^{n-p}$
whenever $r \leq 2\bar{R}$ and $B_r \subset \Omega''$. This implies
$\|\mu\|_{L^{1,\theta}(B_{2\bar R})}\leq (2c_1)^{-1}\ep $ whenever
$B_{2\bar{R}}$ is centered in $\Omega'$. Indeed this and
$\bar{R}\leq \dist(\Omega',\partial \Omega'')/4$ imply
$B_{2\bar{R}}\subset \Omega''$. From now on all the balls considered
will be centered in $\Omega'$. Taking $\bar{k}\equiv
\bar{k}(\ep,\dist(\Omega',
\partial \Omega)) \in \en$, also depending on $n,p,\ratio,q$,
such that $1/\bar{k}\leq \bar{R}$ we have \eqn{servelittle4}
$$
c_1\|\mu\|_{L^{1,\theta}(B_{\bar R+1/k})}\leq \ep/2\;.
$$
This fixes $\bar{k}$ in \rif{finalebmo}. From now on we shall use
\rif{servelittle3} with $R \equiv \bar{R}$. Now take
$\bar{\varrho}\equiv \bar{\varrho}(\ep,\dist(\Omega',
\partial \Omega))\leq \bar{R}$, also
depending on $n,p,\ratio,q,s, \Omega$, in order to have
\eqn{servelittle5}
$$
c_1 \bar{R}^{q-n}\left(\frac{\bar{\varrho}}{\bar{R}}\right)^{\beta q
/2}\left[[|\mu|(\Omega)]^{\frac{q}{p-1}}+s^q|\Omega|\right]\leq
\ep/2\;.
$$
This fixes $\bar{\varrho}$  in \rif{finalebmo}. We finally obtain
\rif{finalebmo} merging \rif{servelittle4}-\rif{servelittle5} to
\rif{servelittle3}, the latter used with $\bar{R}\equiv R$, and
$\varrho \leq \bar{\varrho}$.
\end{proof}
\begin{proof}[Proof of Theorem \ref{mainmorrey}] As usual we shall
proceed deriving a priori estimates, therefore let $u \in
W^{1,p}(\Omega)$ be the solution to \trif{Dirapp2} for a fixed $f\in
L^{\infty}(\Omega)$. We shall use the estimates from the proof of
Theorem \ref{main5} below, as explained at the beginning of the
section, therefore this proof should be read after the one of
Theorem \ref{main5}. Let $B_R \CC \Omega $, with $R \leq 1$. By
Lemma \ref{hollo} with $q \in (1,m)$
\begin{eqnarray}
\nonumber \int_{B_R} |Du|^q \,  dx &\stackrel{\rif{hollore}}{\leq} &
m\left(m-q\right)^{-1}|B_R|^{1-\frac{q}{m}}
\|Du\|_{\MM^m(B_R)}^q\\
\nonumber &\stackrel{\rif{mama}}{\leq} & c R^{n-\frac{q\theta}{m}}\|Du\|_{\MM^{m,\theta}(B_R)}^q\\
&\stackrel{\rif{stimetta4}}{\leq}&
c\left[\|f\|_{L^{1}(\Omega)}^{\frac{q}{p-1}}+
\|f\|_{L^{1,\theta}(\Omega)}^{\frac{q}{p-1}}+s^q
|\Omega|^{\frac{q}{m}}\right]R^{n-\delta(q)}\;,\label{lolost}
\end{eqnarray}
where $c \equiv c(n,p,\ratio,q, \dist(B_R ,
\partial\Omega), \Omega).$ We used that $q\theta/m=\delta(q)$, see
\rif{ildelta}. Therefore taking the supremum over all possible such
balls with $B_R \CC \Omega'$ we have
$$
\|Du\|_{L^{q, \delta(q)}(\Omega')} \leq c
\|f\|_{L^{1}(\Omega)}^{1/(p-1)}+c
\|f\|_{L^{1,\theta}(\Omega)}^{1/(p-1)}+cs |\Omega|^{1/m}\;,
$$
$c \equiv c(n,p,\ratio,q, \dist(\Omega' ,
\partial\Omega), \Omega).$ The assertion follows once again via the approximation scheme of
Section 5, a lower semicontinuity to handle the left hand side of
the latter estimate, and using \rif{trimo} as for Theorem
\ref{main5}.
\end{proof}
\section{Marcinkiewicz estimates} This section contains the proof of
Theorem \ref{main5}. One of our starting points here will be the
brilliant technique for proving
 $\MM^{n}$-estimates introduced in
\cite{DHM2} (case $p=n$, that implies $\theta=p=n$). We shall use a
delicate combination of some the arguments from the latter paper
with the Morrey space estimates of Section 8, a direct comparison
argument on certain Calder\'on-Zygmund type balls, and finally a
modification of some ideas from \cite{Caffpe, KM}. A different,
elegant approach to $\MM^n$ estimates based on a suitable version of
Gehring's lemma in Marcinkiewicz spaces has been recently given in
\cite{KSZ}. Let us emphasize here the fact that our technique is
robust enough to catch the borderline case $\theta =p$, and
therefore to get the limiting regularity \rif{lolo2}.

As everywhere else, we shall derive a priori estimates and in the
following $u \in W^{1,p}_0(\Omega)$ is a solution to \rif{Dirapp2}
for a fixed $f\in L^{\infty}(\Omega)$; we assume of course that
$\|f\|_{L^{1,\theta}(\Omega)}>0$, otherwise all assertions
trivialize. To begin with the proof let us consider two open subsets
$\Omega' \CC \Omega''\CC \Omega$. Take a ball $B_0$ with radius
$R_0\leq 1/2$, such that $2B_0 \CC \Omega''$. We use the restricted
maximal function of $f$ relative to $\QQ$, that is
$$
M(f)(x)\equiv M_{\QQ}(f)(x):= \sup_{\ x \in B,B \subseteq 2B_0}
\mean{B} |f(y)|\ dy\;,
$$
where $B$ is a ball varying amongst all possible ones in $2B_0$. The
weak $(1,1)$ estimate
$$
|\{x \in 2B_0 \ : \ M_{\QQ}(f)(x)> \lambda \}| \leq
\frac{c(n)}{\lambda} \int_{\QQ}|f(y)|\ dy \qquad \forall\ \lambda
>0\;,
$$
holds, see for instance \cite{BI}, and it immediately follows that
 \eqn{weakes}
$$
|\{x \in 2B_0 \ : \ M_{\QQ}(f)(x)> \lambda \}| \leq
\frac{c\|f\|_{L^{1,\theta}(\QQ)}|B_0|^{1-\theta/n} }{\lambda}\qquad
\forall\ \lambda
>0\;.
$$
Let us fix $R_0 < t < \varrho < 2R_0$. With $\lambda \geq 0$ we
shall denote
$$ E_{\lambda}^{t}:=\{x \in
B_{t} \ : \ |Du(x)|> \lambda\}, \quad E_{\lambda}^{\varrho}:=\{x \in
B_{\varrho} \ : \ |Du(x)|> \lambda\}\;.$$ Here the balls $B_t,
B_{\varrho}$ are concentric to $B_0$, and it obviously holds $ B_0
\subset B_{t} \subset B_{\varrho}\subset 2B_0$. We recall that $b$
is in \rif{boccardoex} and $m$ as in \rif{mingione}, while in the
following $q$ and $q_1$ will be fixed numbers such that $p-1\leq q <
q_1 < b.$

{\em Step 1: Calder\'on-Zygmund type decomposition.} Let us
set\eqn{la0}
$$\lambda_0
:=\left(\mean{\QQ} (|Du|^q+s^q)\, dx \right)^{\frac{1}{q}}\;,$$ and
from now on we shall always take $\lambda$ large enough to have
\eqn{lasott} $$\lambda \geq
4^{n/q}(\varrho-t)^{-n/q}\lambda_0=:\lambda_l\;,$$ unless otherwise
specified. Observe that if $x_0 \in B_{t}$ then
$B(x_0,(\varrho-t)R_0) \subset B_{\varrho}\subset 2B_0$ and
therefore \eqn{la02}
$$
\nonumber\mean{B(x_0,(\varrho-t)R_0)} (|Du|^q+s^q)\, dx
\stackrel{\rif{la0}}{\leq }2^{n}(\varrho-t)^{-n}\lambda_0^q \leq
\lambda^q\;;
$$
in particular \eqn{trila}
$$
s \leq \lambda\;.
$$ Now, let $x_0 \in E_{4\lambda}^{t}$ and define
$$
i(x_0):=\min\left\{i \in \en\ : \ \mean{B(x_0,2^{-i}(\varrho-t)R_0)}
(|Du|^q+s^q)\, dx\geq 4^q\lambda^q\right\}\;.
$$
By \rif{la02} and Lebesgue's differentiation theory for a.e.~$x_0
\in E_{4\lambda}^t$ we have $1\leq i(x_0)< \infty$, and the family
$\{B(x_0,2^{-i(x_0)}\dro R_0)\}$ is a covering of $E_{4\lambda}^t$
up to a negligible set. We may apply Besicovitch covering theorem
\cite{AFP} in order to extract from
$\{B(x_0,2^{-i(x_0)}(\varrho-t)R_0)\}$ a finite number $Q(n)$ of
possibly countable families of mutually disjoint balls
$\{\B_j\}_{j\leq Q(n)}$, $\B_j \equiv\{B^j_i\}$, such that
$E_{4\lambda}^{t}$ is covered by the union of the closure of such
balls up to a negligible set. Rename all these balls in order to
have a new, possibly countable family $\{B_k\}$. We need to observe
that $2B_k \subset B_{\varrho}$ for every $k$; this follows from the
construction, since for a.e.~$x_0 \in B_0$ we have $i(x_0)\geq 1$,
therefore the radius of $B_k$ does not exceed $\dro R_0/2$, and
being $B_k$ centered in $B_t$ then $2B_k \subset B_{\varrho}$
follows. All in all, again by construction the following facts hold:
\eqn{czdec}
$$ E_{4\lambda}^t \subset \bigcup_{k} \overline{B_k} \cup \mbox{negligible set}\;, \quad
\sum_k |E_{\lambda}^{\varrho}\cap B_k|\leq
Q(n)|E_{\lambda}^{\varrho}|\,, \quad 2B_k \subset B_{\varrho}\;$$
and, for every $k \in \en$ \eqn{czdec3}
$$4^q\lambda^q \leq \mean{B_k} (|Du|^q+s^q)\, dx\,, \qquad \mean{2B_k} (|Du|^q+s^q)\, dx<
4^q\lambda^q\;.
$$
Denote by $R_k$ the radius of $B_k$, so that $R_k\leq R_0\leq 1$;
using Lemma \ref{mainmorrey2} gives
$$
4^q\lambda^q \leq \mean{B_k} (|Du|^q+s^q)\, dx \leq c
\||Du|+s\|_{L^{q,\delta(q)}(\QQ)}^q
R_k^{-\frac{q(\theta-1)}{p-1}}\;,
$$
and it follows\eqn{rdecay}
$$
R_k \leq c
K^{\frac{1}{\theta-1}}\lambda^{-\frac{p-1}{\theta-1}},\qquad
K:=\||Du|+s\|_{L^{q,\delta(q)}(\QQ)}^{p-1}+\|f\|_{L^{1,\theta}(2B_0)}\;.
$$

{\em Step 2: A density estimate.} Here we single out one generic
ball $B_k$ and argue under the assumption that there exists $x_k \in
\overline{B_k}$ such that\eqn{maxf}
$$M(f)(x_k)\leq T^{-1}K^{1/(1-\theta)}\lambda^{m}\;,
$$ with $T\geq 1$ to be determined later. Using H\"older's inequality and the fact that $B_k
\subset B_{\varrho}$, we start estimating
\begin{eqnarray}
\nonumber 4^q\lambda^q |B_k|  & \stackrel{\rif{czdec3}}{\leq} &
\int_{B_k} (|Du|^q+s^q)\, dx\\ & \stackrel{\rif{trila}}{\leq} &
2\lambda^q
|B_k \setminus E_{\lambda}^{\varrho}| + \int_{B_k \cap E_{\lambda}^{\varrho}} (|Du|^q+s^q)\, dx \nonumber\\
&\leq & 2\lambda^q |B_k \setminus E_{\lambda}^{\varrho}| + (2|B_k
\cap E_{\lambda}^{\varrho}|)^{1-\frac{q}{q_1}}\left(\int_{B_k \cap
E_{\lambda}^{\varrho}}
 (|Du|^{q_1}+s^{q_1})\, dx\right)^{\frac{q}{q_1}}\;.
\end{eqnarray}
Therefore, another elementary estimation gives\eqn{BB1}
$$
2^q \leq \frac{|B_k \setminus E_{\lambda}^{\varrho}|}{|B_k|} +
2\left[ \frac{|B_k \cap E_{\lambda}^{\varrho}|}{|B_k|}
\right]^{1-\frac{q}{q_1}} \lambda^{-q} \left(\mean{B_k}
(|Du|^{q_1}+s^{q_1})\, dx\right)^{\frac{q}{q_1}}\;.
$$
We now estimate the last integral. To this aim, let us introduce the
comparison function $v_k \in u + W^{1,p}_0(2B_k)$ as the unique
solution to \eqn{Dirck}
$$
\left\{
    \begin{array}{cc}
    -\divo \ a(x,Dv_k)=0 & \qquad \mbox{in} \ 2B_k\\
        v_k= u&\qquad \mbox{on $\partial 2B_k$.}
\end{array}\right.
$$
Now \eqn{ancoras}
$$
\mean{B_k} |Du|^{q_1}\, dx\leq c \mean{B_k} |Du-Dv_k|^{q_1}\, dx +c
\mean{B_k} |Dv_k|^{q_1}\, dx\;,
$$
and we estimate the last two integrals. Using Lemma \ref{coco2} we
find
\begin{eqnarray}
 \nonumber \mean{2B_k}|Du-Dv_k|^{q_1} \, dx  &\stackrel{\rif{comp1m}}{\leq}&
 c \|f\|_{L^{1,\theta}(\QQ)}^{\frac{q_1-p+1}{p-1}} \mean{2B_k} |f|\, dx \,
R_k^{\sigma(q_1,\theta)} \nonumber \\
 \nonumber
& \stackrel{\rif{rdecay}}{\leq}&
c\frac{K^{\frac{1}{\theta-1}}}{\lambda^{m-q_1}}
 \mean{2B_k} |f|\, dx\\&\stackrel{\rif{maxf}}{\leq}&
\frac{ c
 \lambda^{q_1}}{T}\;,\label{didi1}
\end{eqnarray}
where $c \equiv \ c(n,p,\ratio,q_1)$; in a completely similar way we
also get \eqn{didi2}
$$
\mean{2B_k}|Du-Dv_k|^{q} \, dx\leq \frac{c  \lambda^{q}}{T}\;.
$$ On the other hand, since $q_1<p$, using
 H\"older's inequality and \rif{zz2}, we have
\eqn{zz222}
$$
\mean{B_{k}} |Dv_k|^{q_1} \, dx\leq  \left(\mean{B_{k}}
|Dv_k|^{p\chi}\, dx\right)^{\frac{q_1}{p\chi}}\leq c
\left(\mean{2B_{k}} (|Dv_k|^q +s^q)\, dx\right)^{\frac{q_1}{q}} \;.
$$
In the last line $\chi \equiv \chi(n,p,\ratio)>1$ is the higher
integrability exponent such that $Dv \in
L^{p\chi}_{\loc}(2B_k,\er^n)$, that has been determined in Lemma
\ref{superrev}. In turn, since $T \geq 1$, \rif{czdec3} and
\rif{didi2} give \eqn{kaka1}
$$ \mean{2B_k}
|Dv_k|^{q}\, dx  \leq   c \mean{2B_k} |Du-Dv_k|^{q}\, dx +c
\mean{2B_k} |Du|^{q}\, dx\leq  c\lambda^q\;.
$$
Merging \rif{kaka1} and \rif{zz222}, and using \rif{trila}, gives
\eqn{kaka2}
$$
\mean{2B_{k}} |Dv_k|^{q_1} \, dx\leq c\lambda^{q_1}\;.
$$
Connecting \rif{didi1}, \rif{kaka2}, to \rif{ancoras}, and using
again \rif{trila}, yields \eqn{UU1}
$$
\mean{B_k} (|Du|^{q_1}+s^{q_1})\, dx \leq c \lambda^{q_1}\;,
$$
where $c \equiv c (n,p,\ratio,q_1)$. Using this last inequality in
\rif{BB1} gives
$$
2^q \leq \frac{|B_k \setminus E_{\lambda}^{\varrho}|}{|B_k|} +
c_1\left[ \frac{|B_k \cap E_{\lambda}^{\varrho}|}{|B_k|}
\right]^{1-\frac{q}{q_1}} \;,
$$
where $c_1 \equiv c_1 (n,p,\ratio,q_1)$, and therefore, since
$q_1>q$, we have\eqn{density}
$$
\frac{|B_k \cap E_{\lambda}^{\varrho}|}{|B_k|}\geq
\left[\frac{1}{c_1}(2^q-1)\right]^{\frac{q_1}{q_1-q}}=:
\frac{1}{c_2}>0\;,
$$
where $c_2 \equiv c_2 (n,p,\ratio,q,q_1)$; this is the density
estimate we were looking for.

{\em Step 3: Estimates on balls.} We take $H\geq 4^{n/q}\geq 4$ to
be chosen, and estimate the measure of $E_{H\lambda}^t$ splitting as
\begin{eqnarray}
\nonumber|E_{H\lambda}^t| & \leq &|E_{H\lambda}^t\cap\{x \in B_t \ :
\ M(f)(x) \leq T^{-1}K^{1/(1-\theta)}\lambda^{m} \}|\\
&& \quad  + |E_{H\lambda}^t\cap\{x \in B_t \ : \ M(f)(x)>
T^{-1}K^{1/(1-\theta)}\lambda^{m} \}|=: I + II\;.\label{livelli}
\end{eqnarray}
By \rif{weakes}, and the definition of $K$ in \rif{rdecay}, we
immediately have \eqn{IIes}
$$
II \leq c(n)TK^{\frac{\theta}{\theta-1}}
\lambda^{-m}R_0^{n-\theta}\;,
$$
and we concentrate on $I$. To this aim, since $H\geq 4$ by
\rif{czdec} we may estimate \eqn{stimaI}
$$I\leq  \sum I_k := \sum |E_{H\lambda}^t\cap\{x \in \overline{B_k} \ : \
M(f)(x) \leq T^{-1}K^{1/(1-\theta)}\lambda^{m}\}|\;,$$ and in turn
we estimate each $I_k$. Fix one; we may assume there exists $x_k \in
\overline{B_k}$ such that \rif{maxf} holds; otherwise $I_k=0$ and we
are done. By definition of $I_k$
\begin{eqnarray}
I_k \leq |E_{H\lambda}^t\cap \overline{B_k}| &\leq & |\{x \in
\overline{B_k} \ : \
|Du(x)|> H\lambda\}|\nonumber\\
&\leq & |\{x \in \overline{B_k} \ : \ |Du(x)-Dv_k(x)|> H\lambda/2\}|
\nonumber
\\ && \quad + |\{x \in \overline{B_k} \ : \
|Dv_k(x)|> H\lambda/2\}|=: III_k + IV_k\;.\label{Ip}
\end{eqnarray}
Then, keeping in mind the definition of $K$ in \rif{rdecay}
\begin{eqnarray*}
 III_k & \leq& \frac{2^q}{H^q\lambda^q}\int_{B_k}|Du-Dv_k|^q\, dx
\\ & \stackrel{\rif{comp1m}}{\leq} &
\frac{c\|f\|_{L^{1,\theta}(\QQ)}^{\frac{q-p+1}{p-1}}}{H^q\lambda^q}
\int_{2B_k} |f|\, dx\, R_k^{\sigma(q,\theta)}\\\\ &
\stackrel{\rif{rdecay}}{\leq} &
\frac{cK^{\frac{q-p+1}{p-1}}}{H^q\lambda^q} \int_{2B_k} |f|\, dx\,
R_k^{\sigma(q,\theta)}\\
&\stackrel{\rif{rdecay}}{\leq}&
\frac{cK^{\frac{1}{\theta-1}}}{H^q\lambda^{m}} \int_{2B_k} |f|\,
dx\\
&\stackrel{\rif{maxf}}{\leq} & \frac{c_3|2B_k|}{H^qT}\\ &
\stackrel{\rif{density}}{\leq} &   \frac{c_3c_22^n|B_k\cap
E_{\lambda}^{\varrho}|}{H^qT}\;.
\end{eqnarray*}
Let $\chi\equiv \chi(n,p,\ratio)>1$ be as in \rif{zz222}, that is
the number determined in Lemma \ref{superrev}; using \rif{trila} we
have
\begin{eqnarray*}
 IV_k& \leq&
\frac{2^{p\chi}}{H^{p\chi}\lambda^{p\chi}}\int_{B_k}|Dv_k|^{p\chi}\,
dx\\ & \stackrel{\rif{zz2}}{\leq} &
\frac{c|2B_k|}{H^{p\chi}\lambda^{p\chi}} \left(\mean{2B_k}
(|Dv_k|^q+s^q)\, dx\right)^{\frac{p\chi}{q}}\\
&\stackrel{\rif{kaka1}, \rif{trila}}{\leq}& \frac{c_4|2B_k|}{H^{p\chi}}\\
&\stackrel{\rif{density}}{\leq}& \frac{c_4c_22^n|B_k\cap
E_{\lambda}^{\varrho}|}{H^{p\chi}}\;.
\end{eqnarray*}
Connecting the estimates found for $III_k,IV_k$ to \rif{Ip} gives
$$I_k
\leq c_5[H^{-q}T^{-1}+ H^{-p\chi}]|B_k \cap
E_{\lambda}^{\varrho}|\;,$$ with $c_5\equiv c_5(n,p,\ratio,q,q_1)$.
Summing up on $k$ using \rif{stimaI}, \rif{czdec} yields
$$I
\leq c_5Q(n)[H^{-q}T^{-1}+ H^{-p\chi}]|E_{\lambda}^{\varrho}|\;.$$
Merging the latter estimate and \rif{IIes} with \rif{livelli} we
finally have \eqn{lll1}
$$
|E_{H\lambda}^t|\leq c_6
\left[H^{-q}T^{-1}+H^{-p\chi}\right]|E_{\lambda}^{\varrho}|+ c_6T
K^{\frac{\theta}{\theta-1}}\lambda^{-m}R_0^{n-\theta}\;,
$$
where $c_6\equiv c_6(n,p,\ratio,q,q_1)$, while $H\geq 4^{n/q}$ and
$T\geq 1$ are still to be chosen.

 {\em Step 4: Iteration and
a priori estimate.} Let us for a moment assume that $R_0=1/2$; we
shall eventually deal with the general case by means of a scaling
argument. We introduce the level function $l(\cdot,\cdot)$ as
\eqn{ella}
$$l(\lambda,\gamma):=\lambda^{m}|E_{\lambda}^{\gamma}|\,, \hspace{1cm} \mbox{for every}\  \gamma  \in [1/2,1], \ \mbox{and}\ \lambda
> 0\;,$$ and observe that \rif{lll1} can be rephrased as
\eqn{quasiq}
$$
l(H\lambda,t)\leq
c_6\left[H^{m-q}T^{-1}+H^{m-p\chi}\right]l(\lambda, \varrho)+
c_6H^mT K^{\frac{\theta}{\theta-1}}\;.
$$
Now observe that $m\leq p<p\chi$, and equality in the first
inequality occurs iff $p=\theta$; therefore we take $H$ large enough
in order to have $c_6H^{m-p\chi}\leq 1/4$; taking into account the
dependence of $m,\chi$, and $c_6$, this fixes $H \equiv
H(n,p,\ratio,q,q_1)$. Next, take $T$ large enough to balance $H$
i.e.~$T:=4c_6H^{m-q}$, recall that $m>q$; therefore $T \equiv
T(n,p,\ratio,q,q_1)$. With such choices \rif{quasiq} gives
\eqn{quasiq2}
$$l(H\lambda, t)\leq (1/2)l(\lambda, \varrho)
+c_7K^{\frac{\theta}{\theta-1}}\;,$$ with $c_7 \equiv
c_7(n,p,\ratio,q,q_1)$. The last inequality holds whenever $\lambda$
satisfies \rif{lasott}, therefore, taking into account the
definition in \rif{dema} to bound the right hand side, with
$\lambda_l$ as in \rif{lasott} we have
$$\sup_{\lambda\geq \lambda_l} l(H\lambda,
t) \leq (1/2)\|Du\|_{\MM^m(B_{\varrho})}^m
+c_7K^{\frac{\theta}{\theta-1}}\;,$$ and obviously, again by the
definition in \rif{dema}, and \rif{la0}-\rif{lasott}, we have
$$\|Du\|_{\MM^m(B_{t})}^m \leq (1/2)\|Du\|_{\MM^m(B_{\varrho})}^m
+ H^m|2B_0|\dro^{-n m/q}\lambda_0^m+
c_7K^{\frac{\theta}{\theta-1}}\;.$$ Observe also that we are proving
a priori estimates for approximate solutions, and therefore we are
assuming that $u \in W^{1,p}$; since $m \leq p$ in any case it
follows that $\|Du\|_{\MM^m(2B_{0})}$ is finite. We can apply Lemma
\ref{simpfun} with $\varphi(t):= \|Du\|_{\MM^m(B_{t})}^m$,
$R_0=1/2$, and $1/2 < t < \varrho < 1$; this yields \eqn{pezzo12}
$$\|Du\|_{\MM^m(B_0)}^m \leq  c\lambda_0^m+ cK^{\frac{\theta}{\theta-1}}\;,$$
with $c \equiv c(n,p,\ratio,q,q_1)$, as $H$ depends on
$n,p,\ratio,q,q_1$ and $|2B_0| \leq c(n)$. Using the definition of
$\lambda_0$ in \rif{la0} and that $R_0=1/2$, we observe that
\eqn{pezzo1}
$$
\lambda_0^m \leq\left(\mean{\QQ} (|Du|^q+s^q)\, dx
\right)^{\frac{m}{q}}\leq  c\||Du|+s\|_{L^{q,\delta(q)}(\QQ)}^{m}
\;.
$$
Merging \rif{pezzo1} with \rif{pezzo12}, taking into account the
definition of $K$ in \rif{rdecay} we easily obtain \eqn{pezzo1000}
$$
\||Du|+s\|_{\MM^m(B_0)}^m \leq c \left[
\||Du|+s\|_{L^{q,\delta(q)}(\QQ)}^{m}+\|f\|_{L^{1,\theta}(2B_0)}^{\frac{m}{p-1}}
\right]\;.
$$
All this holds provided $R_0=1/2$. The general case $R_0 \in
(0,1/2]$ can be dealt with by scaling, that is: first considering a
general ball $B_0$, then from the very beginning of the proof
reducing the problem to the case $R_0 = 1/2$ scaling as in
\rif{scalaraggio}, then obtaining \rif{pezzo1000} for the scaled
solution $\tilde{u}$ with data $\tilde{f}$, and finally scaling back
to $u$; then using also Lemma \ref{scalamorrey} all this yields, for
any $R_0 \in (0,1/2]$ \eqn{stimetta}
$$
\||Du|+s\|_{\MM^m(B_0)}^m \leq c \left[
\||Du|+s\|_{L^{q,\delta(q)}(\QQ)}^{m}+\|f\|_{L^{1,\theta}(2B_0)}^{\frac{m}{p-1}}
\right]R_0^{n-\theta}\;,
$$
where again $c \equiv c(n,p,\ratio,q,q_1)$. For later convenience
let us observe that taking $q=p-1$ in the latter estimate, and using
that $\delta(p-1)=\theta-1$ by \rif{ildelta}, we have
\eqn{stimettapm1}
$$
\||Du|+s\|_{\MM^m(B_0)}^m \leq c \left[
\||Du|+s\|_{L^{p-1,\theta-1}(\QQ)}^{m}+\|f\|_{L^{1,\theta}(2B_0)}^{\frac{m}{p-1}}
\right]R_0^{n-\theta}\;.
$$
Using \rif{stimetta} together with a standard covering argument, and
an elementary estimation involving the definition in \rif{mama}, we
have \eqn{stimetta2}
$$
\||Du|+s\|_{\MM^{m,\theta}(\Omega')} \leq c
\||Du|+s\|_{L^{q,\delta(q)}(\Omega'')}+c\|f\|_{L^{1,\theta}(\Omega'')}^{1/(p-1)}
\;,
$$
where $c \equiv c(n,p,\ratio,q, \Omega' ,\Omega''),$ since
\rif{distanze} holds. Finally using \rif{sstimamo} in the previous
estimate, and via Lemma \ref{stanom} passing again to outer and
inner open subsets to $\Omega'$ and $\Omega''$ respectively, as
everywhere open subsets $\Omega', \Omega''$ are arbitrary, we
conclude with the desired local a priori estimate \eqn{stimetta3}
$$
\||Du|+s\|_{\MM^{m,\theta}(\Omega')} \leq c
\||Du|+s\|_{L^{q}(\Omega'')}+c\|f\|_{L^{1,\theta}(\Omega'')}^{1/(p-1)}
\;,
$$
for any choice $\Omega'\CC \Omega''$, where $c \equiv
c(n,p,\ratio,q, \Omega' ,\Omega'') $. In turn, with $\Omega' \CC
\Omega$ fixed as in the statement of Theorem \ref{main5}, we can
pick $\Omega''$ in \rif{stimetta3} as prescribed in Lemma
\ref{stanom}, and taking into account \rif{apscala} we get
\eqn{stimetta4}
$$
\||Du|+s\|_{\MM^{m,\theta}(\Omega')} \leq
c\|f\|_{L^{1}(\Omega)}^{1/(p-1)}+
c\|f\|_{L^{1,\theta}(\Omega)}^{1/(p-1)}+cs |\Omega|^{1/q} \;,
$$
where now $c \equiv c(n,p,\ratio,q, \Omega' , \Omega).$ Applying the
latter inequality to $u_k$ from \rif{Dirapp} and taking into account
the approximation scheme of Section 5, and in particular
\rif{total}, as in the proof of Theorem \ref{main1} we
get\eqn{stimetta5}
$$
\||Du|+s\|_{\MM^{m,\theta}(\Omega')} \leq
c[|\mu|(\Omega)]^{1/(p-1)}+ cM^{1/(p-1)}+cs |\Omega|^{1/q} \;,
$$
where now $u$ is the solution to the original problem \rif{Dir1}
constructed in Section 5. The assertion of Theorem \ref{main5} with
estimate \rif{apestlo} follow plugging estimate \rif{trimo} in
\rif{stimetta5}. Just one remark about the convergence of the
approximating solutions $u_k$. In the proof of Theorem \ref{main1}
we used the higher (fractional) differentiability of solutions to
pass to the limit via compactness; this information is not available
here since in Theorem \ref{main5} we are just assuming a measurable
dependence of the coefficients, and not \rif{lipi}, which was in
turn necessary to get differentiability of $Du$. In the present case
the converge of the solutions $u_k$ can be nevertheless obtained
exactly as in \cite{boccardo, BG1, DHM1}.\qed
\begin{remark}[A local estimate] Estimate \rif{apestlo} has a local
companion. More precisely, having \rif{stimetta3} at our disposal,
we may apply the usual scaling procedure in \rif{scalaraggio}, as
already done for instance to obtain \rif{stimaloc1}. Using such
estimates for the approximating problems \rif{Dirapp}, and employing
Lemma \ref{scalamorrey} we end up with the natural estimate
$$
\||Du_k|+s\|_{\MM^{m,\theta}(B_{R/2})} \leq
cR^{\frac{\theta-1}{p-1}-\frac{n}{q}}\||Du_k|+s\|_{L^{q}(B_{R})}
+c\|f_k\|_{L^{1,\theta}(B_R)}^{1/(p-1)} \;,
$$
for $q\in [p-1,b)$ and $c\equiv c(n,p,\ratio)$, whenever $B_R \CC
\Omega$. Using \rif{total2}-\rif{apbase}, and letting $k \nearrow
\infty$ we conclude with
$$
\||Du|+s\|_{\MM^{m,\theta}(B_{R/2})} \leq
cR^{\frac{\theta-1}{p-1}-\frac{n}{q}}\||Du|+s\|_{L^{q}(B_{R})}
+c\|\mu\|_{L^{1,\theta}(B_{R+\ep})}^{1/(p-1)},\quad \forall \
\ep>0\;.
$$
\end{remark}
\begin{remark}[On the limit case $\theta=p$] In proving Theorem \ref{main5} we used
\rif{zz2} from Lemma \ref{superrev} to estimate $ IV_k$. In turn
Lemma \ref{superrev} uses Gehring's lemma. The use of Gehring's
lemma is actually needed only in the borderline case $\theta=p$,
when $m=p$. Indeed in \rif{quasiq} we need $m-p\chi< 0$. Now observe
that $m < p$ as soon as $\theta < p$. Therefore in this latter case
we may use inequality \rif{zz2} with $p$ replacing $p \chi$; in
\rif{quasiq} we would have $H^{m-p}$, still small taking $H$ large.
All this does not need Gehring's lemma: indeed for solutions to
\rif{superrev} the basic Caccioppoli's type inequality, together
with Poincar\'e's one (see \cite{G}, Chapters 6-7) give
$$
\left(\mean{B_{R/2}} |Dv|^{p} \, dx\right)^{\frac{1}{p}}\leq c
\left(\mean{B_R} (|Dv|^{\frac{np}{n+p}} +s^{\frac{np}{n+p}})\, dx
\right)^{\frac{n+p}{np}}\;.
$$
From this the new form of \rif{zz2}, with $p$ replacing $p\chi$,
follows by Lemma \ref{revq}.
\end{remark}
\section{The super-capacitary case}
In this section we are going to prove Theorem \ref{main3}. As usual
we shall derive a priori estimates; in the following let $u \in
W^{1,p}_0(\Omega)$ be the solution to \trif{Dirapp2} for a fixed $f
\in L^{\infty}(\Omega)$. Take $B_{4R} \CC \Omega$ with $4R \leq 1$,
and then scale $u(x)$ in $B_R$ as in \rif{scalaraggio}, therefore
obtaining a solution $\tu (y)$ in $B_1$. We fix $d \in (0,1)$. Now
apply \rif{bettereh2} with $\Omega' \equiv B_{1/2}$ and $q=p-1$ to
have
\begin{eqnarray}
\nonumber \sup_{h} \int_{B_{1/2}}
\frac{|\tau_{i,h}D\tu(y)|^{p-1}}{|h|^{1-d}}\, dy &\leq & c
\||D\tu|+s\|_{L^{p-1}(B_1)}^{p-1}+ c\|\tilde{f}
\|_{L^{1}(B_1)} \\
&\leq & c \||D\tu|+s\|_{L^{p-1,\theta-1}(B_1)}^{p-1}+ c\|\tilde{f}
\|_{L^{1,\theta}(B_1)} \;,\label{dariscala} \end{eqnarray} for every
$i \in \{1,\ldots, n\}$, where $0< |h|< 1/4$ and $c\equiv
c(n,p,\ratio,d)$. We have used that $\sigma(p-1,\theta)=1$ and
$\delta(p-1)=\theta-1$ for every $\theta \in [p,n]$, by
\rif{ssiidef} and \rif{ildelta} respectively; recall also
\rif{tilf}. Notice that the application of \rif{bettereh2} to $\tu $
is legitimate since the arguments for Lemma \ref{222it} are local,
making no use of boundary information for the solution. Scaling back
\rif{dariscala} to $B_R$ via Lemma \ref{scalamorrey} gives
\eqn{decni1}
$$
\sup_{h} \int_{B_{R/2}} \frac{|\tau_{i,h}Du(x)|^{p-1}}{|h|^{1-d}}\,
dx\leq c M_{p-1}^{p-1}(B_R)R^{n-\theta+d}\;,
$$
where $c \equiv c(n,p,\ratio,d)$, for every $i \in \{1,\ldots, n\}$
where $0< |h|< R/4$, and where \eqn{emmeq}
$$
M_{q}(B_R):=\||Du|+s\|_{L^{q,\delta(q)}(B_R)}+
\|f\|_{L^{1,\theta}(B_R)}^{1/(p-1)},\qquad q \in [p-1,m)\;.
$$
Now take $q \in (p-1,m)$ and select $\gamma \in (0,1)$ such that
\eqn{interp00}
$$ q
=\frac{(\theta-\gamma)(p-1)}{\theta-1} \qquad \Longleftrightarrow
\qquad \sigma(q, \theta)=\gamma\;.
$$
If $\gamma_0 \in (0,\gamma)$ then \eqn{interp0}
$$
q < m_0 :=\frac{(\theta-\gamma_0)(p-1)}{\theta-1} < m\;,
$$
and write, with $t \in (0,1)$ \eqn{interp1}
$$
q=
(1-t)(p-1)+tm_0=\left(\frac{\gamma-\gamma_0}{1-\gamma_0}\right)(p-1)+\left(\frac{1-\gamma}{1-\gamma_0}\right)m_0\;.
$$
It follows \eqn{mm00}
$$
\frac{m-m_0}{m}=\frac{\gamma_0}{\theta}, \qquad \qquad
\frac{m_0}{m}=\frac{\theta-\gamma_0}{\theta}\;.
$$
Now, by estimate \rif{stimettapm1} we have
$$
\||Du|+s\|_{\MM^m(B_R)} \leq c \left[
\||Du|+s\|_{L^{p-1,\theta-1}(B_{2R})}+\|f\|_{L^{1,\theta}(B_{2R})}^{\frac{1}{p-1}}
\right]R^{\frac{n-\theta}{m}}\;.
$$
Using Lemma \ref{hollo}, the latter estimate and \rif{mm00} we find
\begin{eqnarray}
\nonumber && \int_{B_R} |Du|^{m_0}\, dx \leq m(m-m_0)^{-1}
R^{n-\frac{m_0n}{m}}\|Du\|_{\MM^m(B_R)}^{m_0}\\&& \nonumber \qquad
\leq c \gamma_0^{-1}R^{n-\frac{m_0n}{m}+\frac{m_0}{m}(n-\theta)}
\left[
\||Du|+s\|_{L^{p-1,\theta-1}(B_{2R})}^{m_0}+\|f\|_{L^{1,\theta}(B_{2R})}^{\frac{m_0}{p-1}}
\right]\\ && \qquad \leq c\gamma_0^{-1} R^{n-\frac{m_0\theta}{m}}
M_{p-1}^{m_0}(B_{2R})= c R^{n-\theta+\gamma_0}
M_{p-1}^{m_0}(B_{2R})\;,\label{interp11}
\end{eqnarray}
where $c \equiv c(n,p,\ratio,\gamma_0)$. Now, by \rif{interp1} and
H\"older's inequality
\begin{eqnarray}
\nonumber && \int_{B_{R/2}} |\tau_{i,h}Du(x)|^{q}\, dx =
\int_{B_{R/2}}
|\tau_{i,h}Du(x)|^{(1-t)(p-1)}|\tau_{i,h}Du(x)|^{tm_0}\, dx\\
\nonumber & &\qquad \leq \left(\int_{B_{R/2}}
|\tau_{i,h}Du(x)|^{p-1}\, dx \right)^{1-t}\left(\int_{B_{R/2}}
|\tau_{i,h}Du(x)|^{m_0}\, dx \right)^{t}\\ \nonumber && \qquad
\qquad \leq |h|^{(1-d)(1-t)}\left( \int_{B_{R/2}}
\frac{|\tau_{i,h}Du(x)|^{p-1}}{|h|^{1-d}}\, dx\right)^{1-t}\cdot
\\ && \hspace{4cm}\cdot \left(\int_{B_{R/2}}
(|Du(x)|^{m_0}+|Du(x+he_i)|^{m_0})\, dx\right)^t\;. \label{interp2}
\end{eqnarray}
In turn, taking $|h|\leq R/4$, and using \rif{interp11} we have
\begin{eqnarray}
\nonumber \int_{B_{R/2}} (|Du(x)|^{m_0}+|Du(x+he_i)|^{m_0})\, dx
&\leq &
2\int_{B_{R}} |Du|^{m_0}\, dx \\
&\leq & c M^{m_0}_{p-1}(B_{2R})R^{n-\theta+\gamma_0}\;.
\label{interp3}
\end{eqnarray}
Combining \rif{decni1}, \rif{interp2} and \rif{interp3}, and taking
into account \rif{interp1} we have
$$
 \sup_{h}\int_{B_{R/2}}
\frac{|\tau_{i,h}Du(x)|^{q}}{|h|^{(1-d)(1-t)}}\, dx \leq c
M^{q}_{p-1}(B_{2R})R^{n-\theta+(1-t)d+t\gamma_0}\;,
$$
where $c \equiv (n,p,\ratio,d,\gamma_0)$, and $h$ is a real number
such that $0< |h|< R/4$. Since $i \in \{1,\ldots, n\}$ is arbitrary
the last inequality yields in a standard way \eqn{interp5}
$$
 \sup_{h}\int_{B_{R/2}}
\frac{|Du(x+h)-Du(x)|^q}{|h|^{\sigma}}\, dx \leq c
M^{q}_{p-1}(B_{2R})R^{n-\theta+(1-t)d+t\gamma_0}\;,
$$
where this time $h\in \er^n$ with $|h| \in (0,R/4]$. Here we have
set \eqn{sigma000}
$$
\sigma:=
(1-d)(1-t)\stackrel{\rif{interp00},\rif{interp1}}{=}(1-d)\left(\frac{\sigma(q,\theta)-\gamma_0}{1-\gamma_0}\right)\;.
$$
As $\sigma(q,\theta)=\gamma$ by \rif{interp00}, a direct computation
reveals that $ (1-t)d+t\gamma_0> d\sigma(q,\theta);$ using that
$R\leq 1$ and the latter inequality in \rif{interp5} we have
\eqn{interp6}
$$
 \sup_{h}\int_{B_{R/2}}
\frac{|Du(x+h)-Du(x)|^q}{|h|^{\sigma}}\, dx \leq c
M^{q}_{p-1}(B_{2R})R^{n-\theta+d \sigma(q,\theta)}\;,
$$
with $\sigma$ as in \rif{sigma000} and $c \equiv
c(n,p,\ratio,d,\gamma_0)$. Estimate \rif{interp6} has been proved
for $q \in (p-1,m)$. It actually holds for the case $q=p-1$ too, and
even with $\gamma_0=0$ in \rif{sigma000}. This is just a consequence
of $\sigma(p-1,\theta)=1$ and \rif{decni1}. We are now ready to
conclude the proof. Take $\gamma_1 \in (0,\sigma)$; then changing
variables
\begin{eqnarray}
\nonumber && \int_{B_{R/2}} \int_{B_{R/2}} \! \frac{|Du(x) - Du(y)
|^{q}}{|x-y|^{n+\sigma-\gamma_1}} \, dx \, dy \\&& \leq
\int_{B(0,R)} \frac{1}{|h|^{n-\gamma_1}} \int_{B_{R/2}}
\frac{|Du(x+h)-Du(x)|^q}{|h|^{\sigma}}\, dx\, dh \nonumber\\
& & \leq \left(\int_{B(0,R)}
\frac{dh}{|h|^{n-\gamma_1}}\right)\sup_{h}\int_{B_{R/2}}
\frac{|Du(x+h)-Du(x)|^q}{|h|^{\sigma}}\, dx \nonumber \\
&& \stackrel{\rif{interp6}}{\leq} c(n)\gamma_1^{-1}
M^{q}_{p-1}(B_{2R})R^{n-\theta+d \sigma(q,\theta)}\;,\label{ffrac2}
\end{eqnarray}
valid for any $q \in [p-1,m)$, where $c \equiv
c(n,p,\ratio,\gamma_0)$. Therefore since in \rif{sigma000} and
\rif{ffrac2} $\gamma_0,\gamma_1$ can be picked arbitrarily small,
all in all we have proved that \eqn{interp7}
$$
[Du]_{\sigma/q, q; B_R} \leq c M^{q}_{p-1}(B_{4R})R^{n-\theta+d
\sigma(q,\theta)},\qquad \sigma < (1-d)\sigma(q,\theta)\;,
$$
for all balls $B_{R}$ such that $B_{4R}\CC \Omega$, and $q \in
[p-1,m)$. The constant $c$ depends on $n,p,\ratio,q,d,\sigma$. This
needs an explanation. The constant $c$ blows up when $q \nearrow m$
and/or $\sigma\nearrow (1-d)\sigma(q,\theta)$. Indeed, taking $q$
close to $m$ forces $\gamma_0$ to be small via \rif{interp0}, and
this in turn increases $c$ via \rif{interp11}; on the other hand
taking $\sigma$ close to $ (1-d)\sigma(q,\theta)$ forces
$\gamma_0,\gamma_1$ to be small via \rif{sigma000},\rif{ffrac2}
respectively, and this again increases $c$ via \rif{interp11} and
\rif{ffrac2}.

Now using \rif{interp7} together with a standard covering argument,
and the fact that $d$ can be chosen arbitrarily small, and taking
into account the definition of $M_{q}(B_R)$ in \rif{emmeq}, we
conclude that for every couple of open subsets $\Omega' \CC \Omega''
\CC \Omega$ and every $\sigma < \sigma(q,\theta)$ it holds
\eqn{interp8}
$$
\|Du\|_{W^{\sigma/q, q, \theta}(\Omega')} \leq
c\||Du|+s\|_{L^{p-1,\theta-1}(\Omega'')}+
c\|f\|_{L^{1,\theta}(\Omega'')}^{1/(p-1)}\;.
$$
Finally using \rif{stimamorreypre} with $q=p-1$, and changing
subsets via Lemma \ref{stanom} we gain \eqn{interp9}
$$
\|Du\|_{W^{\sigma/q, q, \theta}(\Omega')} \leq c
\|f\|_{L^{1}(\Omega)}^{1/(p-1)}+c\|f\|_{L^{1,\theta}(\Omega)}^{1/(p-1)}+
cs|\Omega|^{1/q}\;,
$$
and the constant depends on $n,p,\ratio,q,\sigma,
\dist(\Omega',\partial\Omega)$. The assertion of Theorem
\ref{main3}, together with estimate \rif{apestmo} follow via the
approximation scheme of Section 5 as for the other proofs of this
paper. \qed
\begin{remark}[Fractional differentiability vs Morrey regularity]
Let us go back to \rif{interp7}, keep now $d$ fixed in $(0,1)$, not
necessarily ``small" in order to approach $\sigma(q,\theta)$ with
$\sigma$. Then, again via the approximation of Section 5, it easily
follows \eqn{interp10}
$$
Du \in W^{\sigma/q,q,\theta+d\sigma(q,\theta)}_{\loc}(\Omega,
\er^n), \qquad \  \mbox{for every} \ \sigma <
(1-d)\sigma(q,\theta)\,.
$$
With $q \in [p-1,m)$ being fixed, inclusion \rif{interp10} tells us
that if we decrease the rate of differentiability down to
$(1-d)\sigma(q,\theta)$, we gain in the Morrey scale up to $\theta
+d \sigma(q,\theta)$. Observe that inclusion \rif{interp10} realizes
a perfect interpolation between the maximal differentiability in
\rif{fra13mo} that we may obtain taking $d$ close to $0$, and the
maximal Morrey regularity in \rif{stimamorrey} that we may obtain
formally letting $d \nearrow 1$ in \rif{interp10}, as
$\theta+d\sigma(q,\theta)\nearrow \delta(q)$ when $d \nearrow 1$;
look at \rif{ildelta} and \rif{ssiidef}. In other words, with a very
rough but suggestive notation
$$
\lim_{d\searrow 0}
W^{(1-d)\sigma(q,\theta)/q,q,\theta+d\sigma(q,\theta)} =
W^{\sigma(q,\theta)/q,q,\theta},$$and $$  \lim_{d\nearrow 1}
W^{(1-d)\sigma(q,\theta)/q,q,\theta+d\sigma(q,\theta)} =
L^{q,\delta(q)}\;.
$$
More in general, since when considering Morrey decay properties as
\rif{morrey} the exponent $\theta$ replaces $n$ everywhere, the
integer dimension of the space $W^{\alpha, q,\theta}$ should be
defined as $\alpha - \theta/q$, compare with Remark \ref{integerex}.
In this respect, exactly as in Remark \ref{integerex}, all the
spaces $W^{(1-d)\sigma(q,\theta)/q,q,\theta+d\sigma(q,\theta)}$
share the same integer dimension $(\theta-1)/(p-1)$, for every
possible choice of $q \in [p-1,m)$ and $d \in (0,1)$.
\end{remark}
\section{Sharpness, comparisons, extensions}
We hereby discuss the sharpness of some of the foregoing results,
and outline a few extensions and connections.

\subsection{\bf Sharpness of Theorem \ref{main1}} The result in \rif{fra1} is
sharp for every choice of the couple $(q, \sigma(q))$ in the range
\rif{ssii}, and in particular the inclusions \rif{pipi} and
\rif{duedue} are sharp too. Indeed, we cannot have $Du \in
W^{\sigma(q)/q,q}_{\loc}$, as shown by the usual counterexample
\cite{KiXu}. Consider the equation \rif{pmeas} in the ball $B_1
\equiv \Omega$, with $\mu \equiv \delta$, the Dirac measure charging
the origin, with the related zero-Dirichlet condition. The unique
solution to problem \rif{Dir1} is now given by the Green's function
$$u(x):=c(n,p)\left\{\begin{array}{ccc}
|x|^{\frac{p-n}{p-1}}-1 & \quad \mbox{if} & \quad p<n
\\ \log|x| & \quad \mbox{if} & \quad p=n\end{array}\right.\;,$$
where $c(n,p)$ is a suitable re-normalization constant. We have $Du
\in \MM^{b}(B_1)$, but $Du \not\in L^{b}_{\loc}(B_1)$, and crucial
integrability is lost at the origin. Now, assume by contradiction
that $Du \in W^{\sigma(q)/q,q}_{\loc}(B_1)$, then by Theorem
\ref{fraemb} we would have $Du \in
L^{nq/(n-\sigma(q))}_{\loc}(B_1)$, but this is impossible since
$nq/(n-\sigma(q))=b$ by \rif{ssii}. Therefore $Du \not\in
W^{\sigma(q)/q,q}_{\loc}(B_1)$, and this gives the optimality of
Theorem \ref{main1}. On the other hand, as $nq/(n-\sigma(q))=b$,
then assuming Theorem \ref{main1} allows to recover the original
integrability result in \rif{boccardo} in a local form, again via
Theorem \ref{fraemb}.

\subsection{\bf About Theorem \ref{noncz}.}
This is also sharp. In fact assuming that
$$V(Du) \in
W^{\frac{p}{2(p-1)},\frac{2(p-1)}{p}}_{\loc}(\Omega,\er^n), \qquad
\mbox{for every}\ \ep \in (0,1)\;$$ by the fractional Sobolev
embedding theorem \ref{fraemb} we would get
$$
V(Du) \in L^{\frac{n(p-1)2}{(n-1)p}}_{\loc}(\Omega,\er^n)\,,
$$
and in turn this would imply, via \rif{elemV}, that
$$
Du \in L^{\frac{n(p-1)}{n-1}}_{\loc}(\Omega,\er^n)\,,
$$
which is excluded by the discussion of Section 11.1. Theorem
\ref{noncz} can be regarded as a non-linear version of the so called
``uniformization of singularities principle", well-known in Complex
Analysis: raising a function to a suitably large power we get a
function with better regularity properties. In such respect we
conclude with an open problem that for the sake of simplicity we
state for solutions to equations involving the $p$-Laplacean
operator \rif{pmeas}. Take $\gamma \in \er$ such that \eqn{ranga}
$$ \frac{p-2}{2}\leq \gamma \leq p-2
$$
and prove - or disprove - in the spirit of Theorem \ref{noncz},
that, once a Dirichlet class is fixed as boundary datum, there
exists a SOLA solution to \rif{pmeas} such that \eqn{unising}
$$
|Du|^{\gamma}Du \in
W^{\frac{\gamma+1}{p-1}-\ep,\frac{p-1}{\gamma+1}}_{\loc}(\Omega,\er^n)\,,
\qquad \mbox{for every} \ \ \ep >0\;.$$ In the first limit case
$\gamma = (p-2)/2$ this is essentially the content of Theorem
\ref{noncz}, while in the other borderline case $\gamma = p-2$ this
amounts to prove that $|Du|^{p-2}Du \in
W^{1-\ep,1}_{\loc}(\Omega,\er^n)$, for every $\ep>0$. When $p=2$ all
such statements collapse in Theorem \ref{main1}. Observe that,
exactly as for \rif{propro}, for every choice of $\gamma$ in the
range \rif{ranga} the product between the differentiability and the
integrability indexes in \rif{unising} remains constant, up to the
presence of $\ep$.
\subsection{The exponent $m$ in \rif{mingione}} We now
demonstrate the optimality of $m$ in \rif{lolo} in the case $p=2$ by
comparing Theorem \ref{main5} with the optimal ones of Adams
\cite{adamsduke} for the case $\triangle u = \mu$. Since our results
are local, up to a standard localization procedure we shall consider
the latter equation in the whole $\er^n$. We consider the fractional
integral operator defined by
$$I_{\alpha}(\mu)(x):=\int_{\er^n}\frac{d \mu(y)}{|x-y|^{n-\alpha}}, \qquad \qquad \alpha \in (0,n]\;.$$
When $\mu$ has compact support, the unique solution to $\triangle u
= \mu$ is given by $u(x):=c_1I_2(\mu)(x)$, with $c$ being a suitable
re-normalization constant; as a consequence $Du(x)=c_2I_1(\mu)(x)$;
see also \cite{LSW}. Now we recall the following result of Adams
\cite{adamsduke}: \eqn{adamsresult}
$$I_{\alpha}:L^{1,\theta} \to \MM^{\theta/(\theta-\alpha), \theta}\;,$$ that is sharp in the sense that we
cannot expect $I_{\alpha}(\mu) \in L^{\theta/(\theta-\alpha)}$, even
locally, for $\mu \in L^{1,\theta}$, see \cite{adamsduke} page 770,
no.~2. Taking in our case $\alpha= 1$ gives
$\theta/(\theta-\alpha)=m$, and therefore the exponent $m$ is the
natural one for $p=2$.

The case $p>2$ cannot be treated by such an argument since no
explicit representation formula is available for solutions to
\rif{pmeas}. We just remark that in the case $p>2$ the exponent $m$
is obtained by multiplying the one for $p=2$ times $(p-1)$. This
appears to be a natural phenomenon for measure data problems
\cite{BG2}. We hereby conjecture that the exponent $m$ is optimal
for every $p>2$. Finally observe how the fact that $\theta$ replaces
$n$ everywhere when assuming \rif{morrey} is in perfect accordance
with the embedding properties for Sobolev-Morrey spaces. Indeed,
assuming $Du \in L^{p,\theta}$ with $p< \theta$, leads to the
improved embedding $u \in L^{\theta p/(\theta-p)}$ \cite{adamsduke,
Ca63, Ca64}; this covers the usual Sobolev embedding theorem when
$\theta=n$.

\subsection{\bf Sharpness of Theorem \ref{main3}} Here we discuss
the optimality of the choice of the couple $(q, \sigma(q,\theta))$
in \rif{fra13mo} in the range displayed in \rif{ssii2}. The input
here in the Sobolev-Morrey embedding Theorem in the fractional case.
We have that $W^{\alpha, q, \theta}$ embeds in $L^t$ for every $t<
\theta q/(\theta-\alpha q)$ whenever $\alpha q < \theta$; see for
instance \cite{Ross}. Now take $p=2$ and assume that $Du \in
W^{(\sigma(q,\theta)+\ep)/q, q}_{\loc}$ for some $\ep>0$; since $m
=\theta q /(\theta - \sigma(q,\theta))$  we would conclude with $Du
\in L^{m}_{\loc}$, which is impossible at least when $p=2$, as seen
a few lines above. Similarly, as the optimality of $m$ in \rif{lolo}
is expected when $p>2$, the optimality of \rif{fra13mo}-\rif{ssii2}
in the case $p>2$ is expected too. In fact, this is the same
argument used to get the optimality of Theorem \ref{main1} at the
beginning of the section.

\subsection{\bf Lebesgue vs Morrey} Assuming \rif{morrey} improves on
\rif{Marbase} up to \rif{lolo}. Now assume that $\mu \in L^{t}$ for
$t \in [1,(p^*)')$; in this case $Du \in L^{g}$ with $g
=nt(p-1)/(n-t)$ \cite{BG2, KiL}; in particular $Du \in \MM^{g}$. On
the other hand $\mu \in L^t$ implies that $\mu$ satisfies
\rif{morrey} with $\theta=n/t$; in this case Theorem \ref{main5}
gives $Du \in \MM^{m}$ with $m = n(p-1)/(n-t)$, that is worse than
$Du \in \MM^{g}$, but for $t=1$. This does not contradict the
sharpness of \rif{lolo}. Indeed we may find functions $f \in
L^{1,\theta}$, with $\theta$ arbitrarily close to zero, such that $f
\not\in L^{t}$ for any $t>1$, see \cite{G}, comments at Chapter 2.
On the other hand, truncation techniques fully apply in the case of
$L^t$ data \cite{KiXu}, because {\em functions can be truncated,
while measures cannot}, and better integrability of $Du$ follows.

\subsection{\bf Systems} Theorem \ref{main4} extends to
systems, under assumptions \rif{asp} and \rif{lipi}, when obviously
recast for the vectorial case; $u : \Omega \to \er^N$, $z \in \Ma$
and so on. In this case the measure $\mu$ takes its values in
$\er^N$. Indeed for Theorem \ref{main4} we do not need Lemma
\ref{coco1}; this employs the truncation operators \rif{troncamenti}
and they do not work for general elliptic systems. We also do not
need Lemma \ref{superrev}, which under the general assumptions
\rif{asp} only works in the scalar case. The only basic ingredients
are Lemmas \ref{coco5}-\ref{coco52} and \ref{revbase}. The first two
only need monotonicity in \rif{asp}$_1$, while the third one is here
stated directly in the vectorial case $N\geq 1$. Anyway, we are
planning further extensions to certain special classes of systems.

\subsection{\bf Condition \rif{morrey}} This can be relaxed in a local one,
since the results we are giving are local. More precisely, we may
assume that for every $\Omega'' \CC \Omega$ there exists a constant
$M(\Omega'')$ such that \eqn{morreyrel}
$$
|\mu|(B_R) \leq M(\Omega'') R^{n-\theta}, \qquad \  \mbox{for every
ball} \ B_R \CC \Omega''\;.
$$
Roughly, we are considering $\mu \in L^{1,\theta}_{\loc}(\Omega)$
rather than $\mu \in L^{1,\theta}(\Omega)$. When assuming
\rif{morreyrel} instead of \rif{morrey} the inclusions of Theorems
\ref{main5}-\ref{main6} still hold, but the a priori estimates
change. We give the new statement for the estimate of Theorem
\ref{main5}, the others to be modified in a similar fashion. For
every couple of open subsets $\Omega' \CC \Omega'' \CC \Omega$ there
exists a constant $c$ depending on $n,p,\ratio,\Omega',\Omega''$
$$
\|Du\|_{\MM^{m, \theta}(\Omega')} \leq c
\|Du\|_{L^{p-1}(\Omega'')}+c[M(\Omega'')]^{1/(p-1)}
+cs|\Omega''|^{1/m}\;.
$$
Moreover there exists $c$ depending on $n,p,\ratio,\Omega', \Omega$
such that
$$
\|Du\|_{\MM^{m, \theta}(\Omega')} \leq c
[|\mu|(\Omega'')]^{1/(p-1)}+c[M(\Omega'')]^{1/(p-1)}
+cs|\Omega''|^{1/m}\;.
$$

\end{document}